\documentclass[a4paper,12pt]{article}         
\def\bysame{{\bf --- }}                      

\def\~{{\bf --}}

\usepackage{newlfont}
\usepackage{amsmath,amsfonts,amssymb,amscd,epsfig}
\usepackage{latexsym}
\usepackage{amsthm}

\usepackage[dvipdfm,ps2pdf=true,dvips=true, 
hyperindex=false,pagebackref=true,bookmarks=true, 
pdfpagemode=None,colorlinks=true,linkcolor=blue,citecolor=red]{hyperref} 
\def\~{{\rm --}}  

\newcommand{\comment}[1]{}

\newcommand{\BR}{{\mathbb R}}
\newcommand{\BQ}{{\mathbb Q}}
\newcommand{\BC}{{\mathbb C}}
\newcommand{\BP}{{\mathbb P}}
\newcommand{\BZ}{{\mathbb Z}}
\newcommand{\BN}{{\mathbb N}}
\newcommand{\BS}{{\mathbb S}}

\newcommand{\cH}{{\mathcal H}}

\newcommand{\cB}{{\mathcal B}}

\newcommand{\cP}{{\mathcal P}}
\newcommand{\cX}{{\mathcal X}}
\newcommand{\cY}{{\mathcal Y}}

\newcommand{\Z}{{\mathbb Z}}
\newcommand{\Q}{{\mathbb Q}}
\newcommand{\N}{{\mathbb N}}
\newcommand{\C}{{\mathbb C}}
\newcommand{\R}{{\mathbb R}}

\def\HH{\mbox{${\mathcal H}$\kern-5.2pt${\mathcal H}$}}

\def\der{\partial}
\def\tensor{\otimes}

\def\kap{\kappa}
 
\def\Comp{{\mathbb C}}
\def\sM{{\mathcal M}}

\newtheorem{theorem}{Theorem}[section]
\newtheorem{proposition}[theorem]{Proposition}
\newtheorem{definition}[theorem]{Definition}

\newtheorem{corollary}[theorem]{Corollary}

\newtheorem{theorem }{Theorem}[section]
\newtheorem{proposition }[theorem]{Proposition}
\newtheorem{definition }[theorem]{Definition}
\newtheorem{lemma }[theorem]{Lemma}
\newtheorem{corollary }[theorem]{Corollary}
\newtheorem{notation }[theorem]{Notation}
\newtheorem{remark }[theorem]{Remark}
\newtheorem{example }[theorem]{Example}

\newtheorem{ theorem}{Theorem}[section]
\newtheorem{ proposition}[theorem]{Proposition}
\newtheorem{ definition}[theorem]{Definition}
\newtheorem{ lemma}[theorem]{Lemma}
\newtheorem{ corollary}[theorem]{Corollary}
\newtheorem{ notation}[theorem]{Notation}
\newtheorem{ remark}[theorem]{Remark}
\newtheorem{ example}[theorem]{Example}

 \newcommand{\rem}{{\bf Comment.\ }}
 \newcommand{\rmk}{{\bf Comment.\ }}

\newtheorem*{acka}{Acknowledgments}

\def\for{\  \hbox{ for } \ }
\def\iif{ \ \hbox{ if } \ }

\def\where{\  \hbox{ where } \ }
\def\and{\  \hbox{ and } \ }
\def\and{\  \hbox{ and } \ }

\def\equal{\stackrel{\,\mathbf{def}}{= \kern-3pt =}}

\def\la{\lambda}
\def\La{\Lambda}
\def\om{\omega}

\def\th{\theta}
\def\al{\alpha}

\def\ga{\gamma}

\def\de{\delta}
\def\De{\Delta}
\def\ka{\kappa}
\def\kapp{\hbox{\bf \ae}}
\def\si{\sigma}

\def\Ga{\Gamma}
\def\ze{\zeta}

\def\pa{\partial}

\def\vep{\varepsilon}

\def\vth{{\vartheta}}

\def\tal{\tilde{\alpha}}

\def\tV{\tilde{V}}

\def\tGa{\tilde{\Gamma}}

\def\tw{\widetilde w}
\def\tW{\widetilde W}

\def\tV{\tilde V}
\def\tz{\tilde z}
\def\tb{\tilde b}

\def\tR{\tilde R}

\def\hw{\widehat{w}}
\def\hW{\widehat{W}}
\def\hu{\hat{u}}

\def\hv{\hat{v}}

\def\I{\mathbf{I}}

\def\0{\mathbf{0}}
\def\H{\mathbf{H}}
\def\V{\mathbf{V}}

\def\f{\mathcal{F}}
\def\çF{\mathcal{F}}

\def\m{\mathcal{M}}

\def\d{\mathcal{D}}
\def\p{\mathcal{P}}
\def\cP{\mathcal{P}}
\def\a{\mathcal{A}}
\def\h{\mathcal{H}}

\def\y{\mathcal{Y}}
\def\e{\mathcal{E}}

\def\x{\mathcal{X}}

\def\g{\mathcal{G}}

\def\w{\mathcal{W}}
\def\i{\mathcal{I}}
\def\j{\mathcal{J}}
\def\b{\mathcal{B}}

\def\lan{\langle}

\def\ran{\rangle}

 \def\dim{{\hbox{\rm dim}}_{\mathbb C}\,}
\def\lng{\hbox{\tiny lng}}
\def\sht{\hbox{\tiny sht}}

\def\gl{\mathfrak{gl}_N}

\newcommand{\Aut}{\operatorname{Aut}}

\newcommand{\Ind}{\operatorname{Ind}}

\newcommand{\sq}{\phantom{1}\hfill$\qed$}

\newcommand{\lr}{\langle}
\newcommand{\rr}{\rangle}
\newcommand{\eps}{\varepsilon}

\newcommand{\Res}{\mbox{Res}\;}

\def\HH{\mathfrak{H}}

\def\CC{\mathfrak{C}}

\def\HH{\hbox{${\mathcal H}$\kern-5.2pt${\mathcal H}$}}

\font\smm=msbm10 at 12pt 
\def\symbol#1{\hbox{\smm #1}}
\def\lsmash{{\symbol n}}
\def\rsmash{{\symbol o}}
\def\#{\sharp}

\font\eightbf=cmbx8

\begin{document}

\comment{ This paper is based on the introduction 
to the monograph ``Double affine Hecke algebras'' 
to be published by Cambridge University Press. 
The connections with Knizhnik-Zamolodchikov 
equations, Kac-Moody algebras, tau-function, 
harmonic analysis on symmetric spaces, and special functions
are discussed. The rank one case is considered in
detail including the classification of Verlinde algebras
and their deformations, Gauss-Selberg integrals and 
Gaussian sums, a topological interpretation of DAHA, 
a relation of the rational DAHA to sl(2), and applications 
to the diagonal coinvariants. The last three sections are 
devoted to relations of the general DAHAs to the p-adic 
affine Hecke algebras, trigonometric and rational DAHAs, 
and applications to the Harish-Chandra theory. The purpose  
of this introduction is a demonstration that DAHA can be 
considered as a natural formalization of the concept of the  
Fourier transform in mathematics and physics.
}

\def\blfootnote{\xdef\@thefnmark{}\@footnotetext} 

\title{\bf Introduction to \\ Double Hecke algebras}
\author{Ivan Cherednik\,${}^{1,2}$}

\date{\today}

\newcommand{\agree}
{\footnotetext[1]
{Please read this legal notice carefully before accessing 
the content on this work.\\
BY ACCESSING OR USING THIS WORK, YOU AGREE TO BE BOUND BY 
THE TERMS AND CONDITIONS SET FORTH BELOW. 
IF YOU DO NOT WISH TO BE BOUND BY THESE TERMS AND CONDITIONS, 
YOU MAY NOT ACCESS OR USE THE WORK. IVAN CHEREDNIK
MAY MODIFY THIS AGREEMENT AT ANY TIME, AND SUCH 
MODIFICATIONS SHALL BE EFFECTIVE IMMEDIATELY UPON NOTIFICATION.
YOUR CONTINUED USE OF THE WORK SHALL BE DEEMED YOUR CONCLUSIVE
ACCEPTANCE OF THE MODIFIED AGREEMENT.\medskip\\ 
COPYRIGHT INFORMATION\\
\copyright Ivan Cherednik. The entire contents of this work are 
protected by international copyright and trademark laws. 
All rights reserved.\medskip\\
YOU MAY NOT MODIFY, REPRODUCE, 
REPUBLISH, UPLOAD, POST, TRANSMIT, OR DISTRIBUTE, IN ANY MANNER, 
THE MATERIAL OF THIS WORK, INCLUDING TEXT AND/OR GRAPHICS 
WITHOUT THE EXPRESS PERMISSION OF IVAN CHEREDNIK.\\ 
You may print and copy portions of material solely for your own 
non-commercial use provided that you agree not to change or 
delete any copyright or proprietary notices from the materials. 
Violators of this policy may be subject to legal action.}} 

\footnotetext[1]{Partially supported by NSF grant 
DMS-0200276}
\footnotetext[2]{Department of Mathematics, UNC 
Chapel Hill, NC 27599-3250, USA
\newline$^{}$\hskip 20pt chered@email.unc.edu}
\maketitle  

\def\and{\  \hbox{ and } \ }



This paper is based on the introduction
to the monograph "Double affine Hecke algebras"
to be published by Cambridge University Press.
It is based on a series of lectures 
delivered by the author in Kyoto  (1996--1997), 
at University Paris 7 (1997--1998), at Harvard University
in 2001, and in several other places, including recent
talks at the conferences "Quantum Theory and Symmetries 3" 
(Cincinnati, 2003), "Geometric methods in algebra and 
number theory" (Miami, 2003),
and also at RIMS (Kyoto University), 
MIT, and UC at San Diego in 2004.

The connections with Knizhnik\~Zamolodchikov
equations, Kac\~Moody algebras, harmonic analysis on
symmetric spaces, special functions
are discussed. We demonstrate that
the $\tau$\~function from soliton theory is 
a generic solution of the so-called $r$\~matrix KZ
with respect to the Sugawara $L_{-1}$\~operators,
which is an important part of the theory of 
integral formulas of the KZ equations. 

The double affine Hecke algebra (DAHA) of type $A_1$
is considered in detail including the classification of 
the nonsymmetric Verlinde algebras, their deformations, 
Gauss\~Selberg integrals and Gaussian sums, the topological
interpretation, the relation of the rational DAHA 
to $\mathfrak{sl}(2),$ and recent applications
to the diagonal coinvariants. The last three sections of the
paper are devoted to the general DAHA, its origins
in the classical $p$\~adic theory of affine Hecke
algebras, the trigonometric and rational
degenerations, and applications to the Harish-Chandra theory.
\medskip

Affine Hecke algebras appeared as a technical tool in the theory
of automorphic functions, but now are indispensable in
modern representation theory, combinatorics, geometry,
harmonic analysis, mathematical physics, and the theory
of special functions. The Langlands program, the
theory of quantum groups, the modern theory of the symmetric group,
and the hypergeometric functions are well known examples.

A major development in this field
was the introduction of DAHA
in \cite{C15} (the reduced root systems).
A counterpart of DAHA for $C^\vee C$  was defined by 
Noumi and Sahi and, then, was used in \cite{Sa,St} 
to prove the Macdonald conjectures  
for the Koornwinder polynomials and extend
the Mehta\~Macdonald formula to the  $C^\vee C$\~case,
completing the theory from \cite{C2,C3,C5}. 

DAHAs now are a powerful tool in representation theory,
with impressive applications. We demonstrate
in this paper that these algebras can be considered
as a natural formalization of the concept of the 
Fourier transform in mathematics and physics.
\smallskip

\begin{acka}
I would like to thank M.~Duflo, P.~Etingof, A.~Garsia, 
M.~Kashiwara, D.~Kazhdan, T.~Miwa,
M.~Nazarov, E.~Opdam, V.~Ostrik, 
N.~Wallach for discussions and help with
improving the paper.  Some its parts 
are connected with \cite{C14}, published
by MSJ, and my note in the Proceedings of QTS3 to be 
published soon.
\end{acka}
\smallskip

{\normalsize
\tableofcontents
}


\vfill\eject

\setcounter{section}{-1} 
\section{Universality of Hecke algebras} \label{SEC:Heckea}
\subsection{Real and imaginary} 

\index{Q@$\Q$:\ rational numbers|(} 
\index{R@$\R$:\ real numbers|(} 
\index{R-a@$R,r$:\ $r$\~matrices|(} 
\index{C@$\C$:\ complex numbers|(} 
\index{N@$\N$:\ natural numbers|(} 
\index{Z@$\Z$:\ integers|(} 
\index{S@$\BS_n$:\ permutation group|(} 
\index{F@$\Phi$:\ solution of KZ,QKZ|(} 
\index{f@$\varphi$:\ solution of QMBP|(} 
\index{sa@$s_i$:\ simple reflections|(} 
\index{a@$\alpha_i$:\ simple roots|(} 
\index{pa@$\pi_r$:\ elements from $P^\vee/Q^\vee$|(} 
\index{Pa@$P$:\ weight lattice|(} 
\index{Paa@$\p$:\ polynomial rep|(} 
\index{Qa@$Q$:\ root lattice|(} 
\index{Sa@$\Sigma$:\ root system|(} 

Before a systematic exposition, I will try to outline the 
connections of the representation theory of Lie groups, 
Lie algebras, and Kac\~Moody algebras with the Hecke algebras and 
the Macdonald theory. 

\begin{figure}[htbp] 
\begin{center} 
\resizebox{\textwidth}{!}{\input{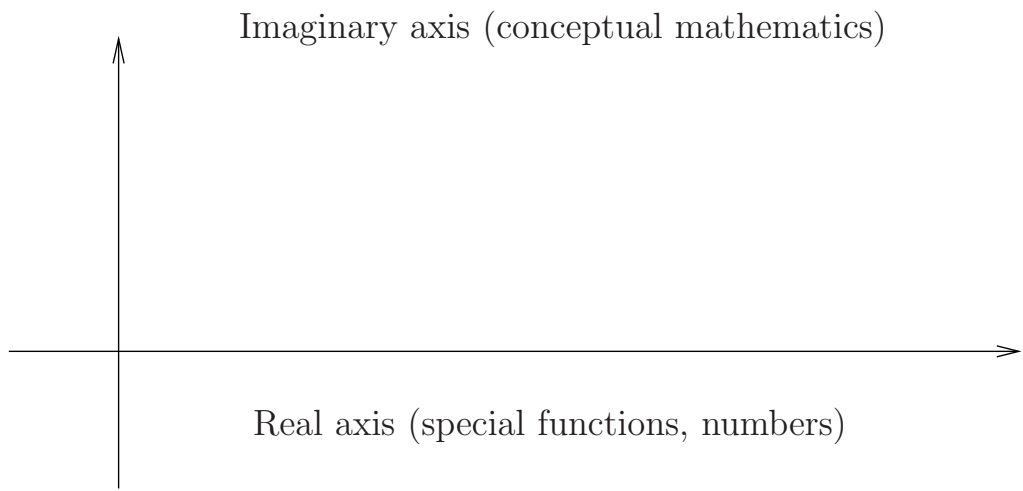}} 
\caption{Real and Imaginary} 
\label{figure1} 
\end{center} 
\end{figure}

The development of mathematics may 
be illustrated by Figure \ref{figure1}.

Mathematics is fast in the imaginary (conceptual) direction 
but, generally, slow in the real direction (fundamental objects). 
Mainly I mean modern mathematics, 
but it may be more universal. 
For instance, the ancient Greeks created a highly conceptual 
axiomatic geometry with modest ``real output." 
I do not think that the ratio $\Re/\Im$ is much higher 
now. 

\begin{figure}[htbp] 
\begin{center} 
\resizebox{\textwidth}{!}{\begin{picture}(0,0)%
\special{psfile=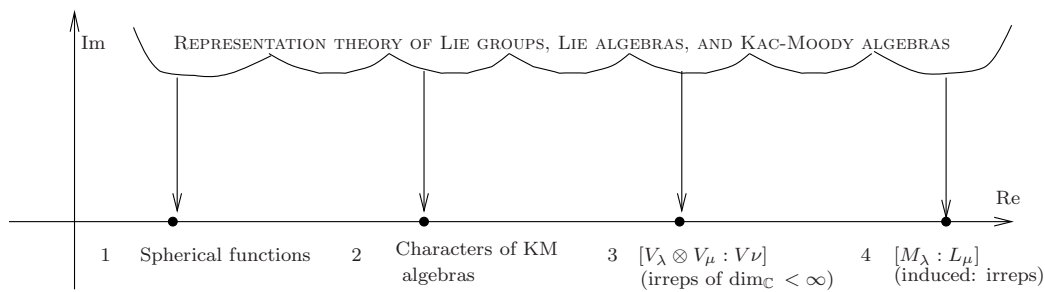}%
\end{picture}%
\setlength{\unitlength}{0.01250000in}%
\begingroup\makeatletter\ifx\SetFigFont\undefined
\def\x#1#2#3#4#5#6#7\relax{\def\x{#1#2#3#4#5#6}}%
\expandafter\x\fmtname xxxxxx\relax \def\y{splain}%
\ifx\x\y   
\gdef\SetFigFont#1#2#3{%
  \ifnum #1<17\tiny\else \ifnum #1<20\small\else
  \ifnum #1<24\normalsize\else \ifnum #1<29\large\else
  \ifnum #1<34\Large\else \ifnum #1<41\LARGE\else
     \huge\fi\fi\fi\fi\fi\fi
  \csname #3\endcsname}%
\else
\gdef\SetFigFont#1#2#3{\begingroup
  \count@#1\relax \ifnum 25<\count@\count@25\fi
  \def\x{\endgroup\@setsize\SetFigFont{#2pt}}%
  \expandafter\x
    \csname \romannumeral\the\count@ pt\expandafter\endcsname
    \csname @\romannumeral\the\count@ pt\endcsname
  \csname #3\endcsname}%
\fi
\fi\endgroup
\begin{picture}(471,255)(22,320)
\put( 67,459){\makebox(0,0)[lb]{\smash{\SetFigFont{8}{9.6}{rm}1}}}
\put(199,462){\makebox(0,0)[lb]{\smash{\SetFigFont{8}{9.6}{rm}
\special{ps: gsave}Characters of KM
\special{ps: grestore}}}}
\put( 82,459){\makebox(0,0)[lb]{\smash{\SetFigFont{8}{9.6}{rm}
Spherical functions}}}
\put( 58,556){\makebox(0,0)[lb]{\smash{\SetFigFont{8}{9.6}{rm}Im}}}
\put( 97,556){\makebox(0,0)[lb]{\smash{\SetFigFont{8}{9.6}{sc}
\special{ps: gsave}Representation theory of Lie groups, 
Lie algebras, and Kac-Moody algebras\special{ps: grestore}}}}
\put(182,459){\makebox(0,0)[lb]{\smash{\SetFigFont{8}{9.6}{rm}2}}}
\put(310,459){\makebox(0,0)[lb]{\smash{\SetFigFont{8}{9.6}{rm}
$[V_\lambda\otimes V_\mu:V\nu]$}}}
\put(415,459){\makebox(0,0)[lb]{\smash{\SetFigFont{8}{9.6}{rm}4}}}
\put(477,485){\makebox(0,0)[lb]{\smash{\SetFigFont{8}{9.6}{rm}Re}}}
\put(427,459){\makebox(0,0)[lb]{\smash{\SetFigFont{8}{9.6}{rm}
$[M_\lambda:L_\mu]$}}}
\put(311,448){\makebox(0,0)[lb]{\smash{\SetFigFont{8}{9.6}{rm}
\special{ps: gsave}(irreps of $\dim<\infty$)
\special{ps: grestore}}}}
\put(427,450){\makebox(0,0)[lb]{\smash{\SetFigFont{8}{9.6}{rm}
\special{ps: gsave}(induced: irreps)
\special{ps: grestore}}}}
\put(299,459){\makebox(0,0)[lb]{\smash{\SetFigFont{8}{9.6}{rm}3}}}
\put(205,450){\makebox(0,0)[lb]{\smash{\SetFigFont{8}{9.6}{rm}algebras}}}
\end{picture}} 
\vspace{-1.5in} 
\caption{Representation Theory} 
\label{figure2} 
\end{center} 
\end{figure}

Let us try to project representation theory 
on the real axis. In Figure \ref{figure2} 
we focus on Lie groups, Lie algebras, and  Kac\~Moody algebras, 
omitting the arithmetic direction (ad\`eles and automorphic forms). 
The theory of special functions, arithmetic, 
and related combinatorics 
are the classical objectives of representation theory.

\smallskip 
Without going into detail and giving 
exact references, the following are brief explanations. 

\smallskip 
\begin{description} 
\item[($1$)] I mean the zonal spherical functions on 
$K\backslash G/K$ 
for maximal compact $K$ in a semisimple Lie group $G$. 
The modern theory was started by 
Berezin, Gelfand, and others in the early 1950s 
and then developed significantly by Harish-Chandra \cite{HC}. 
Lie groups greatly helped to make the classical theory 
multidimensional, although they did not prove to be 
very useful for the hypergeometric function. 
\smallskip 

\item[($2$)] The characters of Kac\~Moody (KM) algebras are not 
far from the products of classical one-dimensional 
$\theta$\~functions 
and can be introduced without representation theory 
(Looijenga, Kac, Saito). See \cite{Lo}. 
However, it is a new and important class of special 
functions with various applications. Representation theory 
explains some of their properties, but not all. 
\smallskip 

\item[($3$)] This arrow gives many 
combinatorial formulas. 
Decomposing tensor products of finite dimensional 
representations of compact Lie groups and related problems 
were the focus of representation theory in the 1970s 
and early 1980s. They are still important, 
but representation theory moved toward 
infinite dimensional objects. 
\smallskip 

\item[($4$)] Calculating the multiplicities of irreducible 
representations of Lie algebras in the BGG\~Verma modules or 
other induced representations belongs to conceptual mathematics. 
The Verma modules were designed as a technical tool for 
the Weyl character formula (which is ``real"). 
It took time to understand that these multiplicities 
are ``real" too, with strong analytic and modular aspects. 
\end{description} 

\smallskip 
\subsection{New vintage} 
\begin{figure}[htbp] 
\begin{center} 
\resizebox{\textwidth}{!}{\begin{picture}(0,0)%
\special{psfile=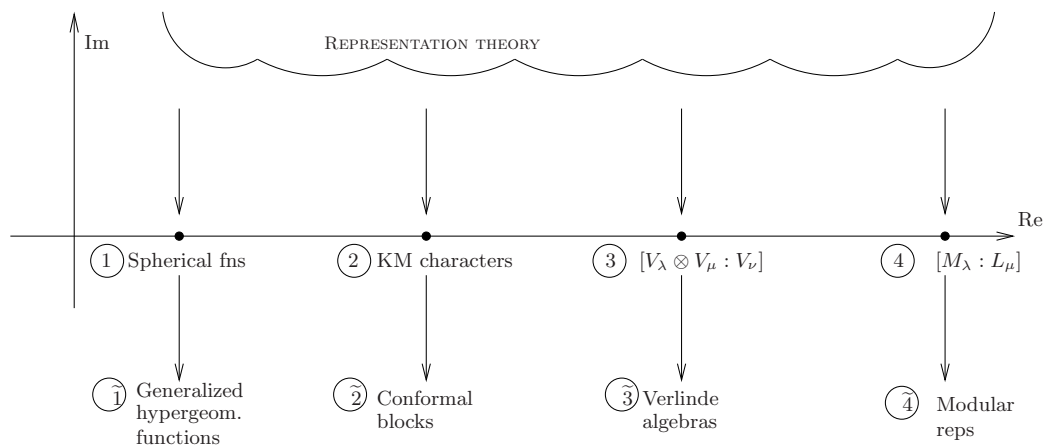}%
\end{picture}%
\setlength{\unitlength}{0.01250000in}%
\begingroup\makeatletter\ifx\SetFigFont\undefined
\def\x#1#2#3#4#5#6#7\relax{\def\x{#1#2#3#4#5#6}}%
\expandafter\x\fmtname xxxxxx\relax \def\y{splain}%
\ifx\x\y   
\gdef\SetFigFont#1#2#3{%
  \ifnum #1<17\tiny\else \ifnum #1<20\small\else
  \ifnum #1<24\normalsize\else \ifnum #1<29\large\else
  \ifnum #1<34\Large\else \ifnum #1<41\LARGE\else
     \huge\fi\fi\fi\fi\fi\fi
  \csname #3\endcsname}%
\else
\gdef\SetFigFont#1#2#3{\begingroup
  \count@#1\relax \ifnum 25<\count@\count@25\fi
  \def\x{\endgroup\@setsize\SetFigFont{#2pt}}%
  \expandafter\x
    \csname \romannumeral\the\count@ pt\expandafter\endcsname
    \csname @\romannumeral\the\count@ pt\endcsname
  \csname #3\endcsname}%
\fi
\fi\endgroup
\begin{picture}(500,213)(3,359)
\put( 155,552){\makebox(0,0)[lb]{\smash{\SetFigFont{8}{9.6}{sc}
\special{ps: gsave}Representation theory 
\special{ps: grestore}}}}
\put( 42,552){\makebox(0,0)[lb]{\smash{\SetFigFont{9}{10.8}{rm}Im}}}
\put( 63,446){\makebox(0,0)[lb]{\smash{\SetFigFont{9}{10.8}{rm}Spherical fns}}}
\put(184,446){\makebox(0,0)[lb]{\smash{\SetFigFont{9}{10.8}{rm}KM characters}}}
\put(452,446){\makebox(0,0)[lb]{\smash{\SetFigFont{9}{10.8}{rm}
$[M_\lambda:L_\mu]$}}}
\put(492,466){\makebox(0,0)[lb]{\smash{\SetFigFont{9}{10.8}{rm}
Re}}}
\put(170,446){\makebox(0,0)[lb]{\smash{\SetFigFont{9}{10.8}{rm}2}}}
\put(295,446){\makebox(0,0)[lb]{\smash{\SetFigFont{9}{10.8}{rm}3}}}
\put(435,446){\makebox(0,0)[lb]{\smash{\SetFigFont{9}{10.8}{rm}4}}}
\put(308,446){\makebox(0,0)[lb]{\smash{\SetFigFont{9}{10.8}{rm}
$[V_\lambda\otimes V_\mu:V_\nu]$}}}
\put( 50,446){\makebox(0,0)[lb]{\smash{\SetFigFont{9}{10.8}{rm}1}}}
\put(184,379){\makebox(0,0)[lb]{\smash{\SetFigFont{9}{10.8}{rm}Conformal}}}
\put(456,376){\makebox(0,0)[lb]{\smash{\SetFigFont{9}{10.8}{rm}Modular}}}
\put( 52,379){\makebox(0,0)[lb]{\smash{\SetFigFont{9}{10.8}{rm}
$\widetilde 1$}}}
\put(168,379){\makebox(0,0)[lb]{\smash{\SetFigFont{9}{10.8}{rm}
$\widetilde 2$}}}
\put(300,379){\makebox(0,0)[lb]{\smash{\SetFigFont{9}{10.8}{rm}
$\widetilde 3$}}}
\put(436,376){\makebox(0,0)[lb]{\smash{\SetFigFont{9}{10.8}{rm}
$\widetilde 4$}}}
\put(313,379){\makebox(0,0)[lb]{\smash{\SetFigFont{9}{10.8}{rm}Verlinde}}}
\put( 67,383){\makebox(0,0)[lb]{\smash{\SetFigFont{9}{10.8}{rm}Generalized }}}
\put( 67,372){\makebox(0,0)[lb]{\smash{\SetFigFont{9}{10.8}{rm}hypergeom.}}}
\put( 67,360){\makebox(0,0)[lb]{\smash{\SetFigFont{9}{10.8}{rm}functions}}}
\put(184,368){\makebox(0,0)[lb]{\smash{\SetFigFont{9}{10.8}{rm}blocks}}}
\put(313,368){\makebox(0,0)[lb]{\smash{\SetFigFont{9}{10.8}{rm}algebras}}}
\put(456,364){\makebox(0,0)[lb]{\smash{\SetFigFont{9}{10.8}{rm}reps}}}
\end{picture}} 
\caption{New Vintage} 
\label{figure3} 
\end{center} 
\end{figure}

Figure \ref{figure3} is an 
update of Figure \ref{figure2}. We add 
the results which were obtained in the 1980s and 1990s, 
inspired mainly by a breakthrough in  mathematical physics, 
although mathematicians had their own strong reasons to 
study generalized hypergeometric functions and modular 
representations. 

\begin{description} 
\item[($\tilde{1}$)] These functions are the key
in the theory of special function, both
differential and difference. 
The interpretation and generalization 
of the hypergeometric functions via representation 
theory was an important problem of the so-called Gelfand program 
and remained unsolved for almost three decades. 
\smallskip 

\item[($\tilde{2}$)] Actually, the conformal blocks belong to the 
(conceptual) imaginary axis as well as their kin, the 
$\tau$\~function. 
However, they extend the hypergeometric functions, theta functions, 
and Selberg's integrals. They attach 
the hypergeometric function to representation theory, 
but affine Hecke algebras serve this purpose better. 
\smallskip 
 
\item[($\tilde{3}$)] The Verlinde algebras were born from 
conformal field theory. They are formed by 
integrable representations of 
Kac\~Moody algebras of a  given level with ``fusion" 
instead of tensoring. These algebras can be also defined using 
quantum groups at roots of unity and interpreted 
via operator algebras. 
\smallskip 

\item[($\tilde{4}$)] Whatever you may think about the ``reality" of 
$[M_\lambda:L_\mu]$, these multiplicities are connected 
with the representations of Lie groups and Weyl groups 
over finite fields (modular representations). 
Nothing can be more real than finite fields! 
\end{description} 

\smallskip 
\subsection{Hecke algebras} 

The Hecke operators and later the Hecke algebras were introduced 
in the theory of modular forms, actually in the theory 
of $GL_2$ over the $p$\~adic numbers. In spite of their $p$\~adic 
origin, they appeared to be directly 
connected with the $K$\~theory of the 
{\em complex} flag varieties \cite{KL1} and, more recently, 
with the Harish-Chandra theory. It suggests that finite 
and $p$\~adic fields are of greater fundamental importance 
for mathematics and physics than we think. 

Concerning the 
great potential of $p$\~adics, let me mention the 
following three well-known confirmations: 
\smallskip 

\noindent (i) The Kubota\~Leopold $p$\~adic zeta function, 
which is a $p$\~adic 
analytic continuation of the values of the classical Riemann 
zeta function at negative integers. 

\noindent (ii) My theorem about ``switching 
local invariants" based on the $p$\~adic uniformization 
(Tate\~Mumford) 
of the modular curves which come from the quaternion algebras. 

\noindent (iii) The theory of $p$\~adic strings due to Witten, 
which is based 
on the similarity of the Frobenius automorphism in arithmetic to 
the Dirac operator. 
\medskip 

{\bf Observation.} {\em 
The real projection of representation theory 
goes through Hecke-type algebras.} 
\smallskip 

\begin{figure}[htbp] 
\begin{center} 
\resizebox{\textwidth}{!}{\begin{picture}(0,0)%
\special{psfile=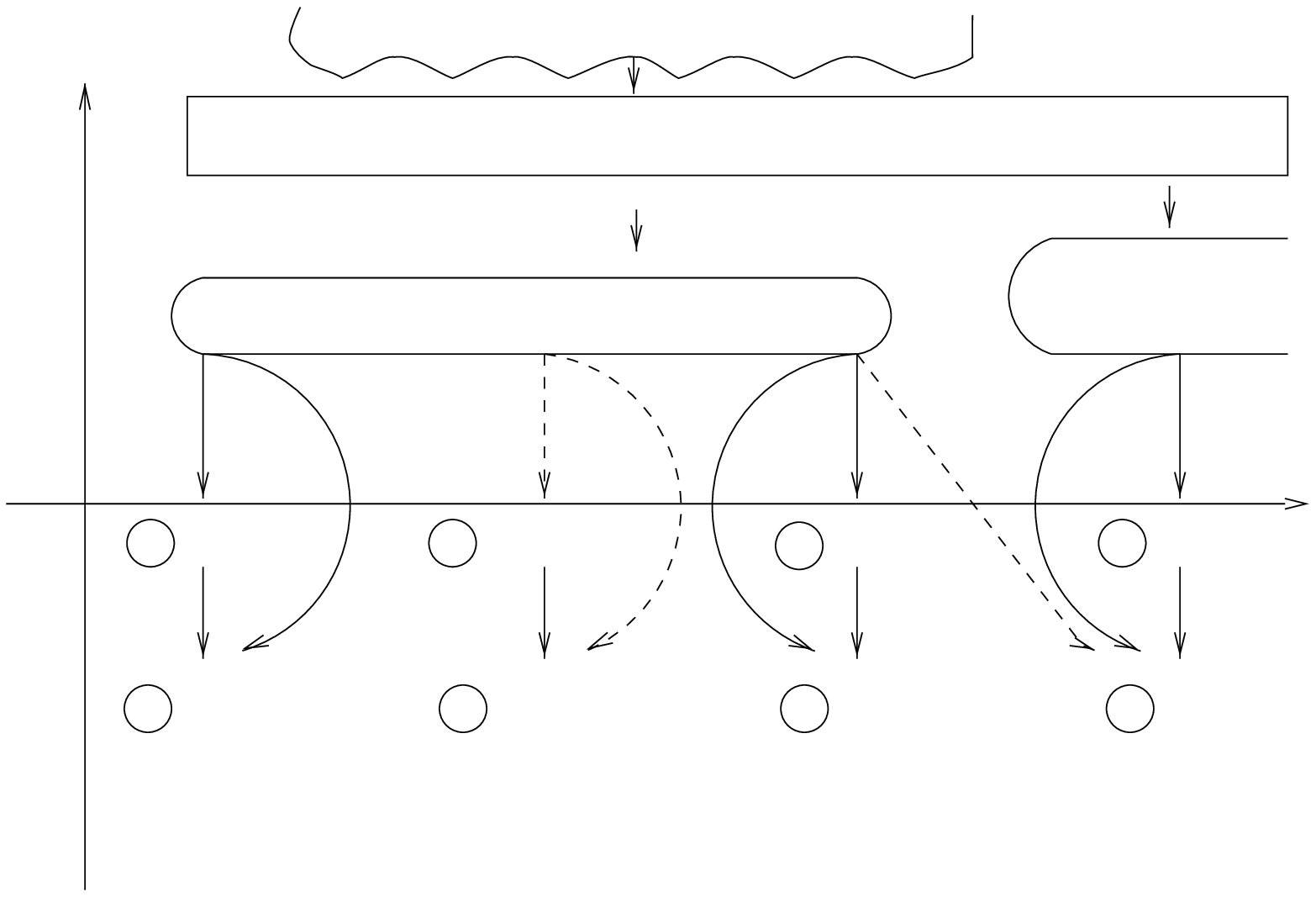}%
\end{picture}%
\setlength{\unitlength}{0.01250000in}%
\begingroup\makeatletter\ifx\SetFigFont\undefined
\def\x#1#2#3#4#5#6#7\relax{\def\x{#1#2#3#4#5#6}}%
\expandafter\x\fmtname xxxxxx\relax \def\y{splain}%
\ifx\x\y   
\gdef\SetFigFont#1#2#3{%
  \ifnum #1<17\tiny\else \ifnum #1<20\small\else
  \ifnum #1<24\normalsize\else \ifnum #1<29\large\else
  \ifnum #1<34\Large\else \ifnum #1<41\LARGE\else
     \huge\fi\fi\fi\fi\fi\fi
  \csname #3\endcsname}%
\else
\gdef\SetFigFont#1#2#3{\begingroup
  \count@#1\relax \ifnum 25<\count@\count@25\fi
  \def\x{\endgroup\@setsize\SetFigFont{#2pt}}%
  \expandafter\x
    \csname \romannumeral\the\count@ pt\expandafter\endcsname
    \csname @\romannumeral\the\count@ pt\endcsname
  \csname #3\endcsname}%
\fi
\fi\endgroup
\begin{picture}(503,339)(18,320)
\put(445,451){\makebox(0,0)[lb]{\smash{\SetFigFont{9}{10.8}{rm}4}}}
\put(415,549){\makebox(0,0)[lb]{\smash{\SetFigFont{9}{10.8}{rm}
Kazhdan-Lusztig}}}
\put(349,492){\makebox(0,0)[lb]{\smash{\SetFigFont{9}{10.8}{rm}$c$}}}
\put(314,497){\makebox(0,0)[lb]{\smash{\SetFigFont{9}{10.8}{rm}$\tilde c$}}}
\put( 74,497){\makebox(0,0)[lb]{\smash{\SetFigFont{9}{10.8}{rm}$a$}}}
\put(143,506){\makebox(0,0)[lb]{\smash{\SetFigFont{9}{10.8}{rm}$\tilde a$}}}
\put(477,497){\makebox(0,0)[lb]{\smash{\SetFigFont{9}{10.8}{rm}$d$}}}
\put(456,451){\makebox(0,0)[lb]{\smash{\SetFigFont{9}{10.8}{rm}
$[M_\lambda:L_\mu]$}}}
\put(192,387){\makebox(0,0)[lb]{\smash{\SetFigFont{9}{10.8}{rm}$\tilde 2$}}}
\put(208,387){\makebox(0,0)[lb]{\smash{\SetFigFont{9}{10.8}{rm}Conformal}}}
\put( 89,387){\makebox(0,0)[lb]{\smash{\SetFigFont{9}{10.8}{rm}Hypergeom.}}}
\put(361,515){\makebox(0,0)[lb]{\smash{\SetFigFont{9}{10.8}{rm}?!}}}
\put(235,501){\makebox(0,0)[lb]{\smash{\SetFigFont{9}{10.8}{rm}?}}}
\put(279,498){\makebox(0,0)[lb]{\smash{\SetFigFont{9}{10.8}{rm}?}}}
\put(302,515){\makebox(0,0)[lb]{\smash{\SetFigFont{9}{10.8}{rm}!}}}
\put(267,515){\makebox(0,0)[lb]{\smash{\SetFigFont{9}{10.8}{rm}$\tilde b$}}}
\put(426,536){\makebox(0,0)[lb]{\smash{\SetFigFont{9}{10.8}{rm}polynomials}}}
\put( 53,630){\makebox(0,0)[lb]{\smash{\SetFigFont{9}{10.8}{rm}Im}}}
\put(204,451){\makebox(0,0)[lb]{\smash{\SetFigFont{9}{10.8}{rm}KM-characters}}}
\put( 89,451){\makebox(0,0)[lb]{\smash{\SetFigFont{9}{10.8}{rm}Spherical fns}}}
\put( 89,372){\makebox(0,0)[lb]{\smash{\SetFigFont{9}{10.8}{rm}fns}}}
\put(218,372){\makebox(0,0)[lb]{\smash{\SetFigFont{9}{10.8}{rm}blocks}}}
\put(139,536){\makebox(0,0)[lb]{\smash{\SetFigFont{9}{10.8}{rm}
Macdonald theory, double Hecke algebras}}}
\put(204,501){\makebox(0,0)[lb]{\smash{\SetFigFont{9}{10.8}{rm}$b$}}}
\put(399,497){\makebox(0,0)[lb]{\smash{\SetFigFont{9}{10.8}{rm}$\tilde d$}}}
\put(462,387){\makebox(0,0)[lb]{\smash{\SetFigFont{9}{10.8}{rm}Modular}}}
\put(341,387){\makebox(0,0)[lb]{\smash{\SetFigFont{9}{10.8}{rm}Verlinde}}}
\put(492,477){\makebox(0,0)[lb]{\smash{\SetFigFont{9}{10.8}{rm}Re}}}
\put(467,375){\makebox(0,0)[lb]{\smash{\SetFigFont{9}{10.8}{rm}reps}}}
\put(341,376){\makebox(0,0)[lb]{\smash{\SetFigFont{9}{10.8}{rm}algebras}}}
\put(187,451){\makebox(0,0)[lb]{\smash{\SetFigFont{9}{10.8}{rm}2}}}
\put(322,387){\makebox(0,0)[lb]{\smash{\SetFigFont{9}{10.8}{rm}$\tilde 3$}}}
\put( 73,389){\makebox(0,0)[lb]{\smash{\SetFigFont{9}{10.8}{rm}$\tilde 1$}}}
\put(332,451){\makebox(0,0)[lb]{\smash{\SetFigFont{9}{10.8}{rm}
$[V_\lambda\otimes V_\mu:V_\nu]$}}}
\put( 74,451){\makebox(0,0)[lb]{\smash{\SetFigFont{9}{10.8}{rm}1}}}
\put(445,388){\makebox(0,0)[lb]{\smash{\SetFigFont{9}{10.8}{rm}$\tilde 4$}}}
\put(320,450){\makebox(0,0)[lb]{\smash{\SetFigFont{9}{10.8}{rm}3}}}
\put(213,646){\makebox(0,0)[lb]{\smash{\SetFigFont{9}{10.8}{rm}
Representation \ \ theory}}}
\put(192,605){\makebox(0,0)[lb]{\smash{\SetFigFont{9}{10.8}{rm}Representation 
 theory of Hecke algebras}}}
\end{picture}} 
\caption{Hecke Algebras} 
\label{figure4} 
\end{center} 
\end{figure}

The arrows in Figure \ref{figure4} are as follows.

\begin{description} 
\item[($a$)] This arrow is the most recognized now. 
Quite a few aspects of the Harish-Chandra theory 
in the zonal case were covered by representation theory of 
the degenerate (graded) affine Hecke algebras, introduced 
in \cite{L} (\cite{Dr} for $GL_n$). 
The radial parts of the invariant differential operators 
on symmetric spaces, the hypergeometric functions 
and their generalizations arise directly from these 
algebras \cite{Ch4}. 

The difference theory appeared even more promising. 
It was demonstrated in \cite{C3} 
that the $q$\~Fourier transform is 
self-dual like the classical Fourier and Hankel transforms, 
but not the Harish-Chandra transform. 
There are connections 
with the quantum groups and quantum symmetric spaces 
(Noumi, Olshansky, and others; see \cite{No}). 
However, the double Hecke algebra technique is simpler and 
more powerful. 
\smallskip 

\item[($b$)] The conformal blocks 
are solutions of the KZ\~Bernard equation 
(KZB). The double Hecke algebras lead to certain elliptic 
generalizations of the Macdonald polynomials \cite{Ch12,Ch13,C1} 
(other approaches are in \cite{EK,Ch12,FV2}, 
and the recent \cite{Rai}). 
These algebras govern the monodromy of the KZB equation 
and ``elliptic" 
Dunkl operators (Kirillov Jr., Felder\~Tarasov\~Varchenko, 
and the author). 

The monodromy map is the inverse of arrow ($\tilde{b}$). 
The simplest examples are directly related to 
the Macdonald polynomials and those at roots of unity. 
\smallskip 

\item[(${c}$)] Hecke algebras and their affine generalizations 
give a new approach to the classical combinatorics, 
including the characters of the compact Lie groups. 
The natural setting here is the theory of 
the Macdonald polynomials, although the analytic theory 
seems more challenging. 

Concerning ($\tilde{c}$), the Macdonald polynomials at 
the roots of unity give a simple approach to the Verlinde algebras 
\cite{Ki,C3,C4}. The use of the nonsymmetric Macdonald polynomials 
here is an important development. Generally, these polynomials 
are beyond the Lie and Kac\~Moody theory, although they 
are connected with the Heisenberg\~Weyl 
and $p$\~adic Hecke algebras. 
\smallskip 

\item[($d$)] This arrow is the Kazhdan\~Lusztig 
conjecture proved by Brylinski\~Kashi\-wa\-ra and 
Bei\-lin\-son\~Bern\-stein 
and then generalized to the Kac\~Moo\-dy case by 
Kashiwara\~Tanisaki. 

By ($\tilde{d}$), I mean the modular Lusztig conjecture 
(partially) proved by Anderson, Jantzen, and Soergel. 
There is recent significant progress due to Bez\-ru\-kav\-ni\-kov. 
\smallskip 

The arrow from the Macdonald theory to modular 
representations is marked by ``\,?!\,." 
It seems to be the most challenging 
now (there are already first steps in this direction). 
It is equivalent to extending the Verlinde algebras 
and their nonsymmetric variants 
from the alcove (the restricted category of representations 
of Lusztig's quantum group) to the parallelogram 
(all representations). 

If such an extension exists, it would give a $k$\~extension of 
Lusztig's conjectures, formulas for the modular characters 
(not only those for the multiplicities), a description of modular 
representations for arbitrary Weyl groups, and more. 
\end{description} 

\setcounter{equation}{0} 
\section{KZ and Kac--Moody algebras} \label{SEC:KZe}
In this section we comment on the role of the Kac\~Moody algebras 
and their relations (real and imaginary) 
to the spherical functions and the double Hecke algebras. 

\smallskip 
\subsection{Fusion procedure} 
I think that the penetration of double Hecke algebras into 
the fusion procedure and related problems 
of the theory of Kac\~Moody algebras is a convincing 
demonstration of their potential. 
The fusion procedure was introduced for the first time 
in \cite{ChFun}. On the physics side, let me also mention 
a contribution of Louise Dolan. 

Given an integrable representation of the $n$\~th power of 
a Kac\~Moody algebra  and two sets of points on a Riemann surface 
($n$ points and $m$ points), I constructed an integrable 
representation of the $m$\~th power of the same Kac\~Moody algebra.  
The construction does not change the ``global" central charge, 
the sum of the local central charges over the components. 
It was named later ``fusion procedure." 

I missed that in the special case of this correspondence, 
when $n=2$ and $m=1,$ 
the multiplicities of irreducibles in the resulting representation 
are the structural constants of a certain commutative algebra, 
the Verlinde algebra \cite{Ve}. 

Now we know that the Verlinde algebra and 
all its structures can be readily extracted from 
the simplest representation of the double affine 
Hecke algebra at roots of unity. Thus the 
Kac\~Moody algebras are undoubtedly connected 
with the double Hecke algebras. 

Double Hecke algebras dramatically simplify and generalize 
the algebraic theory of Verlinde algebras, 
including the inner product and 
the (projective) action of $PSL(2,\Z),$ however, excluding 
the integrality and positivity 
of the structural constants. The latter properties require $k=1$ 
and are closely connected with the Kac\~Moody interpretation 
(although they can also be checked directly). 

\smallskip 
I actually borrowed the fusion procedure from 
Y.~Ihara's papers ``On congruence monodromy problem." 
A similar construction is a foundation of  his theory. 
I changed and added some things 
(the central charge has no counterpart 
in his theory), but the procedure is basically the same. 
Can we go back and define Verlinde algebras in ad\`eles' setting? 

\smallskip 
\subsection{Symmetric spaces} 
The classification of Kac\~Moody algebras very much resembles 
that of symmetric spaces. 
See \cite{K}, \cite{Hel}. It is not surprising, 
because the key technical point in both theories is 
the description of the involutions and automorphisms of finite 
order for the semisimple finite dimensional Lie algebras. 
The classification lists are similar {\em but do not coincide.} 
For instance, the $BC_n$\~symmetric spaces have no 
Kac\~Moody counterparts. 
Conversely, the KM algebra of type, say, $D_4^{(3)}$ is 
not associated 
(even formally) with any symmetric space. Nevertheless 
one could hope 
that this parallelism is not incidental. 

Some kind of correspondence can 
be established using the isomorphism of 
the quantum many-body problem  \cite{C,Su,HO1}, 
a direct generalization of 
the Harish-Chandra theory, and the affine KZ equation. 
The isomorphism was found by A.~Matsuo and developed 
further in my papers. 
It holds when the parameter $k,$ given in terms of 
the root multiplicity in the context of symmetric spaces, 
is an arbitrary complex number. In the Harish-Chandra theory, 
it equals $1/2$ for $SL_2(\R)/SO_2$, $1$ in the so-called group case 
$SL_2(\C)/SU_2$, and $k=2$ for the $Sp_2.$ 
The $k$\~generalized spherical functions are 
mainly due to Heckman and Opdam. 

\smallskip 
Once $k$ was made an arbitrary number, it could 
be expected a counterpart of the central 
charge $c$, the level, 
in the theory of Kac\~Moody algebras.
Indeed, it has some geometric meaning. However, generally,
it is {\em not} connected with the central charge.
Indeed, the number of independent $k$\~parameters 
can be from $1$ ($A,D,E$) to $5$ ($C^\vee C$, 
the so-called Koornwinder case), but 
we have only one (global) central element $c$ 
in the Kac\~Moody theory. Also, the 
$k$\~spherical functions are eigenfunctions of differential operators 
generalizing the radial parts of the invariant operators 
on symmetric spaces. These operators 
have no counterparts for the Kac\~Moody characters. 
Also, the spherical functions are orthogonal polynomials;
the Kac\~Moody charactes are not.
In addition, the latter are of elliptic type, the spherical 
functions are of trigonometric type. 

The elliptic quantum many-body 
problem (QMBP) gives a kind of theory of spherical functions 
in the Kac\~Moody setting (at critical level). 
However, it supports the {\em unification} of $c$ and $k$ 
rather than the {\em correspondence} between them. 

The elliptic QMBP in the $GL_N$\~case was introduced 
by Olshanetsky and Perelomov \cite{OP}. The 
classical root systems were considered in the paper 
\cite{OOS}. The 
Olshanetsky\~Perelomov operators for arbitrary root systems were 
constructed in \cite{Ch12}. 
\smallskip 

We see that an exact match cannot be expected. However, 
a map from the Kac\~Moody algebras to 
spherical functions exists. It is for $GL_N$ only and not 
exactly for the KM characters, but it does exist. 

\smallskip 
\subsection{KZ and {\em r}--matrices} 
The KZ equation  is the system of differential equations 
for the matrix elements (using physical terminology,
the correlation functions) of 
the representations of the Kac\~Moody algebras 
in the $n$\~point case. The matrix elements are simpler to 
deal with than the characters. For instance, they satisfy 
differential equations with respect to the positions 
of the points. 

The most general ``integrable" case, 
is described by the so-called $r$\~matrix Kac\~Moody algebras 
from \cite{ChTau} and 
the corresponding {\em $r$\~matrix KZ equations} 
introduced in \cite{Ch0}. 
 
It was observed in the latter paper that the classical Yang\~Baxter 
equation can be interpreted as the compatibility 
of the corresponding KZ system, which dramatically 
enlarged the number of examples. An immediate application 
was a new class of KZ equations with {\em trigonometric} and 
{\em elliptic} dependence on the points. 
 
It was demonstrated in \cite{Ch0} that the abstract 
$\tau$\~function, also called the coinvariant, 
is a generic solution 
of the $r$\~matrix KZ with respect to the action 
of the Sugawara $(-1)$\~operators. 
\smallskip

More generally, the $r$\~matrices and the corresponding 
KZ equations attached to arbitrary{\em root systems} were 
defined in \cite{Ch0}. For instance, the dependence on 
the points is via 
the differences (the $A$\~case) of the points
and also via the {\em sums} for $B,C,D$.
The $BC$\~case is directly related to the so-called 
{\em reflection equations} introduced in \cite{ChTmp}. 

The results due to Drinfeld\~Kohno on the monodromy 
of the KZ equations (see \cite{Ko}) 
can be extended to the $r$\~matrix equations. 
In some cases, the monodromy can be calculated 
explicitly, for instance, for the 
{\em affine KZ} \cite{Ch0,Ch2,Ch1}. 

\smallskip 
\subsection{Integral formulas for KZ} 
The main applications of the interpretation of 
KZ as a system of equations for the coinvariant were:\, 
(i) a simplification of the algebraic part 
of the Schechtman\~Varchenko construction \cite{VSch} 
of integral formulas for the rational KZ, \,
(ii) a generalization of their formulas 
to the trigonometric case \cite{Ch3}. 
\smallskip 
Paper \cite{VSch} is based on direct algebraic
considerations without using the theory of 
Kac\~Moody algebras. 

There is another important ``integrable" case, 
the so-called Knizh\-nik\~Za\-mo\-lo\-d\-chi\-kov\~Ber\-nard 
equation usually denoted by KZB \cite{Ber, FW}. 
It can be obtain in the 
same abstract manner as a system of differential equations 
for the corresponding ``elliptic" coinvariant. There must be 
an implication of this fact toward the integral formulas for KZB, 
but this has not been checked so far. 

We do not discuss the integral formulas 
for KZB in this paper, as well as the integral formulas for QKZ, 
the quantum Knizhnik\~Zamolodchikov 
equation. See, e.g., \cite{TV}, \cite{FV}, and \cite{FTV}. 

Generally, the KZ equations can be associated 
with arbitrary algebraic curves. Then 
they involve the derivatives with respect to the moduli of 
curves and vector bundles. However, in this generality, the 
resulting equations are non-integrable in any reasonable sense. 

\smallskip 
Summarizing, 
we have the following major cases, when 
the Knizhnik\~Zamolodchikov\ 
equation have integral formulas, 
reasonably simple monodromy representations, special symmetries, 
and other important properties: 

\noindent 
{(a)} the KZ for Yang's rational $r$\~matrix (see \cite{VSch}),\\ 
{(b)} the trigonometric KZ equation introduced in \cite{Ch3},  \\ 
{(c)} the elliptic KZ\~Bernard equation (see \cite{Ber,FW}). 
\smallskip 

Given a Lie algebra ${\mathfrak{g}}$,
one may define the integrand for 
the KZ integral formulas is derectly connected with 
the coinvariants of $U(\widehat{\mathfrak{g}})$ 
for the Weyl modules \cite{Ch3}. The contours (cycles) 
of integration  are 
governed by the quantum $U_q({\mathfrak{g}}).$ 
See \cite{FW1}, \cite{V} and references therein. 
We will not discuss the contours  
and the $q$\~topology of the configuration spaces in this paper. 

The later topic was started by Aomoto \cite{A} and seems an 
endless story. We have no satisfactory formalization of 
the $q$\~topology so far. It is especially needed for 
QKZ. Generally, in mathematics, the contours of integration 
(the homology) must be dual to the differential forms 
(the cohomology). It gives an approach to the problem.

We note that the integral KZ formulas are directly
connected with the equivalence of the 
$U( \widehat{\mathfrak{g}})_c$ and  the 
quantum group $U_q({\mathfrak{g}})$ due 
to Kazhdan, Lusztig, and Finkelberg (see \cite{KL2}). 
It is for a proper relation $c\leftrightarrow q.$ 

\smallskip 
\subsection{From KZ to spherical functions} 
Let us discuss what the 
integral formulas could give for the 
theory of spherical functions and its 
generalizations. There are natural limitations. 
\smallskip 

First, only the spherical functions of type $A$ may apper 
(for either choice of ${\mathfrak{g}}$) if we
begin with the KZ integral formulas of type $A$.

Second, one needs 
an $r$\~matrix KZ of  {\em trigonometric} type because the 
Harish-Chandra theory is on the torus. 

Third, only $\mathfrak{g}={\mathfrak{gl}}_N$ may result in 
{\em scalar} differential 
operators due to the analysis by Etingof and Kirillov Jr. 

Summarizing, the integral formulas
for the affine KZ (AKZ) of type $A$ are the major candidates.
The AKZ is {\em isomorphic} to the 
quantum many-body problem, that
is exactly the $k$\~Harish-Chandra theory \cite{M,C13}.
 
Note that the ``basic" trigonometric 
$n$\~point KZ taking values in the 
$0$\~weight component of 
$(\C^n)^{\otimes n},$ which is isomorphic to the group 
algebra $\C\BS_n$, must be considred for AKZ. 
The integral AKZ formula  
is likely to be directly connected 
with the Harish-Chandra formula. I did not check it, but 
calculations due to Mimachi, Felder, Varchenko confirm this. 
For instance, the dimension of the contours (cycles) 
of integration 
for such KZ is $n(n-1)/2,$ 
which coincides with that in the Harish-Chandra integral 
representation for spherical functions of type $A_{n-1}.$ 
His integral is  over  $K=SO_n\subset SL_n(\R)$. 
\smallskip 

Establishing a direct connection with the Harsh-Chandra integral 
representation for the spherical functions does not seem too 
difficult. However it is of obvious importance, because 
his formula is for {\em all} root systems, and one can use 
it as an initial point for the general theory 
of integral formulas of the KZ equations 
associated with root systems. 
\smallskip 
 
We note that the integral 
KZ formulas can be justified without Kac\~Moody 
algebras. A straightforward algebraic combinatorial
analysis is complicated  but possible \cite{VSch}. 
My proof is based on the Kac\~Moody coinvariant \cite{Ch3}. 
However, I use the Kac\~Moody algebras  
at the critical level only, as a technical tool,
and then extend the resulting 
formulas to all values of the center charge. 
\smallskip

There is another approach 
to the same integral formulas 
based on the coinvariant for the Wakimoto modules 
instead of that for the Weyl modules \cite{FFR}. 
The calculations 
with the  coinvariant are in fact similar to mine, 
but the Kac\~Moody theory are used to greater potential 
and the combinatorial part gets simpler. 
\smallskip 

Thus one may expect the desired relation between 
the conformal blocks and spherical functions at the level of 
integral formulas for the trigonometric 
KZ of type $A,$ i.e., in terms of the differences of the 
points, and with values in the
simplest representations of $\mathfrak{g}={\mathfrak{gl}}_n.$ 

I do not think that this correspondence is really general
and can be extended to arbitrary symmetric spaces, 
though it certainly indicates that there 
must be a {\em unification} 
that combines the Kac\~Moody and Harish-Chandra theories. 

\setcounter{equation}{0} 
\section{Double Hecke algebras} \label{SEC:Doublea}
Double affine Hecke algebras 
(DAHA) were initially designed to clarify the classical and quantum 
Knizhnik\~Za\-mo\-lod\-chi\-kov equation (for the simplest 
fundamental 
representation of $\mathfrak{g}=\mathfrak{gl}_N$) and analogous 
KZ equations for arbitrary Weyl groups. The first applications 
were to the Dunkl operators, differential and difference. 
The most natural way to introduce these operators
is via the induced representations of DAHA. 

Eventually, through applications to 
the theory of Macdonald polynomials, 
DAHAs led to a unification of the 
Harish-Chandra transform in the zonal case and the 
$p$\~adic spherical transform in one general $q$\~theory, 
which is one of the major applications. 

The DAHA\~Fourier transform 
depends on a parameter $k,$ which generalizes the 
root multiplicity in the Harish-Chandra theory. 
This parameter becomes $1/(c+g)$ in the KZ theory, 
if KZ is interpreted via Kac\~Moody algebras of
level $c$ and the isomorphism between
AKZ and QMBP is used. The parameter $q$ comes from 
the Macdonald polynomials and QKZ, the quantum KZ 
equation introduced by F.Smirnov and I.Frenkel\~Reshetikhin
\cite{Sm,FR} and then generalized to all root systems in
\cite{Ch5}.

The limiting cases as $q\to 1$ 
and $q\to \infty$ are, respectively, the Harish-Chandra and 
the $p$\~adic Macdonald\~Matsumoto spherical transforms. 

It is not just a unification of the latter transforms.
The $q$\~transform is {\em self-dual} in 
contrast to its predecessors. The self-duality collapses under 
the limits above. However, there is a limiting procedure 
preserving the self-duality. If $q\to 1$ and 
we represent the torus in the form of $\{q^x\}$ 
instead of $\{e^x\},$ 
then the correspomding limit becomes self-dual. 
It is the multidimensional generalization from \cite{O3,Du,J} 
of the classical Hankel transform in terms of the Bessel functions. 

The new $q$\~transform shares many properties with 
the generalized Hankel transform. For instance, 
the $q$\~Mehta\~Macdonald integrals, generalizing the classical
Gauss\~Selberg integrals, are in the focus 
of the new theory. 
The case when $q$ is a root of unity is of obvious importance 
because of the applications to Verlinde algebras,
Gaussian sums, diagonal coinvariants, and more.
\smallskip 

\subsection{Missing link?} 
There are reasons to consider DAHA as a 
candidate for the ``missing link" between representation theory 
and theory of special functions. Let me explain why something 
seems missing. 
 
The representation theory of finite dimensional
Lie groups mainly 
serves multi\-dimensional functions 
and gives only a little for the one-dimensional functions 
(with a reservation about the arithmetic direction). 
The fundamental objects of 
the modern representation theory are $\widehat{\mathfrak{sl}}_2$ and 
its quantum counterpart. They have important applications 
to the hypergeometric and theta functions, but managing 
these functions via $\widehat{\mathfrak{sl}}_2$ is far from 
simple (especially in the $q$\~case). 

Generally, the Kac\~Moody theory is too algebraic 
to be used in developing the theory of functions. 
There were interesting attempts, but 
still we have no consistent harmonic Kac\~Moody analysis. 
The theory of operator algebras can be considered as an analytic
substitute of the Kac\~Moody theory, but it does not  help much 
with the special functions either. 

On the other hand, the Heisenberg and Weyl algebras are 
directly related 
to the theory of special functions. Unfortunately 
they are too simple, 
as is the affine Hecke 
algebra of type $A_1.$ 
\smallskip 

The double affine Hecke algebra $\HH\,$ 
(``double H") of type $A_1$ seems just right 
to incorporate major classical special functions of
one variable into 
representation theory. 
\smallskip 

This algebra has a simple definition, but its representation theory 
is rich enough. I think that 
if we combine what is already known about $\HH_{A_1},$ 
its applications and representations, 
it would be a book 
as big as ``$SL_2$" by Lang or ``$GL(2)$" by Jacquet\~Langlands. 
One-dimensional DAHA (already) has surprisingly many applications. 

\begin{description} 
\item[($a$)] There are direct relations to 
$\mathfrak{sl}(2)$ and $\mathfrak{osp}(2|1)$ and their quantum 
counterparts. The rational degeneration of $\HH_{A_1}$ 
as $q\to 1$ is a quotient of $U(\mathfrak{osp}(2|1)).$ 
Rational DAHAs are very interesting; they are a kind of Lie 
algebras in the general $q$\~theory. 

\item[($b$)] The Weyl and Heisenberg algebras are its limits 
when $t,$ the second parameter of $\HH\,,$ tends to $1.$ 
For instance, the $N$\~dimensional representation of the 
Weyl algebra, as $q$ is a primitive $N$\~th root of unity, 
has a direct counterpart for $\HH ,$ namely, the nonsymmetric 
Verlinde algebra, with various applications, including 
a new theory of Gaussian sums. 

\item[($c$)] It covers the theory of Rogers' ($q$\~ultraspherical) 
polynomials. Its $C^\vee C_1$\~modi\-fi\-cation go\-verns major 
remarkable families of one-dimensio\-nal orthogonal polynomials 
and has applications to the Bessel functions 
(the rational degeneration) and to the 
classical and basic hypergeometric functions. 

\item[($d$)] Generally, the duality from \cite{VV} 
connects $\HH\,$ of type $A$ with the to\-ro\-i\-dal 
(double $q$\~Kac\~Moody) algebras in a sense of 
Ginzburg, Kapranov, and Vasserot. The $\HH_{A_1}$ 
has important applications to the representation 
theory of $\widehat{\mathfrak{sl}}_2$ and in the higher 
ranks via this duality. 

\item[($e$)] DAHAs also result from the $K$\~theory of 
affine flag varieties and are related to the $q$\~Schur 
algebras. The connections with the so-called double arithmetic 
must be mentioned too. Kapranov and then Gaitsgory\~Kazhdan 
interpreted $\HH_{A_1}$ as a Hecke algebra of 
the ``$p$\~adic\~loop group" of $SL(2).$ 

\item[($f$)] In the very first paper \cite{C15} on DAHA, 
its topological interpretation was used. 
The algebra $\HH_{A_1}$ is directly connected with 
the fundamental group of an elliptic 
curve with one puncture. This establishes a link to the 
Gro\-then\-dieck\~Belyi program in the elliptic case 
(Beilinson\~Levin) and to some other problems of modern arithmetic. 
 
\item[($g$)] The representation theory of $\HH_{A_1}$ is far from 
trivial, especially at roots of unity. For instance, the description 
of its center at roots of unity (and in the $C^\vee C_1$\~case)
leads to the quantization of the cubic surfaces  (Oblomkov). 
Some relation to the Fourier\~Mukai transform are expected.

\item[($h$)] DAHA unifies the 
Harish-Chandra spherical transform and the $p$\~adic 
Macdonald\~Matsumoto transforms, The corresponding 
$q$\~transform is self-dual.
For $A_1,$ it leads to a deep analytic theory of
the basic hypergeometric function.
\end{description} 

The theory is now essentially algebraic.
There are some results in the analytic direction 
in \cite{C5} (the construction of the general 
spherical functions), \cite{C10} (analytic continuations 
in terms of $k$ with applications to 
$q$\~counterparts of Riemann's zeta 
function), and \cite{KS1,KS2}. 
The latter two papers are devoted to 
analytic theory of the one-dimensional 
Fourier transform in terms of the basic hypergeometric 
function, a $C^\vee C_1$\~extension 
of the $q$\~transform considred in the paper.

\subsection{Gauss integrals and sums} 
The starting point of many mathematical and physical theories 
is the celebrated formula: 
\begin{align} 
&2\int_{0}^\infty e^{-x^2}x^{2k}\hbox{d}x\ =\ \Gamma (k+1/2), \ 
\Re k>-1/2. 
\label{(1)} 
\end{align} 
Let us give some examples. 
\vskip 0.2cm 

{(a)} Its generalization with the product of two 
Bessel functions added to the integrand or, equivalently, 
the formula for the Hankel transform of the 
Gaussian $e^{-x^2}$ multiplied by a Bessel function, 
is one of the main formulas in the classical theory of 
Bessel functions. 

{(b)} The following ``perturbation" for the same $\Re k>-1/2,$ 
\begin{align} 
&2\int_{0}^{\infty} 
(e^{x^2}+1)^{-1}x^{2k}\hbox{\, d}x = 
(1-2^{1/2-k})\Gamma (k+1/2)\zeta (k+1/2), \notag 
\end{align} 
is fundamental in analytic number theory. For instance, 
it readily gives the functional equation for $\zeta.$ 

{(c)} The multidimensional extension  due to Mehta \cite{Meh}, 
when we integrate 
over $\R^n$ with the measure 
$\prod_{1\le i<j\le n}(x_i-x_j)^{2k}$ instead of $ x^{2k},$ 
gave birth to the theory of matrix models with 
various applications in mathematics and physics. 
Its generalization to arbitrary roots (Macdonald and Opdam) is 
a major formula in the modern theory of Hankel transform. 

{(d)} Switching to the roots of unity, we have an 
equally celebrated 
Gauss formula: 
\begin{align} 
&\sum_{m=0}^{2N-1} e^{\frac{\pi m^2 }{ 2N}i}\ =\ 
(1+i)\sqrt{N},\ \ N\in \BN. \label{(2)} 
\end{align} 
It is a counterpart of (\ref{(1)}) at $k=0,$ 
although there is no direct connection. 
\smallskip 

To fully employ modern mathematics 
we need to go from the Bessel to 
the hypergeometric functions. 
In contrast to the former, the latter 
can be studied, interpreted, and generalized 
by a variety of methods from 
representation theory and algebraic geometry to integrable models 
and string theory. 

Technically, the measure $x^{2k}dx$ has 
to be replaced by $\sinh(x)^{2k}dx$ and the Hankel transform by 
the Harish-Chandra transform, to be more exact, by its 
$k$\~extension. 
However, the latter is no longer self-dual, 
the formula (\ref{(1)}) has no $\sinh$\~counterpart, and 
the Gaussian looses its Fourier\~invariance. Thus 
a straightforward substitution creates problems. 
We need a more fine-tuned approach. 

\subsection{Difference setup} 
The main observation is that 
the self-duality of the Hankel transform 
is restored for the kernel 
\begin{align} 
&\delta(x;q,k)\equal\prod_{j=0}^\infty 
\frac{(1-q^{j+2x})(1-q^{j-2x})}{ 
(1-q^{j+k+2x})(1-q^{j+k-2x})},\ 0<q<1, \ k\in \BC. 
\notag 
\end{align} 
Here $\delta,$ the Macdonald truncated theta function, is 
a certain unification of 
$\sinh(x)^{2k}$ and the measure in terms of the Gamma function 
serving the inverse Harish-Chandra transform ($A_1$). 
Therefore the self-duality of the resulting transform 
can be expected a priori. 

As to (\ref{(1)}), setting $q=\exp(-1/a),\ a>0,$ we have 
\begin{align} 
&(-i)\int_{-\infty i}^{\infty i}q^{-x^2}\delta(x;q,k) 
\hbox{\, d}x= 
{2\sqrt{a\pi}}\prod_{j=0}^\infty 
\frac{1-q^{j+k}}{ 
 1-q^{j+2k}},\ \ \Re k>0. 
\label{(3)} 
\end{align} 
Here both sides are well defined for all $k$ except for the poles 
but coincide only when $ \Re k>0.$ 

One can make (\ref{(3)}) entirely algebraic by replacing 
the Gaussian $\ga^{-1}=q^{-x^2}$ by its expansion 
$$\tilde{\ga}^{-1}\ =\ \sum_{-\infty}^{+\infty}q^{n^2/4}q^{nx} 
$$ 
and using Const Term$(\sum c_nq ^{nx})=c_0$ instead of 
the imaginary integration: 
\begin{align} 
&\hbox{Const\ Term\ }(\tilde{\ga}^{-1}\delta)= 
2\prod_{j=0}^\infty 
\frac{1-q^{j+k}}{ 
 1-q^{j+2k}}.\notag 
\end{align} 

{\bf Jackson integrals.} 
A promising feature of special  $q$\~functions is 
the possibility of replacing the integrals by sums over $\BZ^n$, 
the so-called Jackson integrals. 

Technically, we switch from the imaginary integration to that 
for the 
path, which  begins at\, $z=\epsilon i+\infty$, moves 
to the left till\, $\epsilon i$, then down 
through  the origin to\,  $-\epsilon i$, and then 
returns down the positive real axis to \, 
$-\epsilon i+\infty$ 
(for small $\epsilon$). Then we apply Cauchy's theorem under 
the assumption that $|\Im k|<2\epsilon,\ \Re k>0$. 

We obtain the following counterpart of (\ref{(3)}): 
\begin{align} 
&\sum_{j=0}^\infty q^{\frac{(k-j)^2}{ 4}} 
\frac{1-q^{j+k}}{ 
1-q^{k}} \prod_{l=1}^j 
\frac{1-q^{l+2k-1}}{ 
 1-q^{l}}\ = 
\notag\\ 
& q^{\frac{k^2}{ 4}}\prod_{j=1}^\infty 
\frac{(1-q^{j/2})(1-q^{j+k})(1+q^{j/2-1/4+k/2})(1+q^{j/2-1/4-k/2})} 
{(1-q^j)}, 
\notag \end{align} 
which is convergent for all $k.$ 
 
When $q=\exp(2\pi i/N)$ and $k$ is a positive integer 
$ \le N/2,$ we come to a  Gauss\~Selberg-type 
sum: 
\begin{align} 
&\sum_{j=0}^{N-2k} q^{\frac{(k-j)^2}{ 4}} 
\frac{1-q^{j+k}}{ 
1-q^{k}} \prod_{l=1}^j 
\frac{1-q^{l+2k-1}}{ 
 1-q^{l}}= 
\prod_{j=1}^k 
(1-q^{j})^{-1}\sum_{m=0}^{2N-1} q^{m^2/4}.\notag 
\end{align} 
The left-hand side 
resembles the so-called modular Gauss\~Selberg sums. 
However, the difference is dramatic. The modular 
sums are calculated in the finite fields and are 
embedded into roots of unity right before the final summation
\cite{Ev}. 
Our sums are defined entirely in cyclotomic fields. 
Substituting $k=[N/2],$ we  arrive at (\ref{(2)}).

\subsection{Other directions} 
There are other projects involving the double Hecke 
algebras. We will mention only some of them. 

\smallskip 
(1) {\em Macdonald's $q$\~conjectures} \cite{M1,M2,M3}. 
Namely, the constant term, norm, duality, and  evaluation 
conjectures \cite{BZ,Kad},\cite{Ch7,C2,C3}. 
See also \cite{Ao,I,M5,C1} about the discrete 
variant of the constant term conjecture, the Aomoto conjecture. 
My proof of the norm formula was based on the shift operators and 
is similar to that from  \cite{O1} in the differential case 
(the duality and evaluation conjectures 
collapse as $q\to 1$). I would add to this list the Pieri rules. 
As to the nonsymmetric Macdonald polynomials, 
the references are \cite{O2,M4,C4}. 
See also \cite{DS,Sa,Ma3}, and the recent 
\cite{St,Rai} about the case of $C^\vee C,$ 
the Koornwinder polynomials and generalizations. 

\smallskip 
(2) {\em $K$\~theoretic interpretation}. See papers 
\cite{KL1,KK}, then \cite{GG,GKV}, 
and the important recent paper \cite{Va}. 
The latter leads to the 
Langlands-type description of irreducible representations of 
double Hecke algebras. 
The case of generic parameters $q,t$ is 
directly connected with the affine theory \cite{KL1}. 
For the special parameters, the corresponding geometry 
becomes significantly more complicated. 
The Fourier transform remains unclear in this approach. 
I also mention here the strong Macdonald conjecture (Hanlon) 
and the recent \cite{FGT}.

\smallskip 
(3) {\em Elementary methods}. 
The theory of induced, semisimple, unitary, and 
spherical representations can be developed successfully without 
$K$\~theory. The main tool is the technique of intertwiners 
from \cite{C1}, which is 
similar to that for the affine Hecke algebras. 
The nonsymmetric Macdonald polynomials generate the 
simplest spherical 
representation, with the intertwiners serving as creation operators 
(the case of $GL$ is due to \cite{KnS}). For $GL,$ 
this technique gives a 
reasonably complete classification of irreducible representations 
similar to the theorem by Bernstein and Zelevinsky in the 
affine case. 
Relations to \cite{HO2} are expected. 

\smallskip 
(4) {\em Radial parts via Dunkl operators}. The main references are 
\cite{Du,H} and \cite{C13}. 
In the latter, it was observed that the trigonometric 
differential Dunkl operators can be obtained 
from the degenerate (graded) affine 
Hecke algebra from \cite{L} (\cite{Dr} for $GL_n$). 
The difference, elliptic, and 
difference-elliptic generalizations were introduced in 
\cite{Ch5,C15,Ch12,Ch13}. 
The connections with the KZ equation play an important role here 
\cite{M,C13,C14,C15}. 
The radial parts of the Laplace operators 
of symmetric spaces and 
their generalizations are 
symmetrizations of the Dunkl operators. 
The symmetric Macdonald polynomials are eigenfunctions 
of the difference radial parts. The nonsymmetric 
Opdam\~Macdonald polynomials appear as 
eigenfunctions of the Dunkl operators. 

\smallskip 
(5) {\em Harmonic analysis}. The Dunkl operators in the simplest 
rational-differential setup lead to the definition of the 
generalized Bessel functions and the generalized 
Hankel transform (see \cite{O3,Du2,J} and also \cite{He}). 
In contrast to the Harish-Chandra and the $p$\~adic 
spherical transforms, 
it is self-dual. The self-duality resumes in the 
difference setting \cite{C12,C5}. 
The Mehta\~Macdonald conjecture, directly related to the transform 
of the Gaussian, was checked in \cite{M1,O1} 
in the differential case, 
and extended in \cite{C5} to the difference case. It was used there 
to introduce the $q$\~spherical functions. Concerning the 
applications 
to the Harish-Chandra theory, see \cite{HO1,O2,C8} and also 
\cite{HS}. 

\smallskip 
(6) {\em Roots of unity}. The construction from \cite{C12} 
generalizes 
and, at the same time, simplifies the Verlinde algebras, including 
the projective action of $PSL(2,\BZ)$ 
(cf. \cite{K}, Theorem 13.8, and \cite{Ki}), 
the inner product, 
and a new theory of Gaussian sums. 
In \cite{C4}, the nonsymmetric Verlinde algebras 
were considered. The symmetric 
elements of such algebras form the 
$k$\~generalized Verlinde algebras. 
The simplest example is the classical $N$\~dimensional 
representation of the Weyl algebra at $q=\exp(2\pi i/N).$ 
Let me also mention the recent \cite{Go,C29} about the Haiman 
conjecture \cite{Hai} on the structure of the 
so-called diagonal coinvariants, 
which appeared to be directly connected with rational DAHAs 
and DAHAs at roots of unity. 

\smallskip 
(7) {\em Topology}. The group $PSL(2,\BZ)$ acts projectively 
on the double Hecke algebra itself. 
The best explanation (and proof) is based on the 
interpretation of this algebra as a quotient of the group algebra 
of the $\pi_1$ of the elliptic configuration space from \cite{C15}. 
The calculation of $\pi_1$ in the $GL$\~case 
is essentially due to \cite{Bi,Sc}. 
For arbitrary root systems, it is similar to that from \cite{Lek}, 
but our configuration space is different. Such $\pi_1$ governs 
the monodromy of the eigenvalue problem for 
the elliptic radial parts, the corresponding 
Dunkl operators, and the KZB equation. 
Switching to the roots of unity, the monodromy 
representation is the nonsymmetric Verlinde algebra; 
applications to the Witten\~Reshetikhin\~Turaev 
and Ohtshuki invariants are expected. 

\smallskip 
(8) {\em $GL$\~Duality}. The previous discussion  was 
about arbitrary root systems. In the case of $GL$, 
the theorem from \cite{VV} establishes the duality 
between the double Hecke 
algebras (actually its extension) and the $q$\~toroidal 
(double Kac\~Moody) algebras in the sense 
of Ginzburg\~Kapranov\~Vasserot. 
It generalizes the classical 
Schur\~Weyl duality, Jimbo's $q$\~duality for the nonaffine Hecke 
algebra \cite{Ji}, and the affine Hecke analogs 
from \cite{Dr,C9}. When the  center charge is nontrivial the duality 
explains the results from \cite{KMS} and 
\cite{STU}, which were extended by Uglov 
to irreducible representations of the Kac\~Moody 
$\hat{\gl}$ of arbitrary positive integral levels. 

(9) {\em Rational degeneration.}. 
The rational degeneration of the double affine Hecke 
algebra with trivial center charge ($q=1$) is directly 
related to the Calogero\~Moser varieties \cite{EG}. 
The rational degenerations of the double 
affine Hecke algebra \cite{CM,EG}, 
in a sense, play the role of the Lie algebras of the $q,t$\~DAHA. 
The trigonometric 
degeneration is also a sort of Lie algebra, but the rational 
one has the projective action of $PSL(2,\BZ)$ and other symmetries 
that make it closer to the general $q,t$\~DAHA. 
The theory of the rational DAHA and its connections 
appeared to be an interesting independent direction (see 
\cite{BEG,GGOR,Go,C29}.

\smallskip 
{\bf Some of the latest developments.} 
There are interesting papers \cite{GK,GK1} that continued 
earlier results by Kapranov \cite{Kap} towards using DAHA 
in the so-called double arithmetic, 
started by Parshin quite a few years ago. 
The general problem is to associate a double 
$p$\~adic Lie group with DAHA. 

Let me also mention the results by Ginzburg, Kapranov and Vasserot
concerning the interpretation of
DAHAs via ``Hecke correspondences" over 
special algebraic surfaces. Presumably  DAHAs 
are somehow related to the Nakajima surfaces;
at least, this is understood for the Calogero\~Moser
varieties \cite{EG}.

An important recent development is the establishment 
of the connection 
of DAHA with the Schur algebra in \cite{GGOR} and then in 
\cite{VV1}. It is proven in the latter paper 
that, under minor technical restrictions, 
DAHA of type $A$ is Morita equivalent to the quantum affine Schur 
algebra at roots of unity. 

The latest project so far employs DAHAs for the global 
quantization of algebraic surfaces. 
Oblomkov used the rank one $C^\vee C$\~DAHA 
to quantize the cubic surfaces. 
Then Etingof, Oblomkov, Rains \cite{EOR}
extended this approach to quantize 
the Del Pezzo surfaces via certain generalizations 
of DAHA. Relations to \cite{Rai} are expected. 

\medskip 

{\bf Related directions.} 
Let me also mention several directions that 
are not based on double Hecke algebras (and their variants) 
but have close relations to these algebras and to the 
theory of the symmetric Macdonald 
polynomials.  Mainly they are about the $GL$\~case and the 
classical root systems. 

(a) The spherical functions on $q$\~symmetric 
spaces  due to Noumi and others, which are related 
to the Macdonald polynomials for certain values of $t,$ 
central elements in quantum groups (Etingof and others). 

(b) The so-called interpolation  
polynomials (Macdonald, Lassalle, 
Knop\~Sahi, Okounkov\~Olshansky, Rains), continuing 
the classical Lagrange polynomials, appeared to 
be connected to double Hecke algebras. 
 
(c) The interpretation of the Macdonald 
polynomials as traces of the vertex operators, including 
applications to the Verlinde algebras 
(Etingof, Kirillov Jr.). It is mainly in the $GL$\~case 
and for integral $k.$ See also \cite{EVa}. 

(d) Various related results on KZB and its monodromy, 
(Etingof, Felder, Kirillov Jr., Var\-chen\-ko). 
The monodromy always satisfies the
DAHA group-type relations, but the quadratic ones 
are valid in special cases only. 
See, e.g., \cite{Ki1,TV,FTV,FV2}. 

(e) There are multiple relations of DAHA to the theory of 
the $W$\~invariant differential operators, including recent 
developments due to Wallach, Levasseur, Stafford, and Joseph. 
In the theory of DAHA, a counterpart is concerning 
the centralizer of the nonaffine Hecke subalgebra. 

(f) There are direct links to 
the theory of quantum groups defined by Drinfeld and Jimbo
and dynamical Yang\~Baxter equations 
through the classical and quantum $r$\~matrices;
we refer to \cite{DrQ} and the books \cite{L1,ES}
without going into detail. 
\medskip 

Strong connections with the 
affine Hecke algebra technique in 
the classical theory of $GL_N$ and $\BS_n$ must be 
noted. I mean \cite{C9,Na,Na2,NT,LNT,C12} and 
promising recent results towards Kazhdan\~Lusztig polynomials 
by Rouquier and others. 

The expectations are that Kazhdan\~Lusztig polynomials 
and the canonical (crystal) bases in quantum groups 
are important for the theory of double Hecke algebras, 
although I don't know a good definition of the 
"double" Kazhdan\~Lusztig polynomials. 

The coefficients of the  symmetric Macdonald polynomials 
in the stable $GL$\~case 
have interesting combinatorial properties 
(Macdonald, Stanley, Hanlon, Garsia, Haiman, and others). 
The most celebrated 
is the so-called $n!$\~conjecture recently proved by Haiman. 
See \cite{GH} and \cite{Ha}. 

Let me mention here that the Macdonald polynomials for the 
classical 
root systems appeared for the first time in a Kadell work. He 
also proved the Macdonald norm conjecture for the $BC$ systems 
\cite{Kad}. The constant term conjecture in the $GL_n$\~case 
was verified in \cite{BZ}; see \cite{C2} and \cite{Ma3} 
for further references. 
The first proof of the norm formula for the $GL$ is due to Macdonald. 
The elliptic counterparts 
of the Macdonald operators in the case of $GL_n$ 
are defined in \cite{R}. 

Quite a few constructions can be extended 
to arbitrary finite groups generated by complex reflections. 
For instance, the Dunkl operators and the 
KZ connection exist in this generality (Dunkl, Opdam, Malle). 
One can also try the affine and even the hyperbolic groups 
(Saito's root systems 
\cite{Sai}). 

\subsection{Classical origins} 
There are deep relations to the theory of special 
functions including the $q$\~functions and the classical 
Fourier analysis. We will not try to reconstruct 
systematically the history of the subject and review the 
connections. See, e.g., Section \ref{sec:DAHAVER}, 
and \cite{M4, O2, C6, C12, EG}. As for the classical 
Fourier analysis, we recommend the book \cite{Ed}, 
although it is not directly related to the topics of 
this paper. The theory of the Riemann zeta function 
is one of the major achievements 
of the classical theory of functions; 
Riemann was always referred to as the greatest master 
of the Fourier analysis. 

The papers \cite{AI}, \cite{AW} are 
a good introduction to the basic 
hypergeometric function and the Rogers polynomials; 
see also \cite{GR}. The book \cite{An} can be 
definitely recommended to those who want to 
understand the theory of $q$\~functions.
The paper \cite{Ba} remains the best on the analytic
theory of the $q$\~Gamma function and the multiple Gamma
functions. See also the \cite{Kos} 
(the Bessel functions via $SL(2,\R)$),
\cite{HO1,He} (spherical functions 
and generalizations), 
and \cite{DV} (about Harish-Chandra, 
the creator of the harmonic analysis 
on symmetric spaces). 
\smallskip 
 
The main object of the paper is the $q$\~Fourier transform 
introduced in \cite{C3} as a 
$q$\~generalization of the classical Hankel and 
Harish-Chandra transforms, and the Macdonald\~Matsumoto 
$p$\~adic transform. It has deep relations to the 
Macdonald orthogonal polynomials, combinatorics, 
the Gaussian sums, conformal field theory, 
and the Verlinde algebras. An important development
was the nonsymmetric theory; see 
\cite{Du, O2, M4, KnS, C1, C4}. 

\setcounter{equation}{0} 
\setcounter{figure}{0} 
\section{Soliton connection} \label{SEC:soliton}
A direct link to soliton theory is 
the {\em $r$\~matrix KZ equation}
introduced in \cite{Ch0}.
We will demonstrate that it is directly connected with
the $\tau$\~function, the object of major importance in
soliton theory. 

\subsection{The {\em r}--matrix KZ}
The {\em $r$\~matrix KZ equation} was defined in \cite{Ch0}
for the usual classical $r$\~matrices and their 
generalizations attached to the root systems.
The key observation in this paper was
that the classical Yang\~Baxter equation (YBE) 
for $r$ is {\em equivalent} 
(under minor technical constraints) 
to the cross-derivative integrability condition for 
the generalized Knizhnik\~Zamolodchikov (KZ) equation 
in terms of $r.$ 
An immediate application was a new broad class 
of KZ-type equations with the trigonometric and 
elliptic dependence on the arguments. 

Another application is the {\em $W$\~invariant KZ,} a
$W$\~equivariant system of the differential equations  
with values in the group algebra 
$\C W$ of the Weyl group. It was defined in \cite{Ch0} 
as an example of $r$\~matrix KZ, but, as a matter of fact,
it is directly based on the so-called ``reflection equations" 
and their generalizations from \cite{C2}, 
the YBE extended to arbitrary Weyl groups.
The paper \cite{C2} readily gives a 
justification of the cross-derivative conditions 
for this system, 
although a straightforward proof via the reduction 
to the rank two consideration is simple as well. 

Generally, the $r$\~matrix must have some symmetries and 
satisfy the rank two commutator relations, that can be
immediately verified (if they hold). 
The problem is with finding examples and the 
representation theory interpretation.
In the case of the $W$\~invariant KZ
from \cite{Ch0}, the 
cross-derivative integrability conditions 
and the $W$\~equivariance 
follow from the properties of the intertwining 
operators of the affine Hecke algebras; 
see \cite{Ch4,C13} and the references therein. 

The $W$\~invariant KZ
is now the main tool in the analytic theory of 
the configuration spaces associated with the 
root systems. It was extended to the 
affine KZ equation (AKZ) for the {\em affine} 
Weyl group in \cite{Ch2,Ch1}, 
and then to the groups generated by complex reflections 
by Dunkl, Opdam and Malle. A significant part 
of the theory of AKZ is the calculation 
of its monodromy representation. 
\smallskip 

Note that one can formally associate Hecke algebras with 
abstract Coxeter groups and similar groups. The most 
universal proof of the fact, if it holds, 
that the resulting Hecke 
algebra satisfies the Poincare\~Birkhoff\~Witt 
property, i.e., that this algebra 
is a flat deformation of the group algebra, goes via 
the consideration of a KZ-type connection and 
obtaining the Hecke algebra as its monodromy. 

For instance, the monodromy of the 
$W$\~invariant KZ  identifies the corresponding 
Hecke algebra with $\C W$ for generic values of the 
parameter $k$ in the KZ equation directly related to 
the parameter $t$ of the Hecke algebra. 
\smallskip 

\subsection{Classical {\em r}--matrices} 
The classical $r$\~matrices were 
derived from the quantum inverse method due 
to Faddeev, Sklyanin, and 
Takhtajan (see \cite{FT}). They allowed us 
to standardize the study of soliton equations, to 
come close to their classification, and to understand 
the genuine role of soliton theory in 
mathematics. Without going into detail, let me give 
the following general references \cite{KuSkl,BD,STS,ChFun} and 
mention the paper \cite{Sk}, which influenced my research 
in soliton theory a great deal. 

In the definition of the $r$\~matrix KZ, I combined 
the {\em $r$\~matrices} with the theory of the 
{\em $\tau$\~function} from 
\cite{DJKM,DJKM1} and some previous papers of the same 
authors. As for the $\tau$\~function, let me also mention 
\cite{Sat} and the exposition \cite{Ver}. 
The third important ingredient was the approach to 
{\em affine flag varieties} 
from the papers \cite{KP,PK}. 

Let $\mathfrak{g}$ be a simple Lie algebra and 
$r(\lambda)$ a function of 
$\lambda \in \Comp$ taking values in $\mathfrak{g} \tensor 
\mathfrak{g}.$ We assume that 
$r(\lambda) = t/\lambda + \tilde r(\lambda)$ for some 
analytic function 
$\tilde r$ in a neighborhood of $\lambda = 0,$ where 
$t = \sum_\alpha I_\alpha \tensor I_\alpha, \{I_\alpha\} \subset 
\mathfrak{g}$ is an 
orthonormal basis with respect to the Killing form $(\,,\,)_K$ 
on $\mathfrak{g}.$ The notation, standard in the theory of
$r$\~matrices, will be used: 
$$ 
{}^1 a = a \tensor 1 \tensor 1 \tensor \dots, 
{}^2 a = 1 \tensor a \tensor 1 \tensor \dots,\dots, 
{}^{ij} (a \tensor b) = {}^i a\, {}^j b, 
((a \tensor b) c)_K = a (b,c)_K 
$$ 
for $a,b,c$ being  from the universal enveloping algebra 
$U(\mathfrak{g})$ 
of $\mathfrak{g}$ (or from any of its quotient algebra). 
We also set
$\{A\stackrel{\tensor}{\hbox{\bf ,}} B\}$ if the 
products of the entries $ab$ in $A\tensor B$
are replaced everywhere by $\{a,b\}$.

A classical $r$\~matrix is a solution of the classical 
Yang\~Baxter equation 
\begin{align} 
&[\,{}^{12}r(\lambda), {}^{13}r(\lambda+\mu) + {}^{23}r(\mu)] 
+[\,{}^{13}r(\lambda+\mu), {}^{23}r(\mu)] 
=0.\label{0.0} 
\end{align} 

\begin{theorem} 
Assuming that ${}^{12}r(\lambda) + {}^{21}r(-\lambda)= 0,$ 
the following three assertions are equivalent to  (\ref{0.0}) and to 
each other: 

(a) The relation 
$$ 
\gathered 
\{M(\lambda) \stackrel{\tensor}{\hbox{\bf ,}} M(\mu)\}\ = \ 
[\,r(\mu - \lambda), {}^1 M(\lambda) + {}^2 M(\mu)] 
\endgathered 
$$ 
defines a Poisson bracket on the 
indeterminate coefficients of $M(\lambda)$ that is the 
generic element of the Lie algebra $\widetilde{\mathfrak{g}}$ 
of $\mathfrak{g}$\~valued 
meromorphic functions. 

(b) The functions in the form 
\begin{align} 
&M_r(\lambda) = \sum \Res ( 
{}^{12}r(\lambda - \mu) {}^2 M(\mu))_K d\mu), 
\label{0.1} 
\end{align} 
where the sum is over the set of poles of 
$M \in \widetilde{\mathfrak{g}},$ 
constitute a Lie subalgebra $\widetilde{\mathfrak{g}}_r \subset 
\widetilde{\mathfrak{g}}.$ 

(c) The equations 
\begin{align} 
&\der G / \der\lambda_i = \kap \sum_{j\ne i} {}^{ij} 
r(\lambda_i - \lambda_j) G 
\label{0.2} 
\end{align} 
for a function $G$ of $\lambda_1, \lambda_2, \dots$ taking values in 
$U(\mathfrak{g}) \tensor \dots \tensor U(\mathfrak{g})$ 
satisfy the cross-derivative integrability conditions. 
\end{theorem} 

The references are \cite{ChFun,Ch0}. 
It readily follows from (\ref{0.1}) that the formal 
differentiations 
\begin{align} 
&D_g F = (F g F^{-1})_0 F, \quad M_0 = M - M_r, \quad 
g \in \widetilde{\mathfrak{g}} 
\label{0.3} 
\end{align} 
satisfy the relations 
\begin{align} 
&[D_g, D_f] = D_{[g,f]}, \quad g,f \in \widetilde{\mathfrak{g}}. 
\label{0.4} 
\end{align} 
Here $F = \exp(M_0)$ for the generic element $M_0$ of the Lie 
subalgebra $\widetilde{\mathfrak{g}}_0 \subset 
\widetilde{\mathfrak{g}}$ 
of holomorphic functions; 
differentiations act in the algebra of functionals of the 
coefficients of $M_0.$ These 
$\{D_g, g \in \widetilde{\mathfrak{g}}\}$ 
are the generalized B\"acklund\~Darboux infinitesimal 
transformations, 
while formula (\ref{0.4}) and its group analog embrace many 
concrete results (from the 1850 to the present) on composing 
B\"acklund transforms for nonlinear differential equations; 
see \cite{ChTau,ChFun} and the references therein. 

Now we can introduce the $\tau$\~{\em function} 
\index{tau@$\tau$\~function}  as a solution of the 
system ($g \in \widetilde{\mathfrak{g}}, \kap \in \Comp$) 
\begin{align} 
&(D_g \tau) \tau^{-1} = \kap (\sum \Res(F^{-1} dF,g)_K). 
\label{0.5} 
\end{align} 
Formula (\ref{0.1}) induces the decomposition 
$$ 
\widehat{\mathfrak{g}} = \mathfrak{g} \oplus \Comp c, \quad 
[x + \zeta c, y + \xi c] = [x, y] + \sum \Res(dx, y)_K c. 
$$ 
In this decomposition, $c$ acts as a multiplication by 
$\kap\;$ and $\widetilde{\mathfrak{g}}_0 
\subset \widehat{\mathfrak{g}}$ 
annihilates $\tau.$ If 
$\kap \in k {\BN}$ for an appropriate integer $k > 0$ (depending on 
the $r$\~matrix under consideration; see \cite{ChFun}), 
then the $\tau$\~function is ``integrable," 
for instance, satisfies infinitely 
many Hirota-type relations. 

All these are very close to the theory of modular forms and 
representations of {\em ad\`ele groups} in arithmetic, with 
$\widetilde{\mathfrak{g}}_r \subset \widetilde{\mathfrak{g}} 
\supset \widetilde{\mathfrak{g}}_0$ playing the 
role of principal and, 
respectively, integer ad\`eles, and $\tau$ being something like the 
Tamagawa measure. 
\smallskip 

\subsection{Tau function and coinvariant} 
System (\ref{0.2}) for $1 \leq i,j \leq n$ is a natural 
generalization of the Aomoto\~Kohno system with 
$r(\lambda) = t/\lambda,$ 
which, in its turn, is a direct generalization of the 
Knizhnik\~Zamolodchikov equation in the two-dimensional 
conformal field 
theory (with $\mathfrak{g} = {\mathfrak{gl}}_N$). I came to 
(\ref{0.2}) while performing the 
following calculation in \cite{Ch0}. 

Let us generalize the definition of $\tau$ to allow arbitrary 
representations of $\widehat{\mathfrak{g}}$ 
(not only the representations with the 
trivial action of $\widetilde{\mathfrak{g}}_0$ on $\tau$). 
We also allow 
the elements from 
$\widetilde{\mathfrak{g}}$ to have poles at 
(pairwise distinct) points 
$\lambda_1,\dots,\lambda_n.$ 

Given a finite dimensional 
representation $V$ of $\mathfrak{g},$ 
$\widetilde{\mathfrak{g}}_0$ naturally acts on 
$V^{\tensor n}$ by the projection onto the Lie algebra 
$\mathfrak{g} \times \dots \times \mathfrak{g}$ ($n$ times). 
Let us introduce the 
\index{Verma module} {\em BGG\~Verma module} 
$\sM$ as the universal 
$\widehat{\mathfrak{g}}$\~module generated by 
$V^{\tensor n},$ where $c$ acts by multiplication by $\kap\;$ 
and 
$\widetilde{\mathfrak{g}}_0$ acts on $V^{\tensor n}$ as above. 
For every $x \in \sM,$ 
there exists a unique element $\tau(x) \subset V^{\tensor n} 
\subset \sM$ 
such that $\tau(x) - x \in \widetilde{\mathfrak{g}}_r \sM.$ 
Here we use the decomposition 
$\widehat{\mathfrak{g}} = 
\Comp c \oplus \widetilde{\mathfrak{g}}_0 \oplus 
\widetilde{\mathfrak{g}}_r.$ 
\smallskip 

This construction (from \cite{ChTau}) 
gives an important interpretation of the 
$\tau$\~function as the Kac\~Moody {\em coinvariant}. 
In soliton theory, we examine its transformations 
with respect to Kac\~Moody commutative subalgebras 
(the $D_g$\~flows for pairwise commutative $g$). 
In conformal field theory (CFT), the tau function 
appears as the {\em $n$-point correlation function} 
in a way similar to the above calculation. The 
dependence on the positions of the points, 
the celebrated Knizhnik\~Zamolodchikov equation, is of key 
importance for CFT. As a matter of fact, the following 
reasoning is similar to the original deduction of KZ 
in \cite{KZ}. 
\smallskip 

We set formally 
\begin{align} 
&(l_i \tau)(x) = \tau(l_i x), \quad 
l_i = \sum_{\alpha,j \ge 0} (I_\alpha \tilde{\lambda}_i^{-1-j}) 
(I_\alpha \tilde{\lambda}_i^j),\notag 
\end{align} 
where $\tilde{\lambda}_i = \lambda - \lambda_i$ is a local 
parameter in a 
neighborhood of $\lambda_i\,$ and 
$l_i$ is the Sugawara element of degree $-1,$ 
namely, the generator $L_{-1}$ of the 
Virasoro algebra embedded into 
the $i$\~th component of $\widehat{\mathfrak{g}}.$ 
Then one can verify that 
$$ 
l_i \tau +\kap^\prime \der\tau/\der u_i = 
(\rho_i + \sum_{j \ne i}{}^{ij}r_V(u_i - u_j))\tau 
$$ 
for a proper $\kap^\prime$ and some $\rho$ (trivial 
in the most interesting cases),  where ${}^{ij}r_V$ is 
the image of $r$ in 
${\operatorname{End}}_\Comp(V^{\tensor n}).$ 
This leads to the compatibility condition for equations (\ref{0.2}) 
and eventually provides the integral formulas for the 
$\tau$\~function. 

This derivation of the $r$\~matrix KZ equation 
was used in \cite{Ch3} to simplify the algebraic part 
of the Schechtman\~Varchenko construction \cite{VSch} 
of the integral formulas for the rational KZ 
(\cite{DJMM} for $SL_2$) 
and generalize their formulas 
to the trigonometric case. 

\smallskip 
\setcounter{equation}{0} 
\section{DAHA in harmonic analysis}\label{sec:DAHAVER} 
\medskip 
This section is an attempt to outline the core 
of the DAHA theory in the special case of $A_1$ 
and discuss connections between 
harmonic analysis and mathematical physics. 
There are quite a few projects where DAHA is 
involved (see above). However, I think the Verlinde 
algebra is the most convincing demonstration of 
the power of new methods. It is the main object 
of this section. 
\smallskip 

\subsection {Unitary theories} \label{ORIGINS} 
Generally, the problem is that we do not have the 
Kac\~Moody harmonic analysis. 
For instance, we do not have a good definition 
of the category of $"\,L^2\,"$\~representations in this case. 
The theory of Kac\~Moody algebras, playing a well-known role in 
modern theoretical physics, is too algebraic 
for this, and, presumably, these algebras 
must be developed to more analytic objects. 
Actually, the von-Neumann algebras are such objects, 
but they are, in a sense, too analytic. We need something 
in between. Hopefully the Verlinde algebras and DAHA can help. 
\smallskip 

{\bf Physical connection.} 
Concerning the classical roots of the 
harmonic analysis on symmetric spaces, 
the corresponding representation theory was greatly 
stimulated by (a) physics, (b) the theory of special 
functions, and (c) combinatorics, historically, in the 
opposite order. In my opinion, the demand from physics 
played the major role. 

(a) Harish-Chandra was Dirac's assistant 
for some time 
and always expressed unreserved admiration for Dirac, 
according to Helgason's interesting recent note 
"Harish-Chandra." The Lorentz group led him to the 
theory of infinite dimensional representations of 
semisimple Lie groups. 
However, later the mathematical goals like 
the Plancherel formula became preponderant 
in his research. 

(b) The theory of special and spherical functions was the 
main motivation for Gelfand in his studies of 
infinite dimensional representations, 
although he always emphasized 
the role of physics (and physicists). 
It was reflected 
in his program (1950s) aimed at ``adding" the 
spherical and hypergeometric functions to the Lie theory. 

(c) Before the Lie theory, the symmetric group was the main 
``representation" tool in the theory of functions. 
It still remains of fundamental importance. 
However, using the symmetric group 
only is not sufficient to introduce and understand properly 
the differential equations and operators needed in the theory. 

\smallskip 
As for (a), not all representation theories are of 
physical importance. It is my understanding, that ``unitary" 
representations are of major importance, although 
modern theoretical physics uses all kinds of representations. 
For instance,``massive" quantum theories must have 
a positive inner product. 
Mathematically, the unitarity is need\-ed 
to decompose some natural spaces of functions, for instance, 
the spaces of square integrable functions, $L^p\,$\~spaces, 
and so on and so forth. 
\smallskip 

{\bf Main unitary lines.} 
There are four main sources 
of unitary theories: 

\noindent
(A) compact and finite groups, 

\noindent
(B) noncompact Lie groups, 

\noindent
(C) Heisenberg and Weyl algebras, 

\noindent
(D) operator and von-Neumann algebras. 

The Clifford algebra, super Lie algebras, and ``free fermions" 
are important too, but require a special discussion. 
The theory of automorphic forms and 
the corresponding representation theory, including that of 
the $p$\~adic groups and Hecke algebras, definitely
must be mentioned too. 
The arithmetic representation theory 
is an important part of modern mathematics 
attracting increasing attention in physics now. 

\smallskip 

The above unitary theories have merits and demerits. 
\smallskip 

The Heisenberg\~Weyl algebra 
is heavily used in physics (``bosonization"), 
but it has essentially 
a unique unitary irreducible representation, 
the Fock representation. 
The von-Neumann factors actually 
have the same demerit. Only the pair 
"factor\~subfactor" appeared good enough for 
combinatorially rich theory. Note that the 
operator algebras give 
another approach to the Verlinde algebras. 
\smallskip 

{\bf Spherical functions.} 
The theory (A) is plain and square, but only finite dimensional 
representations can appear in this way. 
The representation 
theory of noncompact Lie groups is infinite dimensional 
(which is needed in modern physics), 
but the Harish-Chandra transform is not self-dual and 
its analytic theory is far from being 
complete. Also (B) is not very fruitful from the viewpoint 
of applications in the theory of special functions. For instance, 
the spherical functions in the so-called group case ($k=1$) 
"algebraically" coincide with 
the characters of compact Lie groups. 
The other two values $k=1/2,2$ in the Figure \ref{diagram} 
below correspond, respectively, to the 
orthogonal case ($SL(n,\BR)/SO(n)$) and the symplectic case. 

The left column of the top block of Figure \ref{diagram} 
shows the classical theory of 
characters and spherical functions 
of compact and noncompact Lie groups 
extended towards the orthogonal polynomials, 
Jack\~Heckman\~Opdam polynomials. 
The latter are generally 
beyond the Lie theory. We can define them as orthogonal 
polynomials ($k$ must be assumed real positive or 
even ``small" negative), or as eigenfunctions of the 
Sutherland\~Heckman\~Opdam operators, generalizing the 
radial parts of the Laplace operators on symmetric 
spaces. The latter definition works for arbitrary 
complex $k$ (apart from a series of special values where the 
complete diagonalization is impossible). 

\smallskip 
\subsection{From Lie groups to DAHA} 
\begin{figure}
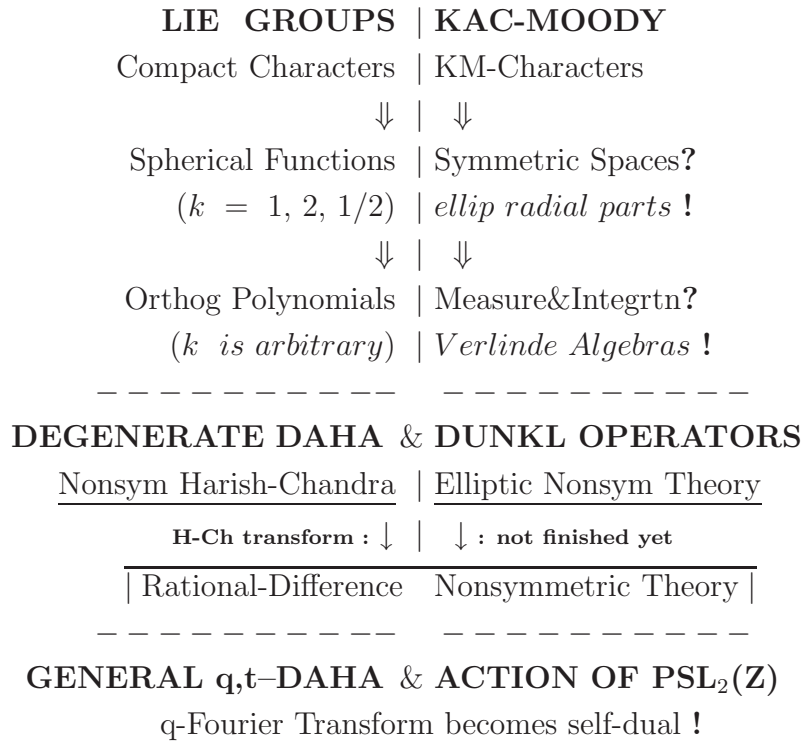
 
\begin{align} 
 \hbox{\bf LIE \ GROUPS} \  \mid\        & 
\hbox{\bf KAC-MOODY} \notag\\ 
 \hbox{Compact\ Characters} \ \mid\ 
&\hbox{KM-Characters}\notag\\ 
 \  \Downarrow \ \ \mid\  &\  \Downarrow\notag\\ 
 \hbox{Spherical\ Functions}  \ \mid\ 
&\hbox{Symmetric\ Spaces} 
\hbox{\bf ?}\notag\\ 
 (k\ =\ 1,\,2,\,1/2)   \ \mid\ &ellip\ radial\ parts\ 
\hbox{\bf !}\notag\\ 
 \  \Downarrow \ \ \mid\  &\  \Downarrow\notag\\ 
\hbox{Orthog\ Polynomials}  \ \mid\ 
&\hbox{Measure\&Integrtn} 
\hbox{\bf ?}\notag\\ 
 (k\ \ is \ arbitrary)  \ \mid\ &Verlinde\ Algebras\ 
\hbox{\bf !}\notag\\ 
      ----------  \ \, 
\ \ & ----------  \notag\\ 
\hbox{\bf DEGENERATE DAHA}\ \,\& \ & 
\hbox{\bf DUNKL OPERATORS}\notag\\ 
\underline{\hbox{Nonsym Harish-Chandra}} \ \mid\ & 
\underline{\hbox{Elliptic Nonsym Theory}}\notag\\ 
\hbox{\eightbf H-Ch transform :}  \downarrow \ \ \mid\ 
&\  \downarrow\hbox{\eightbf : not finished yet}\notag\\ 
\overline{\mid\hbox{Rational-Difference\ \ \ }} 
 & \overline{\hbox{Nonsymmetric Theory}\mid} 
\notag\\ 
      ----------  \ \, 
\ \ & ----------  \notag\\ 
\hbox{\bf GENERAL q,t--DAHA}\ \,\& \ & 
\hbox{\bf ACTION OF PSL}_2\hbox{\bf (Z)} 
\notag\\ 
{\hbox{\rm q-Fourier Transform}} 
&{\hbox{ becomes self-dual \bf !}}\notag 
\notag 
\end{align} 
\vskip -0.2cm 
\caption{Harmonic Analysis and DAHA} 
\label{diagram} 
\end{figure}

Considering the left column of Figure \ref{diagram} as a 
sample harmonic analysis program, 
not much is 
known in the Kac\~Moody case (the right column). 

At the level of ``Compact Characters" in the classical 
harmonic analysis, 
we have the theory of Kac\~Moody 
characters which is reasonably complete, in spite of 
combinatorial difficulties with 
the so-called string functions. The next level is 
supposed to be the theory of KM spherical functions. 
There are several approaches (let me mention Dale Peterson, 
Lian and  Zuckerman, and the book 
of Etingof, I.Frenkel, Kirillov Jr.), but we have 
no satisfactory 
general theory so far 
with a reservation about the group case, where the technique 
of vertex operators and conformal blocks can be used. 
Extending the theory from the group case  ($k=1$) to 
arbitrary root multiplicities is a problem. 

However, in my opinion, the key problem is that 
the Kac\~Moody characters 
are not pairwise orthogonal functions in a way 
that one can expect taking the classical theory 
(the left column) as a sample. 
Analytically, they are given in terms of theta functions and 
cannot be integrated over noncompact regions in any direction 
unless special algebraic tools are used. It is exactly how 
Verlinde algebras enter the game. 
\smallskip 

{\bf Verlinde algebras.} 
We can identify the characters of integrable representations 
of the Kac\~Moody algebra of a given level (central charge) $c$ 
with the corresponding classical characters treated as 
functions at a certain alcove of $P/NP$ for the weight lattice $P$ 
and $N=c+h^\vee$; $h^\vee$ is the dual Coxeter number. 
The $P/NP$ is naturally a set of vectors with the components
in the $N$\~the roots of unity.
The images of the characters 
form a linear basis in the algebra of all functions on this 
alcove, called the Verlinde algebra. 
The fusion product corresponds to the pointwise multiplication, 
and the images of the characters become 
pairwise orthogonal with respect to the Verlinde inner product. 

The identification of the KM characters with the classical 
characters at roots of unity was used by Kac for the first time, 
when he calculated the action of the $SL(2,\BZ)$ 
on the KM-characters. 
The automorphism $\tau\mapsto -1/\tau$ 
transforms the images of the characters 
into the delta functions of the points of the alcove. 

We mention that the interpretation of 
the Verlinde algebras in terms of quantum groups due 
to Kazhdan\~Lusztig\~Finkelberg leads to 
the classical characters at roots of unity as well. 

The general drawback of this approach is that the $c$, 
the levels, 
must be positive integers. The corresponding Verlinde 
algebras are totally disconnected for different levels, 
unless a special $p$\~adic limiting procedure is used, 
similar to the one due 
to Ohtshuki in the theory of invariants of knots and links. 
It is why a uniform theory for all (unimodular) $q$
is needed. 

Another drawback is of a technical nature. The Verlinde 
algebras are too combinatorial apart 
form the $\mathfrak{sl}_2$\~case.
For instance, there is a very difficult problem of finding
the $PSL(2,\Z)$\~invariants in the tensor product
of the Verlinde algebra and its dual.
There is a solution of this problem, important
in physical applications, only in the cases of  
$\mathfrak{sl}_2$ and $\mathfrak{sl}_3$. The list
of the ``invariants" is sophisticated in the
latter case. 
 
It is expected that 
the theory at generic $q$ would be simpler to deal with.
Many physicists now are working in the setting of 
unimodular $q$ ($|q|=1$) that are {\em not} roots of unity. 

{\bf Difference theory.} 
The bottom block of Figure \ref{diagram} shows the 
difference theory and the general double affine Hecke algebra 
(DAHA). In contrast to the previous ``Harish-Chandra level," 
the $q$\~Fourier transform is self-dual, 
as holds for the classical Fourier transform and 
in the group case. Thus we are back to ``normal" 
Fourier theory. 

This, however, does not mean that the analytic 
difficulties dissappear. 
{\em Six}\, different analytic 
settings are known in the $q,t$\~case. Namely, 
there are one compact and two 
noncompact theories, and each exists in two variants: 
for real $q$ and unimodular $q,$ not counting 
the choice of the analytic spaces that can be used for 
the direct and inverse transforms. 
\smallskip 

There are two other theories of an algebraic nature. 
The {\em seventh} setup is the theory at the roots of 
unity. All irreducible representations 
are finite dimensional as $q$ is a root of unity. 
The major example 
is the nonsymmetric Verlinde algebra, which is Fourier-invariant 
and, moreover, invariant with respect to the projective 
action of $PSL(2,\BZ).$ The {\em eighth} is 
the theory of $p$\~adic integration, 
where the connection with DAHA has not yet been established.

The projective action of $PSL(2,\BZ)$ 
is one of the most important parts of the theory 
of DAHA. The Fourier transform corresponds to the matrix 
{\small $\Bigl(\begin{array}{cc} 
\hbox{0} & \hbox{1}\\ 
\hbox{-1}& \hbox{0}\end{array}\Bigr)$ } and 
is related to the transposition of the 
periods of an elliptic curve. 
\smallskip 

{\bf Special functions.} 
There is a fundamental reason 
to expect much from  DAHA in analysis: it is 
a {\em direct} generalization of
the Heisenberg\~Weyl algebra. DAHA was designed to 
fulfill a gap between the representation 
theory and classical theory of the special 
functions. It naturally 
incorporates the Bessel function, the 
hypergeometric and basic (difference) 
hypergeometric functions, 
and their multidimensional generalizations into the 
representation theory. This is especially 
true for the basic hypergeometric 
function. 

DAHA is heavily involved in the {\em multi\-dimen\-sional} 
theory of hyper\-geo\-metric-type functions. 
This direction is relatively recent. The 
Dunkl operators and the rational DAHA are the main tools 
in the theory of the multidimensional Bessel functions.
This theory is relatively recent; the first reference seems
\cite{O3}. Respectively, 
the degenerate DAHA and the trigonometric Dunkl operators 
serve the multidimensional hypergeometric functions. 
\smallskip

The one-dimensional DAHA is closely related 
to the super Lie algebra $\mathfrak{osp}(2|1)$ and 
to its ``even part," $\mathfrak{sl}(2)$ via 
the ``exponential map" from DAHA to its rational 
degeneration.  The one-dimensional 
Dunkl operator is the square root of the radial 
part of the Laplace operator in the rank one case, 
similar to the Dirac operator, although it is 
a special feature of the one-dimensional case. 

For instance, the finite dimensional representations 
of $\mathfrak{osp}(2|1)$ have a natural 
action of DAHA of type $A_1$. 
Their even parts, the classical 
finite dimensional representations of 
$\mathfrak{sl}(2),$ are modules 
over the subalgebra of symmetric elements of DAHA: 

\smallskip 
\centerline{\em  DAHA can be viewed as 
a natural successor of $\mathfrak{sl}(2).$} 

\smallskip 
\subsection{Elliptic theory} 
A solid part of the Kac\~Moody harmonic analysis 
is the construction of the 
"elliptic radial parts" due to 
Olsha\-net\-sky\~Perelo\-mov (the $GL_N$\~case), 
Ochiai\~Oshima\~Sekiguchi 
(the $BC$\~type), and from \cite{Ch12} 
for all reduced root systems. 
They exist only at the {\em critical level}, 
which does not make them 
too promising. We certainly need the theory for an arbitrary 
central charge. The degenerate 
(``trigonometric," to be more exact) 
double affine Hecke algebra presumably provides such a theory. 
\smallskip 

{\bf Nonsymmetric theory.} 
The most recognized applications of the degenerate DAHA 
so far are in classical harmonic analysis, namely, this 
algebra is very helpful in the Harish-Chandra theory of 
spherical functions and, more importantly, 
adds the nonsymmetric spherical functions and polynomials 
to the Lie\~Harish-Chandra theory. 

The {\em same} degenerate DAHA also serves the 
elliptic nonsymmetric Harish-Chandra theory. 
The ``nonsymmetric" is the key here, because 
the ``elliptic" Dunkl operators (infinite trigonometric, to 
be more exact), commute at {\em arbitrary levels}, 
in contrast to the elliptic radial parts. The latter 
can be always defined as symmetrizations of the Dunkl operators, 
but commute only at the critical level \cite{Ch12}. 

We also mention the recent \cite{Rai} devoted to the 
interpolation and orthogonal 
``elliptic polynomials" in the $C^\vee C$\~case. 
The connection with the DAHA of type $C^\vee C$ 
(defined by Noumi\~Sahi) has not been established yet. 

\smallskip 
In the theory of DAHA, 
the ``transforms" become homomorphisms 
of representations. For instance, 
the nonsymmetric Harish-Chandra 
transform appears as an analytic 
homomorphism from the trigonometric-differential 
polynomial representation of the degenerate DAHA 
(``Nonsym Harish-Chandra") to its 
rational-difference polynomial representation. 

There exists a 
third elliptic-differential polynomial representation of DAHA 
(``Elliptic nonsym theory") 
which is analytically isomorphic to the 
rational-difference representation. 
The corresponding isomorphism can be called a 
nonsymmetric elliptic Harish-Chandra transform. 
See the middle block of Figure \ref{diagram}. 

The analytic part of the elliptic 
theory is not finished yet. The general
$q,t$\~DAHA governs the monodromy
of the eigenfunctions of the 
elliptic Dunkl eigenvalue problem, corresponding to 
the action of the Weyl group and the translations 
by the periods of an elliptic curve. Therefore, 
besides the ``elliptic 
transform," shown in the figure as the arrow 
to the ``difference-rational 
theory," there is another arrow from 
the elliptic theory to the $q,t$\~block. 

Similar to the Verlinde algebras, that describe
the action of  $PSL(2,\BZ)$ on the Kac\~Moody characters,
the monodromy representation 
can be used to describe the projective action of the 
$PSL(2,\BZ)$ on the eigenfunctions of  
the elliptic Dunkl eigenvalue problem.
\smallskip 

{\bf Deformation.} 
A natural objective is to define the generalized 
Verlinde algebras 
for {\em arbitrary} unimodular $q.$ Since DAHA describes 
the monodromy of the elliptic nonsymmetric theory 
(the previous block), such an extension of the 
Verlinde theory is granted. 
However, what can be achieved via DAHA appears more 
surprising. 

\smallskip 
Under minor 
technical restrictions, DAHA gives a {\em flat} 
deformation of the Verlinde algebra to arbitrary 
unimodular $q.$  To be exact, {\em little} Verlinde 
algebras must be considered here, associated with the 
root lattices instead of the weight lattices. 

The dimension of the little Verlinde algebra remains the same 
under this deformation. 
All properties are preserved 
but the integrality and positivity 
of the structural constants, which are lost 
for generic $q.$ 

These deformations were constructed 
in \cite{C7,C12} using {\em singular} $k,$ 
rational numbers with the Coxeter number as the 
denominator. See Section \ref{sec:FLAT} 
below.

It triggers the question 
about the classification of all representations of 
DAHA for such special fractional 
$k.$ Indeed, if a ``unitary" 
Kac\~Moody representation theory exists, it could 
be associated with such representations, 
by analogy with the 
interpretation of the Verlinde algebras via the integrable 
representations of the Kac\~Moody algebras. 
We also note a relation of the finite dimensional
representations of DAHAs (for singular $k$) 
to the generalized Dedekind\~Macdonald
$\eta$\~identities from \cite{M6}; see \cite{C12}.

There is an important application. We establish a 
connection between little 
Verlinde algebras of type $A_1$ and 
finite dimensional representations of $\mathfrak{sl}(2),$ 
since the deformation above followed by the 
limit $q\to 1$ to the rational DAHA 
equip the little Verlinde algebras with the action of 
$\mathfrak{sl}(2).$ 

\smallskip 
{\bf Toward KM harmonic analysis.} 
We conclude, that generalizations of the Verlinde algebra
are reasonable candidates for categories of representations
of unitary type, especially for those
associated with the Kac\~Moody 
algebras. If we expect the reduced category of
integrable Kac\~Moody representations 
to be a part of such a category, then it must
be an extension of the Verlinde algebra, or, more likely,
must have the reduced category as a quotient.
I think that it is unlikely that 
the $L^2$\~KM theory, if it exists, can be assumed to
contain ``\,$1$\,", the basic representation, and other
integrable representations from \cite{K}. 

Using DAHA for this project, its unitary representations 
at $|q|=1$ are natural candidates for such categories, 
especially for singular $k,$ although 
generic $q,k$ are interesting too. 
They are expected to have some extra structures. The most
important are a) ``fusion multiplication", 
b) hopefully, the action of $PSL(2,\BZ)$, and 
c)``good restrictions" to the affine Hecke subalgebras
similar to the classical decomposition with respect to the maximal 
compact subgroup. 

The existence of multiplication indicates that 
the quotients (if we include $1$)
and constituents (generally) of the polynomial representation 
of DAHA or those for the  functional 
variants of the polynomial representation 
must be examined first. All of them are commutative 
algebras. The most interesting ``AHA\~restriction theory" 
occurs at the singular $k,$ when the polynomial representation 
becomes reducible. 
\smallskip 

\setcounter{equation}{0} 
\section{DAHA and Verlinde algebras} \label{DAHAVer}  
The Lie groups formalize the concept of symmetry 
in the theory of special functions, combinatorics, 
geometry, and, last but not  least, 
physics.  In a similar way, {\em abstract} Verlinde algebras 
``describe" Fourier transforms, 
especially the theories with the Gaussian, satisfying
the following fundamental property.
The Fourier image of the Gaussian must be 
to be proportional to its inverse or the
Gaussian itself, depending on the
setting, similar to the properties
of the Laplace and Fourier
transforms in the classical theory. 
Thus, the Verlinde algebras (are expected to)
formalize an important portion of the 
classical Fourier analysis.

Note that the Fourier transform has an interpretation
in the Lie theory as a reflection from the Weyl group.
There many such reflections apart from the
$\mathfrak{sl}(2)$\~case. This does not match well
the natural expectation that there must be a unique
Fourier transform, or, at least, the major one if there
are several. It is really unique in the DAHA theory.

\smallskip 
\subsection {Abstract Verlinde algebras} 
In the finite dimensional semisimple variant, 
the {\em abstract Verlinde algebra} is 
the algebra of $\BC$\~valued functions $V=$Funct$(\,\bowtie\,)$ 
on a finite set $\,\bowtie\,$ equipped with a 
linear automorphism $\si,$ {\em the Fourier transform}. 
The algebra $V$ has a unit, which is $1$ considered as 
a constant function. 
Note that $\si$ 
is not supposed to preserve the (pointwise) multiplication. 
As a matter of fact, it never does. 

The space $V$ has a natural basis of the 
{\em characteristic functions} $\chi_i(j)=\de_{ij},$ where 
$i,j\in\,\bowtie\,,\ \de_{ij}$ is the Kronecker delta. 

The first two assumptions are that 
\smallskip 

\noindent 
{\bf (a)} $\si^{-1}(1)=\chi_o$ for the 
{\em zero-point} 
$o\in\,\bowtie\,,\,$ and 

\noindent 
{\bf (b)} the numbers $\mu_i\equal\si(\chi_i)(o)$ are 
nonzero. 

Since the latter constants are interpreted as masses of particles 
in conformal field theory, let us impose further 
conditions  $\mu_i>0.$ 
\smallskip 

Then the {\em spherical functions} are 
$p_i\equal\si(\chi_i)\mu_i^{-1}$ for $i\in\,\bowtie\,.$ 
In other words, $p_i$ is proportional 
to $\si(\chi_i)$ and satisfies the normalizing 
condition $p_i(o)=1.$ 

Introducing the {\em inner product} as 
$\langle f, g\rangle\equal\sum_i \mu_i f(i)\overline{g(i)},$ 
the corresponding delta functions are $\de_i= 
\mu_i^{-1}\chi_i=\si^{-1}(p_i).$ Indeed, they 
are obviously dual to the 
characteristic functions with respect to the inner product. 

Concerning taking the weights 
$\mu_i=\si(\chi_i)(o)$ in the inner product, such choice 
is equivalent 
to the following classical property, which holds in all 
variants of the Fourier theory: 

{\em the Fourier-images of the 
delta functions are spherical functions.} 
\smallskip 

The next assumption is that 

\noindent 
{\bf (c)} $\si$ 
is unitary up to proportionality with respect to 
$\langle\, ,\, \rangle.$ 
\smallskip 

It readily gives the 
{\em norm formula:} 
$$\langle p_i,p_j\rangle=\de_{ij}\mu_i^{-1}\lan 1,1\ran. 
$$ 

In this approach, the latter formula is a result of a simple 
sequence of formal definitions. However, 
it is really fruitful. 
It leads to the best-known justification of the 
norm formula for the Macdonald polynomials, including 
the celebrated constant term conjecture, which is the 
formula for $\langle 1,1\rangle.$ 

To be more exact, it provides a deduction of the 
norm formula from the so-called evaluation formula 
for the values of the Macdonald polynomials at 
the ``zero-point."  In its turn, 
the latter formula results 
from the self-duality of the double affine 
Hecke algebra (the existence of the self-dual Fourier 
transform). 

The self-duality and 
the evaluation formula are directly connected with the 
following important symmetry: 

\smallskip 
\noindent 
{\bf (d)} $p_i(j)=p_j(i) \for i\,\in\,\bowtie\,\ni\, j.$ 
\smallskip 

Note that the theory of Macdonald polynomials 
requires an infinite dimensional variant of 
the definitions considered above. 
Namely, $\si^{-1}$ becomes an 
isomorphism from the algebra of Laurent 
polynomials to its dual, the space of 
the corresponding delta functions. Hence it 
becomes a map from one algebra to another algebra. 

Nevertheless it is possible to deduce the 
major formulas for the Macdonald polynomials 
within the finite dimensional ``self-dual" setting above, 
by using the consideration at the roots of unity. 
Without going into detail, it goes as follows. 

The number of Macdonald polynomials, which 
are well defined when $q$ is a {\em root of unity}, 
grows together with the order 
of $q.$ One checks, say, the Pieri rules 
for such Macdonald polynomials using the duality
argument, then tends the order 
of $q$ to $\infty,$ and, finally, 
obtain the desired formula 
for arbitrary $q$ and {\em all} Macdonald 
polynomials. The Pieri rules
for multiplication of the Macdonald polynomials by
the monomial symmetric functions were, indeed,
justified in  \cite{C3,C4} by using the roots if unity.
\smallskip 

Concluding the discussion of abstract Verlinde algebras,
the {\em Gaussian} appears as a function 
$\ga\in V$ such that 

\smallskip 
\noindent 
{\bf (e)} 
$\si(\ga)=$const$\cdot\ga^{-1},\ 
\ga(i)\neq 0\for i\in \,\bowtie\,.$ 
\smallskip 

The constant here is the {\em abstract Gauss\~Selberg sum}. 
It is necessary to fix the 
normalization of $\si$ and $\ga$ 
to make this sum well defined. The normalization of 
$\si$ has already been fixed by the condition $\si(\chi_o)=1.$ 
The natural normalization of the Gaussian is $\ga(o)=1.$ 
\smallskip 

The assumptions (a-e) are more than sufficient to make the 
Verlinde algebras rigid enough. However, 
in my opinion, the key axiom is PBW, the 
Poincare\~Birkhoff\~Witt property, 
which requires the operator approach, to be discussed next. 

\smallskip 
\subsection{Operator Verlinde algebras} 
The above discussion is actually about a ``good self-dual" 
invertible linear operator acting in a commutative algebra. 
It is obviously too general, and more structures have to be added. 
We need to go to the operator level, switching from the 
characteristic functions and the spherical functions 
to the corresponding {\em commutative} algebras $\x,\y$ of 
the operators, which are diagonal at the corresponding sets of 
functions. 

The {\em operator Verlinde algebra} 
$\a$ (the main example will 
be the double affine Hecke algebra) is generated by 
commutative algebras $\x$ and $\y,$ and the algebra $\H$ 
controlling the symmetries of the $X$\~operators and 
the $Y$\~operators. In the main examples, 
the Weyl groups (or somewhat more 
general groups) are the groups of symmetries upon 
degenerations; $\H$ are the corresponding Hecke algebras 
serving the non-degenerate $q,t$\~case. 

The key and the most restrictive 
assumption is the {\em PBW property}, which states that 

\smallskip 
\noindent 
{\bf (A)} the natural map from the tensor product 
$\x\otimes\y\otimes \H$ to $\a$ is an isomorphism of 
the linear spaces, as well as 
the other five maps corresponding to the 
other orderings of $\x,\y,\H.$ 
\smallskip 

The Fourier transform and the Gaussian are formalized 
as follows. The 
{\em projective $PSL_2(\Z)$} must act in $\a$ by 
outer automorphisms, i.e., 

\smallskip 
\noindent 
{\bf (B)} $\tau_+,\tau_-$ 
act in $\a$ as algebra automorphisms
and satisfy the Steinberg relation 
$\tau_+\tau_-^{-1}\tau_+$ $=\tau_-^{-1}\tau_+\tau_-^{-1}$; 

\noindent 
{\bf (C)} 
they preserve the elements from $\H,\,$ the element
$\si\equal \tau_+\tau_-^{-1}\tau_+$ maps $\x$ onto $\y,\,$ and 
$\tau_+$ is the identity on $\x.$ 
\smallskip 

Note that  (B)\~(C) provide that 
$\tau_-=\si \tau_+^{-1} \si^{-1}$ 
is identical on $\y.$ The automorphism $\tau_+$ becomes 
the multiplication by the Gaussian in Verlinde algebras, 
which connects (C) with the assumption (e) above.
\smallskip 

The remaining feature of abstract Verlinde algebras 
to be interpreted using 
the operator approach 
is the existence of the inner product. 
We will postulate the existence of the corresponding 
involution. 
We assume that: 

\noindent 
{\bf (D)} $\a$ is a quotient of the group algebra of the 
group $\b$ such that the anti-involution 
$\b\ni g\mapsto g^{-1}$ of $\b$ becomes an anti-involution 
of $\a;$ 

\noindent 
{\bf (E)} $\tau_{\pm}$ and $\si$ come from automorphisms 
of the group $\b$ and therefore commute with the 
anti-involution from (D). 
\smallskip 

We note that the inner products can be associated
with somewhat different involutions, for instance,
in the case of {\em real} harmonic analysis ($q>0$). 
The assumption that there 
is a system of {\em unitary} generators 
is a special feature of the unimodular 
theory ($|q|=1$). 
The anti-involution in the Verlinde case 
and those in the main generalizations  
do satisfy (D)\~(E). 

\smallskip 
The Verlinde algebras can be now {\em re-defined} as 
$\si$\~invariant unitary irreducible representations $V$ 
of $\a$ that are $\x$\~spherical, i.e., are some  
quotients of the (commutative) algebra $\x.$ 
A representation is called unitary if it has a hermitian 
inner product inducing the anti-involution of $\a$ from (D). 
Note that (D)\~(E) 
{\em automatically} guarantee that 
the Fourier transform $\si$ is ``projectively unitary" 
in $\si$\~invariant unitary irreducible representations. 

Thanks to 
irreducibility, both $\si$ and $\tau_+$ 
are fixed uniquely in the group Aut$_{\BC}(V)/\BC^*$ 
(if they act in $V$) and induce the corresponding 
automorphisms of $\a.$ 
\smallskip 

Strictly speaking, assumptions (D)\~(E) can be replaced 
by a more general property that $\a$ has an anti-involution 
which commutes with the $\tau_{\pm}$ and $\si.$ 
However the double Hecke algebras (the main examples so far)
really appear as quotients of the 
group algebras and $\tau_{\pm}$ and $\si$ act in this group,
the elliptic braid group.

An important advantage
of the operator approach is that we can 
relax the constraints for the abstract Verlinde algebras 
defined above by 
considering {\em non-unitary} and even {\em non-semisimple} 
projective $PSL(2,\Z)$\~invariant $X$\~spherical 
irreducible representations of $\a.$ They do appear 
in applications. 

Another advantage is that it is not necessary to 
impose the Fourier-invariance. 
Generally, the main problem of the Fourier analysis 
is in calculating the Fourier 
images $\si(V)$ of {\em arbitrary} $\a$\~modules $V$ and the 
corresponding transforms 
$V\ni v$ $\mapsto \si(v)\in \si(V),$ which 
induce $\si$ in $\a.$ 

Similarly, $\tau_+$ becomes  multiplication by the Gaussian 
given by a variant of the classical formula 
$e^{x^2}$ in the main examples. The operator approach 
makes it possible to define the Gaussian for any irreducible 
$\a$\~module $V$; it as an operator with values in $\tau_+(V).$ 
\smallskip 

Let us discuss now the DAHA of type $A_1$ in detail, 
where the theory is already quite interesting. 
The transition to arbitrary root system is 
sufficiently smooth. 

\smallskip 
\subsection {Double Hecke Algebra} 
The most natural definition goes through the 
{\em elliptic braid group\,} 
$\mathcal{B}_q\equal \langle T,X,Y,q^{1/4}\rangle/$ 
with the relations 
\begin{align} 
&TXT=X^{-1},\ TY^{-1}T=Y,\notag\\ 
&Y^{-1}X^{-1}YXT^2\ =\ q^{-1/2}. 
\notag \end{align} 

\noindent 
Here 
$$\mathcal{B}_1= 
\pi_1^{\hbox{\small orb}}(\{E\setminus 0\}/ 
\BS_2) 
$$ 
$\cong\pi_1(\{ E\times E\setminus \hbox{diag}\}/\BS_2), 
$ 
where $E$ is an elliptic curve. 
Using the orbifold fundamental group here makes it possible 
to ``divide" by the symmetric group $\BS_2$ 
without removing the ramification completely, 
i.e., removing all four points of second order. 
Only one puncture is needed to obtain the above relations.

Actually, 
in the case of $A_n,$  
it is sufficient to consider 
the product of $n+1$ copies of $E$ and remove the 
``diagonal" before dividing by $\BS_{n+1}$ 
instead of using the orbifold group. 
The corresponding braid group is {\em isogenous} to the 
one with the relations above; it was calculated 
by Birman and Scott. 

To complete the definition of DAHA,
we impose the quadratic $T$\~relation: 

$$\HH\equal\C[\mathcal{B}_q]/((T-t^{1/2})(T+t^{-1/2})).$$ 

If $t=1,$ $\HH$ becomes the 
 Weyl algebra extended by the reflection;
$T$ becomes $s$ satisfying $sXs=X^{-1},\, sYs=Y^{-1},\, s^2=1$.

\smallskip 

The {\em Fourier transform,} 
which plays a major role in the theory, 
is the following outer automorphism of DAHA: 

$$\sigma: X\mapsto Y^{-1}, Y\mapsto XT^2,\ 
T\mapsto T.$$ 
Thus the DAHA Fourier transforms finds a 
conceptual interpretation as the transposition 
of the periods of the elliptic curve, This 
is not surprising from the viewpoint of 
CFT, KZB, and the Verlinde algebras, but 
such a connection with topology still 
remains challenging when we deal with the 
applications of DAHA in harmonic analysis. 

The representations where $\si$ acts (i.e., 
becomes inner) are 
called {\em Fourier-invariant} or 
self-dual. The nonsymmetric Verlinde algebra and the 
Schwartz space are examples. 


More generally, the topological interpretation above 
readily gives that 
the group $PSL(2,\Z)$ acts projectively 
in $\HH$: 
\begin{align}\label{taupsi} 
&\tau_+: Y\mapsto q^{-1/4}XY,\ X\mapsto X,\ 
T\mapsto T,\\ 
&\tau_-:  X\mapsto q^{1/4}YX,\ \, Y\mapsto Y,\ T\mapsto T,\notag\\ 
&\si=\tau_+\tau_-^{-1}\tau_+=\tau_-^{-1}\tau_+\tau_-^{-1}, 
\notag \end{align} 
$$\hbox{where\ \ }\binom{11}{ 01} \mapsto \tau_+,\ \ 
\binom{10}{  11} \mapsto \tau_-.$$

We will use $k$ such that $t=q^k$ and 
(sometimes) set $\HH(k)$ as $t=q^k.$ 
The algebra $\HH(k)$ acts in 
$\cP=$ {\em Laurent polynomials}  in terms of $X=q^x,$ 
namely, 
\begin{align} 
&T\mapsto t^{1/2}s+ 
\frac{t^{1/2}-t^{-1/2}}{ q^{2x}-1}(s-1),\notag\\ 
&Y\mapsto spT,\ sf(x)\equal f(-x),\notag\\ 
&pf(x)\equal f(x+1/2), \ t=q^k. 
\notag \end{align} 

\noindent 
The operator $Y$ is called the {\em difference Dunkl operator.} 

It is important to note that 
$\tau_-$ preserves $\cP.$ On the other hand, 
$\tau_+$ does not act there for a very 
simple reason. If it acts then it must be 
multiplication by the Gaussian $q^{x^2},$ 
which does not belong to $\cP.$ 

The ``radial part" appears as follows. 
One checks that 
$Y+Y^{-1}$ preserves $\cP_{sym}\equal$ 
symmetric (even) Laurent polynomials 
(the Bernstein lemma in the theory of 
affine Hecke algebras). The restriction 
$H=Y+Y^{-1}\mid_{sym}$ is the $q$\~{\em radial part} 
and can be readily calculated: 

$$H=\frac{t^{1/2}X-t^{-1/2}X^{-1}}{ X-X^{-1}} p + 
\frac{t^{1/2}X^{-1}-t^{-1/2}X}{  X^{-1}-X} p^{-1}.$$ 

\subsection {Nonsymmetric Verlinde algebras} 
The key new development in the theory of orthogonal 
polynomials, algebraic combinatorics, and related harmonic 
analysis is the definition of the 
nonsymmetric Opdam\~Macdonald polynomials. 
The main references are \cite{O2,M4,C4}. Opdam mentions 
in \cite{O2} that this definition 
(in the differential setup) was given in Heckman's unpublished 
lectures. 

These polynomials are expected to be, generally speaking, 
beyond quantum groups and Kac\~Moody algebras because of 
the following metamathematical reason. 
Major special functions in the Lie and Kac\~Moody 
theory, including the classical and Kac\~Moody characters, 
spherical functions, and conformal blocks, 
are $W$\~invariant. 

However, this definition is not quite new in 
representation theory. 
The limits of the nonsymmetric 
polynomials as $q\to \infty$ are 
well known. They are the spherical functions due to Matsumoto. 
The nonsymmetric Verlinde algebras 
are {\em directly} connected to
the nonsymmetric polynomials evaluated at the roots of unity,
\smallskip

Let $q=\exp(2\pi i/N),\ 0<k< N/2,\ k\in\Z.$ 
The {\em nonsymmetric Verlinde algebra} $V$ 
is defined as the algebra of functions of the set 
\begin{align}\label{introbowtie}
&\bowtie\,=\,\{-\frac{N-k+1}{ 2},...,-\frac{k+1}{ 2}, 
-\frac{k}{ 2},\, 
\frac{k+1}{ 2},...,\frac{N-k}{ 2}\}.
\end{align} 
It has the unique structure of a $\HH\,$\~module 
that makes the map $q^x(z)=q^{z}$ 
from $\cP$ to $V$ a $\HH$\~homomorphism. 
The above set is not $s$\~invariant; it is $sp$\~invariant. 
Nevertheless, the formula for $T$ 
can be used in $V,$ because the contributions 
of the ``forbidden" points $k/2=s(-k/2)$ and 
$(N-k+1)/2=s(-(N-k+1)/2)$ come 
with zero coefficients. 

The operators 
$X,Y,T$ are unitary in $V$ with respect to 
the positive hermitian form which will not be 
discussed here (it generalizes the inner product 
for conformal blocks). 
The positivity requires 
choosing the ``minimal" primitive $N$\~root of unity $q$ above. 

\smallskip 
The whole $PSL(2,\Z)$ acts in $V$ projectively 
as well as in the image $V_{sym}$ of $\cP_{sym}.$ 
The latter image can of course be 
defined without any reference to 
the polynomial representation. 
A general definition is as follows: 
$V_{sym}=\{f\mid Tf=t^{1/2}f\}.$ 
Here it simply means that the function $f(z)$ must be 
$s$\~invariant (even) for the points $z$ that 
do not leave the above set under the action of $s$ 
(recall that it is not $s$\~invariant). 
Therefore: 

$$\dim V=2N-4k,\ \dim V_{sym}=N-2k+1.$$ 
\smallskip 

We call $V$ a {\em perfect representation.} 
By perfect, we mean that it has all 
major features of the irreducible representations 
of the Weyl algebras at roots of unity, i.e., cannot be better. 
Formally, it means that it has a perfect duality 
pairing, directly realted to the Fourier transform,
and a projective action of $PSL(2,\BZ).$ 

The nonsymmetric {\em characters} in $V$ 
and the symmetric ones in $V_{sym}$ 
are, respectively, the eigenfunctions of $Y$  
and  $Y+Y^{-1}.$ 
\smallskip 

When $k=0,$ we come to the well-known definitions 
in the theory of Weyl algebra. As $k=1,$ 
$V_{sym}^{k=1}$ is the usual {\em Verlinde algebra} 
\cite{Ve}: 

\centerline{\em $\tau_+$ becomes the Verlinde $T$\~operator 
and $\si$ becomes the $S$\~operator.} 

\smallskip 
\setcounter{equation}{0} 
\section{Topological interpretation} 
We use that the relations of an $\HH\,$ 
are mainly of group nature and introduce the 
\index{braid group (elliptic)} 
{\em elliptic braid group} 
$$ 
\cB_q=\lr T, X, Y, q^{1/4}\rr / \left \lr 
\begin{array}{c}TXT=X^{-1}, 
TY^{-1}T=Y,\\ Y^{-1}X^{-1}YXT^2q^{1/2}=1 
\end{array}\right \rr. 
$$ 
Now $T,X,Y,q^{1/4}$ are treated as group generators and 
$q^{1/4}$ is assumed central. 
The double affine Hecke algebra $\HH\,$ is the 
quotient of the group algebra of $\cB_q$ by the quadratic 
Hecke relation. It is easy to see that the change of 
variables $q^{1/4}T\mapsto T$, $q^{-1/4}X\mapsto X$, 
$q^{1/4}Y\mapsto Y$ defines an isomorphism 
$\cB_q\cong \cB_1\times \BZ,$ where the generator of 
$\BZ$ is $q^{1/4}.$ 

In this section, we give 
a topological interpretation 
of the group 
$$\cB_1=\lr T, X, Y\rr / \lr TXT=X^{-1}, TY^{-1}T=Y, 
Y^{-1}X^{-1}YXT^2=1\rr.$$ 

\subsection{Orbifold fundamental group}
Let $E$ be an elliptic curve over $\BC$, 
i.e., $E=\BC/\Lambda$, where $\Lambda=\BZ+\BZ\imath.$ 
Topologically, the lattice can be arbitrary. 
Let $o\in E$ be the zero point, and 
$-1$ the automorphism $x\mapsto -x$ of $E$. 
We are going to calculate the fundamental group of 
the space $(E\setminus o)/\pm 1=\BP_{\BC}^1\setminus o.$ 
Since this space is contractible, its usual 
fundamental group is trivial. We can take the quotient after 
removing all (four) ramification points of $-1.$ However, 
it would enlarge the fundamental group dramatically. 
Thus we need to understand 
this space in a more refined way. 

\begin{figure}[htbp] 
\begin{center} 
\vskip 0.95in 
\includegraphics[scale=0.4]{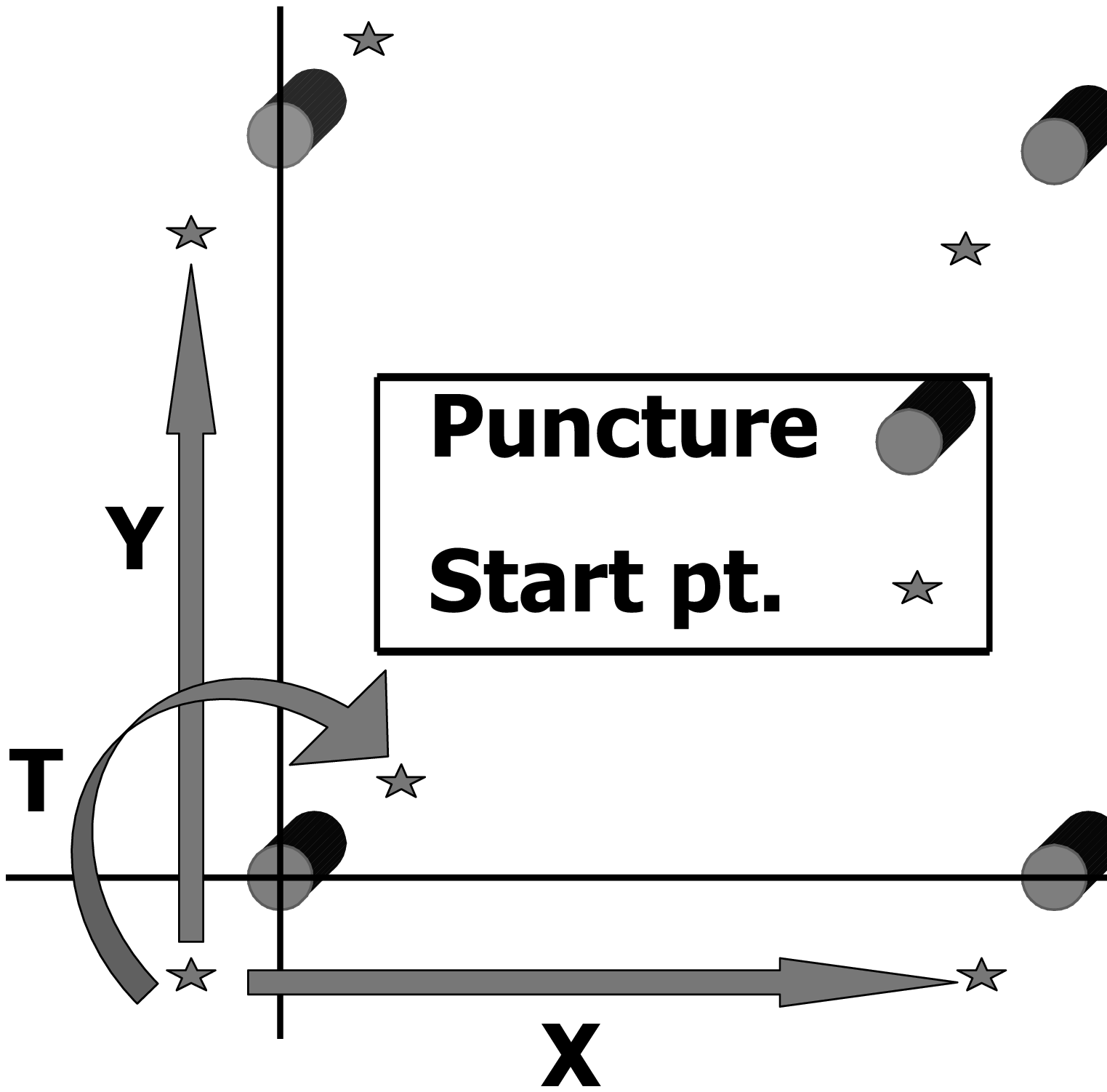} 
\vskip -1.25in 
\caption{Generators of $\cB_1$} 
\label{figxty1} 
\end{center} 
\end{figure} 

Let us fix the {\em 
base (starting) point} $\star=-\eps-\eps\imath\in \BC$ 
for small $\eps>0.$ 

\begin{proposition}\label{t82} 
We have an isomorphism $\cB_1\cong 
\pi_1^{orb}((E\setminus o)/\pm 1)$, where 
$\pi_1^{orb}(\cdot)$ is the orbifold 
fundamental group, which will be defined in 
the process of proving the proposition. 
\end{proposition} 
{\em Proof.} The projection map $E\setminus o\to 
(E\setminus o)/\pm 1=\BP_{\BC}^1\setminus o$ 
has three branching points, which 
come from the nonzero points of order 2 on $E.$ 
So by definition, 
$\pi_1^{orb}((E\setminus o)/\pm 1)$ is 
generated by three {\em involutions} $A, B, C,$ namely, 
the clockwise loops from $\star$ around the branching 
points in $\BP_{\BC}^1.$ 
There are no other relations. We 
claim that the assignment $A=XT$, $B=T^{-1}Y$, $C=XTY$ 
defines a homomorphism 
$$\pi_1^{orb}((E\setminus o)/\pm 1)\to \cB_1.$$ 
Indeed, $A$ and $B$ are obviously involutive. 
Concerning $C,$ the image of its square is 
$$XTYXTY=T^{-1}X^{-1}YXTY=T^{-1}YT^{-1}Y=1. 
$$ 
This homomorphism is an 
isomorphism. The inversion is given by the 
formulas  $ACB=T$, $ABCA=X$ 
and $AC=Y$. \sq 

This approach can hardly be generalized to arbitrary root 
systems; the following (equivalent) constructions can. 
The definition of the orbifold
fundamental group has to be modified as follows. 
We follow \cite{C15}. See also paper \cite{Io} (in this 
paper the orbifold group is not involved and the 
construction due to v.d.~Lek is used). 

We switch from $E$ to its universal cover $\BC$ 
and define the 
{\em paths} 
as curves $\gamma\in \BC\setminus \Lambda$ 
from $\star$ to $\widehat{w}(\star),$ where 
$\widehat{w}\in \widehat{W}=\{\pm 1\}\lsmash \Lambda.$ 
The generators $T,X,Y$ 
are shown as the arrows in Figure \ref{figxty1}. 

\begin{figure}[htbp] 
\begin{center} 
\vskip 0.95in 
\includegraphics[scale=0.4]{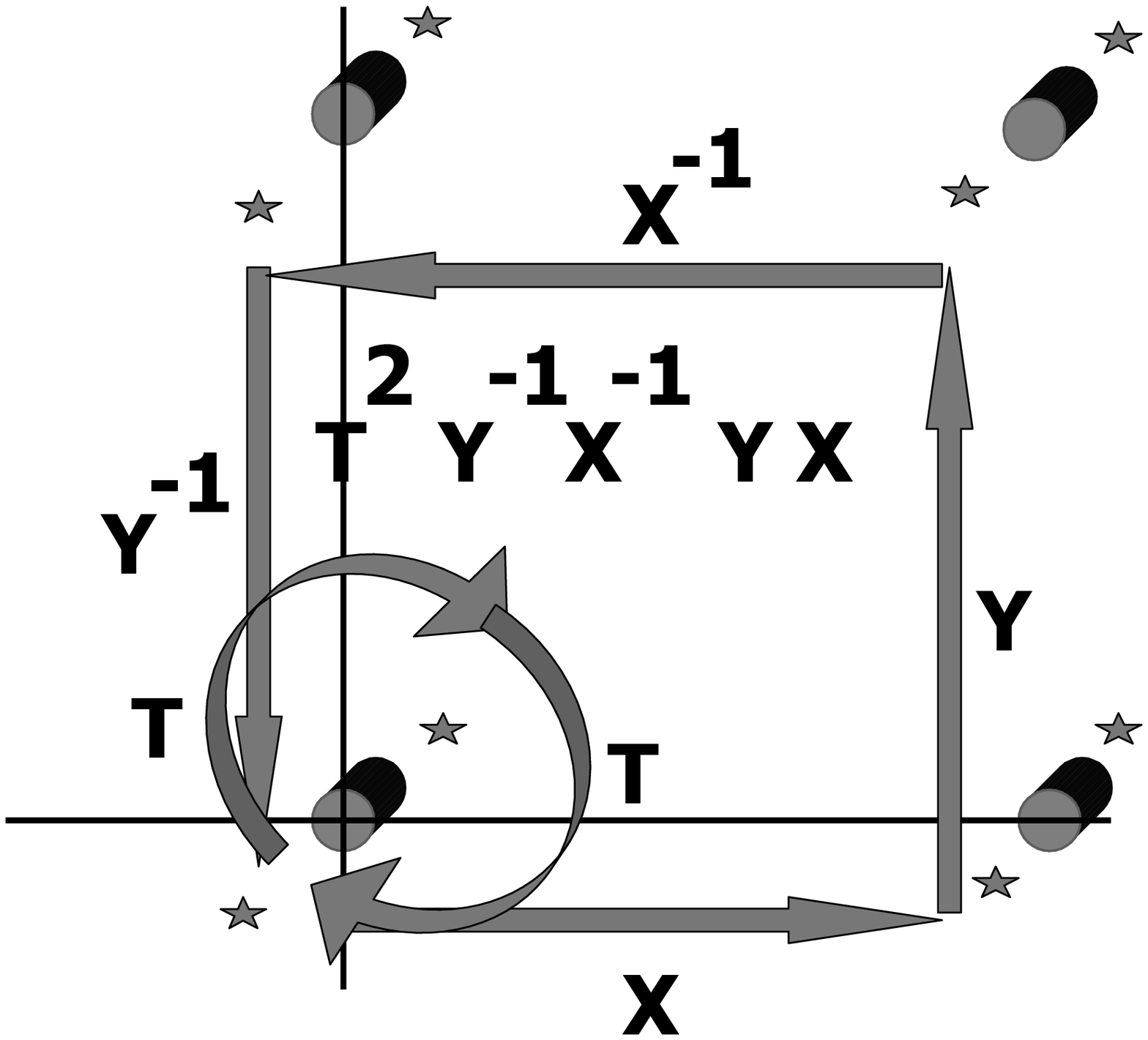} 
\vskip -1.25in 
\caption{Relation $T^2\,Y^{-1}\,X^{-1}\,Y\,X\,=\,1$} 
\label{figxty2} 
\end{center} 
\end{figure} 

The \index{composition of paths} 
{\em composition of the paths} is via $\widehat{W}$: 
we add the image of the second path under $\widehat{w}$ 
to the first path if the latter ends at $\widehat{w}(\star).$ 
The corresponding variant of Proposition \ref{t82} 
reads as follows. 
\smallskip 

\begin{proposition}\label{t82a} 
The fundamental group of the above paths modulo homotopy 
is isomorphic to $\cB_1,$ where 
$T$ is the half-turn, i.e., the clockwise half-circle from 
$\star$ to $s(\star),$ and $X,Y$ are $1$ and $\imath$ 
considered as vectors in $\R^2(=\C)$ 
originated at the base point $\star.$ 
\end{proposition} 
{\em Proof} is in Figures \ref{figxty1} and 
\ref{figxty2}. Figure \ref{figxty1} gives that 
when we first use $X$ then $T$ and 
then again $X$ (note that after $T$ the direction of $X$ 
will be opposit!), and then again $T,$ the loop 
corresponding 
to the product 
$TXTX$ will contain no punctures inside. Thus it equals $id$ 
in the fundamental group. The reasoning for 
$TY^{-1}TY^{-1}T$ is the same. 
Concerning the 
``commutator" relation, see Figure \ref{figxty2}. 
\sq 

Actually, this definition is close to the calculation 
of the fundamental group of $\{E\times E\setminus $\,diagonal$\}$ 
divided by the transposition of the components. See 
\cite{Bi}. 
However, there is no exact coincidence. 
Let me mention the relation to the elliptic braid group due to 
v.d.~Lek, although he removes all points of second order and his 
group is significantly larger. 
\smallskip

\subsection{The action of {\em PSL(2,Z)}} 
A topological interpretation is the best way 
to understand why the group $SL_2(\BZ)$ acts in $\cB_1$ 
projectively. 

Its elements act in $\BC$ natuarally, 
by the corresponding real linear transformations. 
On $E,$ they commute with the reflection $-1,$ 
preserve $o$, and permute the other three points of second order. 
Given $g\in SL_2(\BZ),$ we set $g=\exp(h), g_t=\exp(th)$ 
for the proper $h\in \mathfrak{sl}_2(\BR),$ $0\le t\le 1.$ 

The position of the base point $\star$ 
will become $g(\star),$ 
so we need to go back, i.e., 
connect the image with the base point by 
a path. To be more exact, the $g$\~image of 
$\gamma\in \pi_1^{orb}$ will be the union 
of the paths 
$$\{g_t(\star)\}\ \cup\ g(\gamma)\ \cup\ 
\{\hw(g_{-t}(\star)\}, 
$$ where the path for 
$\gamma$ goes from $\star$ to the point $\hw(\star).$ 

\begin{figure}[htbp] 
\begin{center} 
\vskip 1in 
\includegraphics[scale=0.4]{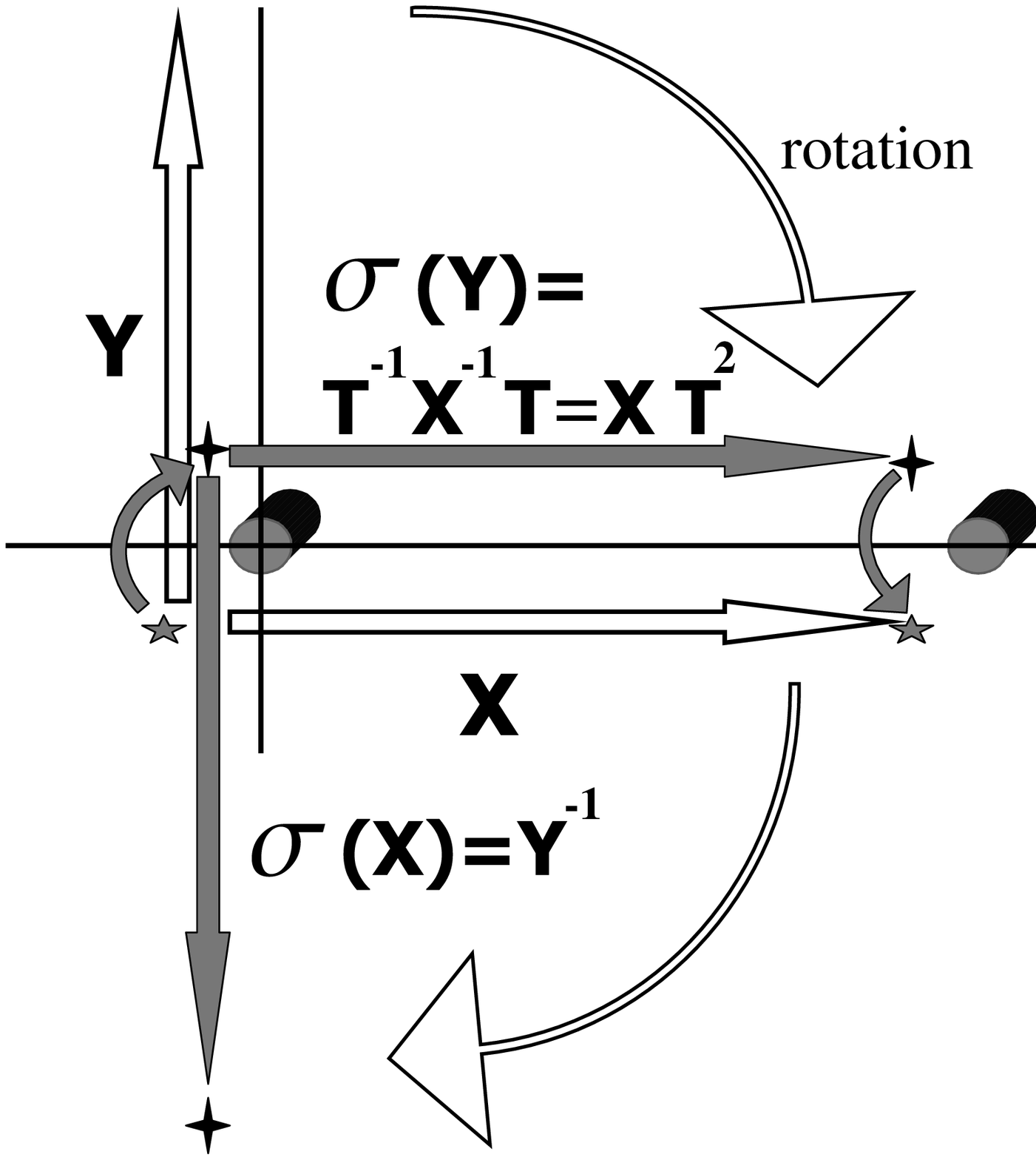} 
\vskip -1.25in 
\caption{Relations $\si(X)=Y^{-1},\ \si(Y)=XT^2$} 
\label{figxty3} 
\end{center} 
\end{figure} 

Figure \ref{figxty3} shows the action of the 
automorphism $\sigma$ corresponding to the 
rotation of the periods and $\star$ by $90^{\circ}$ 
with the origin taken as the center. 
Here  the dark straight arrows show the 
images of $X,Y$ (straight white arrows) with respect to 
this rotation. The quarter of a turn from the point 
$\star$ is the rotation path $\{g_t(\star)\}$ 
of this point as $0\le t\le 1.$ The other quarter of a turn 
is its $X$\~image with the opposit orientation. 
\smallskip 

We can always choose the base point sufficiently close 
to $0$ and connect it with its $g$\~image in 
a small neighborhood of zero. 
This makes the corresponding automorphism of $\cB_1$ 
unique up to powers of $T^2.$ All such automorphisms 
fix $T,$ because they preserve zero and the orientation. 

Thus we have constructed a homomorphism 
$$\alpha: SL_2(\BZ)\to 
\Aut_T(\cB_1)/T^{2\BZ}, 
$$ 
where $\Aut_T(\cB_1)$ is the 
group of automorphisms of $\cB_1$ 
fixing $T.$ 
The elements from $T^{2\BZ}=\{T^{2n}\}$ are 
identified with the corresponding inner automorphisms. 

Let  $\tau_+$, 
$\tau_-$ be the $\alpha$\~images 
of the matrices 
$\left(\begin{array}{cc}1&1\\0&1\end{array} 
\right)$ and $\left(\begin{array}{cc}1&0\\1&1 
\end{array}\right).$ 
Then $\sigma=\tau_+\tau_-^{-1}\tau_+$ 
corresponds to $\left(\begin{array}{cc}0&1\\-1&0 
\end{array}\right),$ 
and $\sigma^2$ has to be the conjugation by 
$T^{2l-1}$ for some $l.$ Similarly, 
$$ \tau_+\tau_-^{-1}\tau_+=T^{2m} \tau_-^{-1}\tau_+\tau_-^{-1}.$$ 
Using the rescaling $\tau_{\pm}\mapsto T^{2m_{\pm}}\tau_\pm$ 
for $m_-+m_+=m,$ we can eliminate $T^{2m}$ and make 
$0\le l\le 5.$ 
\smallskip 

Note that, generally, $\,l\,$mod $6\,$ is the invariant of the 
action, that is due to Steinberg. 
\smallskip 

Taking the ``simplest" pullbacks for $\tau_{\pm},$ 
we easily check that $l=0$ 
and calculate the images of the generators under 
$\tau_\pm$ and $\sigma.$ 
We arrive at the relations from (\ref{taupsi}). 

\smallskip 
\subsection{Topological Kodaira--Spencer map} 
Generalizing, let $\e$ be an algebraic, or 
complex analytic, or symplectic, or real analytic 
manifold, or similar. It may be noncompact and singular. 
We assume that there is a continuous family of {\em topological 
isomorphisms} $\e\to\e_t$ for manifolds $\e_t$ as $0\le t\le 1,$ 
and that $\e_1$ is isomorphic to $\e_0=\e.$ 

The path $\{\e_t\}$ in the moduli space $\m$ of $\e$ 
induces an outer automorphism $\eps$ of the fundamental 
group $\pi_1(\e,\star)$ defined as above. The procedure 
is as follows. We take 
the image of $\gamma\in \pi_1(\e,\star)$ in 
$\pi_1(\e_1,\star_1)$ for the image $\star_1$ of the 
base point $\star\in \e$ and conjugate it by the path 
from $\star_1$ to $\star.$ Thus, 
{\em the fundamental group $\pi_1(\m)$} (whatever it is) 
{\em acts in $\pi_1(\e)$ by outer automorphisms modulo inner 
automorphisms.} 

The above considerations correspond to the case 
when a group $G$ acts in $\e$ preserving 
a submanifold $D.$ Then  $\pi_1(\m)$ acts in 
$\pi_1^{orb}((\e\setminus D)/G)$ by outer automorphisms. 

Another variant 
is with a Galois group taken instead of $\pi_1(\m)$ 
assuming that $\e$ is an algebraic variety over a 
field that is not algebraically closed. 

The action of $\pi_1(\m)$ on an {\em individual} 
$\pi_1(\e)$  generalizes, in a way, the celebrated 
{\em Kodaira\~Spencer map} and 
is of obvious importance. 
\index{Kodaira\~Spencer map} 
However, calculating the fundamental 
groups of algebraic (or similar) varieties, generally 
speaking, is difficult. The main examples are 
the products of algebraic curves and related 
configuration spaces. 
Not much can be extracted from the action above 
without an explicit description of the fundamental 
group. 
\smallskip

\setcounter{equation}{0} 
\section{Applications} \label{SEC:Appl}
The Verlinde algebras and, more generally, the 
finite dimensional representations of DAHA, have 
quite a few applications. We will mainly discuss 
the non-cyclotomic Gaussian sums and the diagonal 
coinvariants. Both constructions are based on the 
DAHA deformation\~degeneration technique. 
We also discuss a classification of the Verlinde algebras. 
\smallskip 

\subsection{Flat deformation}\label{sec:FLAT} 
Recall that we set $\HH(k)$ when $t=q^k.$ 
We will need the 
{\em little nonsymmetric Verlinde algebra} 
$\widetilde{V}^k$ which is the unique nonzero  
irreducible quotient of $V^k$ upon the restriction to 
the {\em little DAHA} 
$\widetilde{\HH}\equal\langle T,X^2,Y^2\rangle,$ 
corresponding to the subspace
of Laurent polynomials in terms of $X^2.$ 
Its dimension is 
dim$\widetilde{V}^k=$ dim$V^k/2=N-2k.$ 

Generally, 
$\widetilde{V}_{sym}^{k=1}$ is a subalgebra of 
the Verlinde algebra defined for the radical weights 
(from the root lattice) instead of all weights. 
\smallskip 

Let $N=2n+1,\  q^{1/2}=-\exp(\frac{\pi i}{ N}).$ 
Such a choice of $q$ is necessary for the 
positivity of the corresponding inner product. 

We set 
$m= n-k, \, \dim\widetilde{V}^k=N-2k=2m+1.$ 
The aim is to deform $\widetilde{V},$ 
$\widetilde{V}_{sym},$ the projective action 
of $PSL(2,\Z),$ and all other structures making
$q$ arbitrary unimodular.
The construction is as follows. 
\smallskip 

{\em 
For any $q,\,\bar{k}= -\frac{1}{ 2}-m,\,m\in\Z_+,$ 
$\cP$ considered as an $\HH(\bar{k})$\~module 
has a unique irreducible quotient $\overline{V}\neq 0$: 
$$\overline{V}=\hbox{Funct} 
(\frac{\bar{k}+1 }{ 2},...,-\frac{\bar{k}+1 }{ 2}, 
-\frac{\bar{k}}{ 2}),$$ 
$$\dim\overline{V}=2m+1,\  \dim\overline{V}_{sym}=m+1.$$ 
It is unitary as $|q|=1, \, 0<arg(q)<\pi/m.$ 

The representation 
$\overline{V}$, the desired deformation, 
becomes $\widetilde{V}^k$ as 
$$q=e^\frac{2\pi i}{ N},\ \hbox{arg}(q)=\frac{2\pi}{ 2n+1}, 
\ \bar{k}= k-\frac{1}{2}-n,\hbox{\ i.e.\ },\ m=n-k.$$ 
} 
\smallskip 

Concerning the positivity of 
the inner product, note that $\frac{2\pi}{ 2n+1}<\pi/m,$ 
since $k\ge 0.$ 

The deformation construction makes it possible to connect 
the Verlinde algebras and their nonsymmetric $k$\~generalizations 
with the classical representation theory. Since the parameter 
$q$ is now generic, one may expect relations
to $\mathfrak{sl}(2).$ 

Indeed, the representation 
$\overline{V}$ is a $q$\~deformation of an irreducible 
representation 
of $\mathfrak{osp}(2|1)$.
Respectively, 
$\overline{V}_{sym}$ deforms the irreducible representation 
of $\mathfrak{sl}(2)$ of spin=${m/2}.$ We need the rational 
degeneration of DAHA to clarify it.

\subsection {Rational degeneration} 
The {\em rational DAHA} 
is the following quasi-classical limit: 
$$Y=e^{-\sqrt{\kappa} y/2},\ X=e^{\sqrt{\kappa} x},\ 
q=e^{\kappa},\  \kappa\to 0$$ 
of $\HH\,.$ 
Explicitly,  
$$\HH''\equal\langle x,y,s\rangle/\hbox{relations:}$$ 
$$[y,x]=1+2k s,\ s^2=1,\ 
sxs=-x,\ sys=-y.$$ 
\medskip 
Note that the notation $\HH'$ is reserved for 
trigonometric degeneration. 

The algebra $\HH''$ 
acts in the polynomial representation, which is 
$$\C[x],\ s(x)=-x,\ \ 
y\,\mapsto\, D\,=\,\frac{d}{ dx}-\frac{k}{ x}(s-1).$$ 
The square of the {\em Dunkl operator} $D$ is the 
radial part of the Laplace operator: 
$$D^2\mid_{sym}= 
\frac{d^2}{ dx^2}+\frac{2k}{ x}\frac{d}{ dx}.$$ 

The DAHA Fourier automorphism 
$\si$ becomes the outer automorphism corresponding 
to the 
{\em Hankel transform} in the theory of Bessel functions. 
\smallskip 

Let $\bar{k}=-\frac{1}{ 2}-m.$ Then 
$$\lim_{\kappa\to 0}\overline V=\overline V'' 
\equal\C[x]/(x^{2m+1}).$$ 
We call this module a {\em perfect rational} representation. 
The automorphism $\si$ acts there. 
Its symmetric part 
$\overline V''_{sym}$ is nothing else but 
the irreducible representation of 
$\mathfrak{sl}(2)$ of dim$=m+1.$ 
The formulas for the generators of $\mathfrak{sl}(2)$ 
in terms of $x,y$ and the action of $\si$ are 
$$h=(xy+yx)/2,\ e=x^2, f=-y^2, 
\ \  \si 
\hbox{\ becomes\ } w_0. 
$$

\subsection {Gaussian sums} 
The new approach to the Gaussian sums based on 
the deformation construction above is as follows. 

The 
$\tau_-$ acts in $\overline V$ and 
$\overline V_{sym}$, since it 
acts in the polynomial representation 
$\cP.$ In contrast to the polynomial representation, 
$\tau_+$ also acts in 
$\overline V$ and its 
symmetric part. It is the multiplication by 
$$q^{x^2}\equal 
q^\frac{(\pm\bar{k}+j)^2 }{ 4},$$ 
where $\pm=$plus for $j>0$ and minus otherwise. 
Similarly, without going into detail, 
$\si$ is essentially the matrix 
$\Bigl( \hbox{char}_i ( \frac{\pm\bar{k}+j }{ 2}) \Bigr)$ 
for the nonsymmetric characters char$_i\,$ 
($Y$\~eigenfunctions), where 
$i$ and $j$ belongs to the same set (\ref{introbowtie}).

The interpretation of the
Gausssian as $\tau_+$ is of key importance for the 
calculation of the {\em non-cyclotomic Gaussian sum.} 
Generally speaking, a Gauss\~Selberg sum is a 
summation of the Gaussian with respect to a certain 
``measure." Here it becomes
$$\sum_{j=0}^{m} q^{j^2-\bar{k}j} 
\frac{1-q^{2j+\bar{k}}}{ 
1-q^{\bar{k}}} \prod_{l=1}^{2j} 
\frac{1-q^{l+2\bar{k}-1}}{ 
 1-q^{l}} 
$$ 
for {\em any} $q,\ \bar k=-1/2-m$. 
It equals 
$$ \prod_{j=1}^{m} 
\frac{1-q^{2\bar{k}+2j}}{ 1+q^{\bar{k}+2j}}. 
$$ 
\smallskip 

Under the reduction considered above, 
$$N=2n+1,\ m=n-k,\ 0\le k\le n,$$ 
$$q^{1/2}=-\exp(\pi i/ N),\ 
\overline V=\widetilde{V}^k,$$ 
and the product can be somewhat simplified: 
$$\prod_{j=1}^{m} 
\frac{1-q^{2\bar{k}+2j}}{ 1+q^{\bar{k}+2j}}= 
q^{-\frac{m(m+1)}{ 4}}\prod_{j=0}^{m-1}(1-q^{n-j}). 
$$ 

\smallskip 
Now let $m$=$n,$ i.e., let $m$ be the maximal possible. 
Then $k=n-m=0,$ 
$\bar k=-1/2-n=0\, \hbox{mod\,}N,\,$ and 
the ``measure" in the Gaussian sum becomes trivial. 
Setting $l=\frac{n}{ 2}\,\hbox{mod\,} N,$ 
we arrive at 
the identity: 
$$\sum_{j=0}^{N-1} q^{j^2} 
=q^{l^2}\prod_{j=1}^{n}(1-q^{j}). 
$$ 

The product can be readily calculated using 
Galois theory. It equals 
$\sqrt{N}$ for $n=2l$ and $i\sqrt{N}$ otherwise. 
It gives a new proof of the classical Gauss formulas. 
\smallskip 

\subsection{Classification} 
We are going to describe {\em all} nonsymmetric Verlinde 
algebras of type $A_1,$ 
i.e., irreducible quotients of the polynomial 
representation $\cP$ that 
are $PSL(2,\Z)$\~invariant. The 
assumption is that $q^{1/2}$ is a primitive $2N$\~th root 
of unity, $t=q^k.$ All $\si$\~invariant quotients 
of $\cP$ must be through the $\HH\,$\~module 
$V_{4N}=\cP/(X^{2N}+X^{-2N}-(t^{N}+t^{-N})).$ 

Indeed, the invariance gives that 
the element $X^{2N}+X^{-2N},$ central 
in $\HH,$ must act as its $\si$\~image 
$Y^{2N}+Y^{-2N},$ which is $t^N+t^{-N}$ in $\cP.$ 
 
The module $V_{4N}$ is irreducible unless $k$ is 
integral or half-integral. Note that for generic $k$, it 
is {\em not} $\si$\~invariant. Instead, the involution 
$$ 
X\leftrightarrow Y,\ T\leftrightarrow -T^{-1}, 
\ t^{1/2}\mapsto t^{1/2},\ q^{1/2}\mapsto q^{-1/2} 
$$ 
is an inner automorphism of  $V_{4N}$. 
Let us discuss special (``singular") $k$ in detail. 

The module $V_{4N}$ becomes, respectively, 
$$ 
V^{-2}=\cP/(X^{2N}+X^{-2N}-2)\hbox{\ or\ } 
V^{2}=\cP/(X^{2N}+X^{-2N}+2), 
$$ 
as $k\in\Z$ or $k\in 1/2+\Z.$ Respectively, 
$t^{N}=q^{kN}=\pm 1.$ 
The module $V^2$ is an extension of 
$V_{2N}\equal\cP/(X^N+X^{-N})$ 
by its $\varsigma_y$\~image $\varsigma_y(V_{2N}).$ 

The representations $V=V_{2N-4k}$ of dimension $2N-4k$ 
defined above for the integral 
$0<k<N/2$ is a quotient of $V^{-2}.$ The representation 
$\overline{V}=\overline{V}_{2|k|}$ 
of dimension $2|k|$ for the half-integral 
$-N/2<k=-1/2-m<0$ (the notation was $\bar{k}$ was such $k$) 
is a quotient of $V^2.$ 

The construction 
of $V_{2N-4k}$ holds without changes for the positive 
{\em half-integral} $k<N/2,$ 
but then $V_{2N-4k}$  becomes a quotient of 
$V^2.$ The same notation $V_{2N-4k}$ will be used. 
\smallskip 

Let $k\in \BZ/2, |k|<N/2.$ 
The substitution 
$T\mapsto -T,\ t^{1/2}\mapsto -t^{1/2}$ 
identifies the polynomial 
representations for $t^{1/2}$ and $-t^{1/2}.$ 
Thus it is sufficient 
to decompose $\cP$ upon the transformation $k\mapsto k+N,$ 
and we can assume that 
$\-N/2\le k <N/2.$ We will also 
need the outer involutions of $\HH$: 
\begin{align*} 
&\iota:\, T\mapsto -T,\, X\mapsto X,\, Y\mapsto Y, 
\, q^{1/2}\mapsto q^{1/2},\, t^{1/2}\mapsto t^{-1/2},\\ 
&\varsigma_x:\, T\mapsto T,\, X\mapsto -X,\, Y\mapsto Y, 
\, q^{1/2}\mapsto q^{1/2},\, t^{1/2}\mapsto t^{1/2},\\ 
&\varsigma_y:\, T\mapsto T,\, X\mapsto X,\, Y\mapsto -Y, 
\, q^{1/2}\mapsto q^{1/2},\, t^{1/2}\mapsto t^{1/2}. 
\end{align*} 
Let $V^{+}_{2|k|}\equal\overline{V}_{2|k|},\ 
V^{-}_{2|k|}\equal\varsigma_x(\overline{V}_{2|k|}).$ 
\smallskip 

Finally, up to $\iota,\varsigma,$ 
there are {\em three 
different series of nonsymmetric Verlinde 
algebras,} namely, 
$V_{2N-4k}$ (integral $N/2>k>0$), 
$V_{2|k|}$ (half-integral $-N/2<k<0$), and 
$V_{2N+4|k|}$ for {\em integral} $-N/2<k<0.$ 
\smallskip

The latter module is defined as a unique nonzero 
irreducible quotient of $\cP$ for such $k$ 
and is also isomorphic to 
the $\iota\varsigma_y$\~image of the  kernel 
of the map $V^{-2}\to V_{2N-4|k|}.$ It is 
{\em non-semisimple}. The previous two (series of modules) 
are semisimple. All three are projective $PSL(2,\Z)$\~invariant. 
\smallskip
 
The modules $V^{\pm 2}$ are decomposed as follows. 
There are four exact sequences 
\begin{align*} 
&0\to\iota\varsigma_y(V_{2N+4k})\to V^{-2}\to V_{2N-4k}\to 0 
\for k\in \Z_+,\\ 
&0\to \iota(V^{+}_{2k}\oplus V^{-}_{2k}) 
\to V_{2N} \to V_{2N-4k}\to 0 
\for k\in 1/2+\Z_+. 
\end{align*} 
The arrows must be reversed for $k<0$: 
\begin{align*} 
&0\to \iota\varsigma_y(V_{2N-4|k|})\to V^{-2}\to 
V_{2N+4|k|}\to 0 \for k\in -1-\Z_+,\\ 
&0\to \iota(V_{2N-4|k|})\to V_{2N}\to 
V^{+}_{2|k|}\oplus V^{-}_{2|k|}\to 0 
\for k\in -1/2-\Z_+. 
\end{align*} 
\noindent 
Otherwise $V^{-2}$ and $V_{2N}$ are irreducible. 
\medskip 

Recently, a certain non-semisimple variant of the Verlinde algebra 
appeared in \cite{FHST} 
in connection with the fusion procedure for 
the $(1,p)$ Virasoro algebra, although this connection 
is still not justified in full. 
Generally speaking, the fusion procedure for the 
Virasoro-type algebras and the so-called $W$\~algebras 
can lead to non-semisimple Verlinde algebras. 
There are no reasons to expect the existence of a 
positive hermitian inner product there like the Verlinde 
pairing for the conformal blocks, because the corresponding 
physics theories are {\em massless}. 
Surprisingly, the algebra from \cite{FHST} is defined 
using the usual ({\em massive}) Verlinde algebra under certain 
degeneration. Presumably it coincides with the algebra of 
even elements of $V_{2N+4|k|}$ for $k=-1$ or, at least, 
is very close to this algebra. 
\smallskip 

Mathematically, 
the module $V_{2N+4|k|}$ and its multidimensional 
generalizations could be expected to be connected 
with the important problem of describing the 
complete tensor category of the representations 
of Lusztig's quantum group at roots of unity \cite{L2}. 
The Verlinde algebra, the symmetric part of 
$V_{2N-4|k|}$ for $k=1$, describes the so-called 
{\em reduced category} and, in a sense, corresponds to the 
Weyl chamber. The non-semisimple modules of type 
$V_{2N+4|k|}$ are expected to appear in the so-called 
{\em case of the parallelogram}. 
\smallskip 

\subsection{Weyl algebra} 
\vskip 0.3cm 
Before turning to the diagonal invariants, let us 
first discuss the specialization $t=1.$ 
One has: 
$$\HH^{(t=1)}=\w\rsmash\BS_2,\, T\mapsto s\in \BS_2,\, s^2=1,$$ 
where the Weyl algebra, denoted here by $\w,$ is 
a quotient of the algebra of noncommutative polynomials 
$\BC[X^{\pm 1},Y^{\pm 1}]$ with the relations 
$$ 
sXs=X^{-1},\,sYs=Y^{-1},\,XY=YXq^{1/2}. 
$$ 

Given $N\in\BN,$ we set 
$q^{1/2}=\exp(\frac{2\pi i}{N}),$ 
$$ 
\w^\bullet=\w/(X^N=1=Y^N),\ 
\HH^\bullet\equal\HH^{(t=1)}/(X^N=1=Y^N)=\w^\bullet\rsmash\BS_2. 
$$ 
 
The algebra $\w^\bullet$ has a unique irreducible representation 
$$V^\bullet\equal\C[X,X^{-1}]/(X^N=1),\, 
Y(X^m)=q^{-m/2}X^m, m\in \BZ, 
$$ 
which is also a unique irreducible $\HH^\bullet$\~module. 

Moreover, $\tau$ and $\si$ act in $V^\bullet.$ 
Recall that the action of $\si$ in $\HH^{(t=1)}$ is as follows: 
$$X\mapsto Y^{-1},\ Y\mapsto X,\ s\mapsto s. 
$$ 

The problem is that the above description 
of $V^\bullet$ doesn't make the $\si$\~invariance 
clear and is inconvenient when studying the $\si$\~action, 
playing a major role in the theory of theta functions 
and automorphic forms. 

To make the $X$ and the $Y$ on equal footing, we can
take a $\si$\~eigenvector $v\in V^\bullet$ and set 
$$ 
V^\bullet=(\w/\j_v)(v)\hbox{\ for\ the\ }\si-\hbox{invariant\ 
ideal\ } \j_v=\{H \mid H(v)=0,H\in\w\}. 
$$ 
This presentation makes the $\si$\~invariance of $V^\bullet$ 
obvious, but finding ``natural" 
$\si$\~eigenvectors in $V^\bullet$ with reasonably 
explicit $\j_v$ is not an easy problem. 
The simplest case is $N=3,$ dim$V^\bullet=3,$ where 
we can proceed as follows. 

The space $\{v\in V^\bullet\mid s(v)=-v\}$ is one-dimensional 
for $N=3.$ It is nothing but the space 
of odd characteristics in the classical theory of one-dimensional 
theta functions. 
It equals $\BC d$ for  $d=X-X^{-1}.$ 
Since $\si(s)=s,$ this space must be $\si$\~invariant. 
Then, calculating $\j_d$ is simple: 
$$ 
\j_d=\w J_o, \where 
J_o=\{H\in \w^\bullet\, \mid sH=Hs,\, 
H(d)=0\in V^\bullet. 
$$ 
One has: 
\begin{align*} 
\{Y+Y^{-1}\}(d)&= \{(Y+Y^{-1})(X-X^{-1})\}(1)= 
\{q^{-1/2}(XY-X^{-1}Y^{-1})\\ 
+&q^{1/2}(XY^{-1}-X^{-1}Y)\}(1)=(q^{1/2}+q^{-1/2})(d)=-d,\\ 
\{X+X^{-1}\}(d)&=X^2-X^{-2}=X^{-1}-X=-d. 
\end{align*} 
Thus, $Y+Y^{-1}+1,\, X+X^{-1}+1\, \in J_o,$ 
and we can continue: 
$$\{YX+Y^{-1}X^{-1}\}(d)=q^{-1}X^2-1+1-q^{-1}X^{-2}=-q^{-1}d.$$ 
Actually, this calculation is not needed, since 
$Y+Y^{-1}+1$ and $X+X^{-1}+1$ algebraically generate $J_o,$ 
(it is not true in the commutative polynomials!) and 
we can use the Weyl relations.
 
Finally, we arrive at 
the equality $\j_d=\w J_o.$ 
The next section contains a generalization of this 
construction for $\HH\,.$ 

\smallskip 
Using the non-commutativity of $X$ and $T$ in
the latter argument is directly connected with 
the recent results on $W$\~invariant polynomial differential 
operators due to Wallach, Levasseur, Stafford, and Joseph. 
They prove that the invariant operators are generated 
by the $W$\~invariant polynomials and 
the $W$\~invariant differential operators with constant 
coefficients, using that they {\em do not} commute. 
Of course it is not true for the polynomials in terms 
of two sets of variables. 
\smallskip 

\subsection {Diagonal coinvariants} 
The previous construction can be naturally generalized to 
the case of an arbitrary Weyl group 
$W$ acting in $\R^n.$ We set
$$
w^{-1}f(x,y)=f(wx,wy)\hbox{ for a function } f
\hbox{ in terms of } x,y\in \R^n,
$$
$$\I_o\equal \hbox{Ker}(\C[x,y]_{sym}\ni f\mapsto f(0,0)),\and 
\nabla\equal\C[x,y]/(\C[x,y]\I_o).
$$ 

Haiman conjectured that 
$\nabla$ has a natural quotient $\V$ 
of dim= $(h+1)^n$ for the Coxeter number $h.$ 
It must be isomorphic to the space Funct($Q/(h+1)Q$) 
as a $W$\~module for the root lattice $Q$ 
and have a proper character with 
respect to the degree in terms of $x$ and $y.$ 
Note that $Q/(h+1)Q$ is isomorphic to $P/(h+1)P$ 
for the weight lattice $P$, since the order $[P:Q]$ is 
relatively prime with $(h+1).$ 

He proved that 
$\nabla$=$\V$ and dim$\nabla=(n+2)^{n}$ in the case of 
$A_{n}\,$ for $W$=$\BS_{n+1}.$ 
The coincidence is a very special feature of $A_{n}.$ 
As far as I know, no general uniqueness 
claims about $\V$ as a graded vector space and 
as a $W$\~module were made. 
One needs the double Hecke algebra to make this quotient 
``natural." 

Haiman's conjecture was recently 
justified by Gordon. He proves that
$$\V=gr \overline V'' \ \hbox{for\ } \bar k=-1/h-1,$$ 
where $\overline V''$ is a $W$\~generalization 
of the {\em perfect rational} representation $\overline V''$ 
considered above as $m=1.$ By $gr,$ we mean taking the 
graded vector space of $\overline V''$ 
with respect to the degree 
in terms of $x$ and $y;$  see 
Section \ref{sec:DIAGCOI}. 

Recall that $h=2$ and, generally, 
$\dim \overline V''=2m+1$ in the case of $A_1.$ 
Here we make $m=1,$ so the representation 
considered by Gordon is of dimension $3$ 
for $A_1.$ 
Even in this case, the coincidence $gr \overline V''=\nabla$ 
is an instructional exercise. It is not immediate 
because $\nabla$ is given in terms of {\em double} 
polynomials and $\overline V''$ was defined as a quotient 
of the space of {\em single} polynomials; cf. the previous 
section. 
\smallskip 

Adding $q$ to the construction, 
we obtain the following theorem, describing a 
$W$\~genera\-li\-zation of the perfect module $\overline V$ 
considered above for the simplest nontrivial $\bar k=-1-1/h.$ 

\smallskip 
{\em 
(a) There exist Lusztig-type 
(exp-log) isomorphisms 
$\overline{V}''\rightleftharpoons \overline{V} 
\rightleftharpoons V^\bullet$ for 
$\bar k$=$-{1/ h}-1.$ Here 
$V^\bullet$ is a unique 
$\HH^\bullet$\~quotient of $\cP=\C[X,X^{-1}]$ 
for the reduction: 
$$\bullet= 
\{q^{1/[P:Q]}=\exp(\frac{2\pi i}{ N}), \, N\equal h+1,\ 
t=q^{\bar k}=1\}.$$ 

(b) Concerning the coinvariants, 
$$\overline{V}=\HH/\HH I_o(\De),\ 
I_o=\hbox{Ker}(\HH_{sym}\ni H\mapsto H(\De) \in 
\overline{V}),$$ 
where $\HH_{sym}$ is a subalgebra of the elements 
of $\HH\,$ commuting with $T.$ Here 
$\De$ is the discriminant 
$\De=t^{1/2}X-t^{-1/2}X^{-1}$ for $n=1$; generally, it is a 
product of such binomials over positive roots. 
} 
\smallskip 

The second part demonstrates how double polynomials 
appear (the general definition of $\overline V$ is given 
in terms of single polynomials). Using that $\De$ 
generates the one-dimensional ``sign-representation" 
of the nonaffine Hecke algebra $\lan T \ran,$ we can 
naturally identify $\overline V$ with quotients of the 
space of double Laurent polynomials. 

\smallskip 

Since $t=1$ for the $\bullet$\~reduction, one has:  
$\HH^\bullet$=(Weyl algebra)$\rsmash W.$ Thus 
$V^\bullet$ is its unique irreducible representation 
under the relations $X^N$=1=$Y^N$ imposed, 
namely, 
$V^\bullet\simeq\C[X,X^{-1}]/(X^N=1).$ 
Generally, its dimension is 
dim=$N^{rank}$ (rank = the number of $X$\~generators), 
exactly what Haiman conjectured. Note that given $n,$ 
part (b) holds only thanks to a very special 
choice of $N.$ 

In Gordon's proof, the $(h+1)^{n}$\~formula requires 
a construction of the resolution of $\overline V''$ and 
more. The theorem above gives an explanation of why 
the dimension is so simple. The dimension looks like that 
from the theory of Weyl algebra, although the 
definition of the space of coinvariants has nothing 
to do with the roots of unity, and Gordon's theorem 
really belongs to this theory! 

\smallskip 
\setcounter{equation}{0} 
\section{DAHA and {\em p}--adic theory} \label{SEC:DAHAp}
This section is devoted to double  Hecke 
algebras in the general setting and their connection 
with the classical $p$\~adic spherical transform. 
See paper \cite{C12} for 
a complete theory. 

\subsection{Affine Weyl group} Let $R\subset \BR^n$ be a 
simple reduced root system, $R_+\subset R$ the 
set of positive roots, $\{ \alpha_1, \ldots, 
\alpha_n\} \subset R_+$ the corresponding set 
of simple roots, and $\vartheta \in R_+$ 
the {\em maximal coroot}
\index{maximal short root $\vth$} 
that is the longest positive root in $R^\vee.$ 
We normalize the scalar product on $\BR^n$ by 
the condition $(\vartheta,\vartheta)=2,$ 
so $\vartheta$ belongs to $R$ as well and it
is the longest {\em short} root there. 
For any $\alpha \in R,$ its dual coroot is 
\index{coroot $\alpha^\vee$} 
$\alpha^\vee $ $ =\frac 
{2\alpha}{(\alpha,\alpha)} $ $=\alpha/\nu_{\alpha}.$ 
So $\nu_\alpha=(\alpha,\alpha)/2=1,2,3.$ 
We set $\rho=(1/2)\sum_\al \al,\,$ 
$\rho^\vee=(1/2)\sum_\al \al^\vee.$

The affine roots are 
$$ 
\tilde{R}=\{\tilde{\alpha}=[\alpha,\nu_\alpha j]\},\ j\in \BZ. 
$$ 
Note the appearance of $\nu_\alpha$ here. It is because of our 
nonstandard choice of $\vartheta.$ 
We identify nonaffine roots $\alpha$ 
with $[\alpha,0]$ 
and set $\alpha_0=[-\vartheta,1].$ 
For $\tilde{\alpha} \in \tilde{R},$ 
$s_{\tilde{\alpha}} \in $ $ \Aut(\BR^{n+1})$ is the 
reflection 
\begin{align}\label{fwaffine} 
&s_{\tilde{\alpha}}([x,\zeta])=[x,\zeta]-
2\frac{(x,\alpha)}{(\alpha,\alpha)}\tilde{\alpha}, 
\hbox{\ and\ }\notag\\ 
&W=\lr s_\alpha\, \mid \, \alpha \in R\rr,\ 
\tW=\lr s_{\tilde{\alpha}}\, \mid \, \tilde{\alpha} \in \tilde{R}\rr. 
\end{align} 

We set $s_i=s_{\alpha_i}.$ It is well known 
that $\tW$ is a Coxeter group with the  generators 
$\{s_i\}$. 
Let $\omega_i \subset \BR^n$ be the 
fundamental weights: $(\omega_i ,\alpha^\vee_j)=\delta_{ij}$ 
for $1\le i,j\le n,$  $P= 
\oplus_{i=1}^n\BZ \omega_i$ the weight 
lattice, and  $P_+ \equal \oplus_{i=1}^n\BZ_+\omega_i$ the 
cone of dominant weights. \index{extended Weyl group} 
\index{affine Weyl group $\widehat{W}$} 
We call $\hW\equal W\lsmash P$ the 
{\em extended affine Weyl group}: 
$$ 
wb([x,\zeta])=[w(x),\zeta-(b,x)]\hbox{\ for\ } b\in P,x\in \BR^n. 
$$ 

\noindent 
The \index{length $l(\hw)$} {\em length function} 
$l: \hW\to \BZ_+$ is given by the formula 
\begin{equation}\label{falength} 
l(wb)= 
\sum_{\begin{array}{c}\alpha \in R_+\\ 
w(\alpha^\vee)\in R_+\end{array}}|(b,\alpha^\vee)|+ 
\sum_{\begin{array}{c}\alpha \in R_+\\ 
w(\alpha)\in -R_+\end{array}}|(b,\alpha^\vee)+1|, 
\end{equation} 
where $w\in W$, $b\in P$. Let $Q = 
\oplus_{i=1}^n\BZ \alpha_i$ be the root 
lattice. Then $\tW=W\lsmash Q \subset \hW$ is a normal 
subgroup and $\hW/\tW=P /Q.$ Moreover, 
$\hW$ is the semidirect product $\Pi \lsmash \tW,$ where 
$$ 
\Pi =\{ \pi\in \hW\,\mid\, l(\pi)=0\}= 
\{ \pi\in \hW\,\mid\, \pi:\{\alpha_i\}\mapsto 
\{\alpha_i\}\, \},\ 0\le i\le n. 
$$ 
It is isomorphic to $P/Q$ and acts naturally on 
the affine Dynkin diagram for $R^\vee$ 
with the reversed arrows. It is not 
the standard Dynkin diagram for $R$ 
because of our choice of $\vartheta.$ 
We set $\pi(i)=j$ as $\pi(\alpha_i)=\alpha_j.$ 

\subsection{Affine Hecke algebra} 
We denote the affine Hecke algebra by $\cH.$ 
Its generators are $T_i$ for  $i=0, \ldots, n$ and 
 $\pi \in \Pi$; the relations are 
\begin{align}\label{ftt} 
&\underbrace{T_iT_jT_i\ldots}_{\mbox{order of}\; s_is_j} 
= \underbrace{T_jT_iT_j\ldots}_{\mbox{order of}\; s_is_j},\ \, 
\pi T_i\pi^{-1}=T_{\pi (i)},\\ 
&(T_i-t^{1/2})(T_i+t^{1/2})=0 \hbox{\ for\ } 
0\le i\le n,\ \pi\in \Pi.\notag 
\end{align} 
If $\hat w=\pi s_{i_l}\ldots s_{i_1}\in \hW$ is a reduced 
expression, i.e., $l(\pi^{-1}w)=l,$ 
we set $T_{\hat w}=\pi T_{i_l} 
\ldots T_{i_1}$. The elements 
$T_{\hat w}$ are well defined and form a basis of $\cH$. 

Let $G$ be the adjoint split $p$\~adic 
simple group corresponding to $R$. Then $\cH$ is 
the convolution algebra of compactly supported 
functions on $G$ that are 
left-right invariant with respect to the Iwahori subgroup $B,$ 
due to Iwahori and Matsumoto. Namely, the $T_i$ are the 
characteristic functions of the double cosets 
$Bs_iB,$ where we use a natural embedding 
$W\to G.$ Generally, it is not a homomorphism. 

To be more exact, the $p$\~adic quadratic equations 
are represented in the form $(T_i-1)(T_i+t_i)=0$ for 
the standard normalization of the Haar measure. Here 
the $t_i$ may depend on the length of $\alpha_i$ 
and are given in terms of the cardinality 
of the residue field. We 
will stick to our normalization of the $T$ and assume that 
the {\em parameters $t$ coincide} to simplify the formulas 
of this section. 
We will also  use 
$$\delta_{\hat w}\equal t^{-l(\hat w)/2}T_{\hat w}, 
$$ 
which satisfy the quadratic equation with $1.$ 

Let $\Delta$ be the left regular representation 
of $\cH$. In the basis $\{\delta_{\hat w}\},$ the representation 
$\Delta$ is given by 
\begin{equation}\label{fthecke} 
T_i\delta_{\hat w}=\left\{ \begin{array}{cc}t^{1/2} 
\delta_{s_i\hat w}& \mbox{if\ } l(s_i\hat w)=l(\hat w)+1,\\ 
t^{-1/2}\delta_{s_i\hat w}+ 
(t^{1/2}-t^{-1/2})\delta_{\hat w}& 
\mbox{if\ } 
l(s_i\hat w)=l(\hat w)-1\end{array}\right. 
\end{equation} 
and the obvious relations $\pi \delta_{\hat w}=\delta_{\pi\hat w},$ 
where $\pi\in \Pi$ and  $\hat w\in\hW.$ 

The \index{spherical rep $\Delta^+$} 
{\em spherical representation} 
appears as follows. 
Let 
$$\delta^+\equal (\sum_{w\in W}t^{l(w)})^{-1} 
\sum_{w\in W}t^{l(w)}\delta_w\in \cH. 
$$ 
One readily checks 
that $T_i\delta^+=t^{1/2}\delta^+$ for $i=1, \ldots, n,$ 
and $(\delta^+)^2=\delta^+$. We call $\delta^+$ 
the \index{symmetrizer} $t$\~{\em symmetrizer}. 
Then 
$$\Delta^+\equal \Delta \delta^+=\oplus_{b\in P}\BC 
\delta^+_{b},\ \delta^+_{\hat w}=\delta_{\hat w}\delta^+, 
$$ 
is an $\cH$\~submodule of $\Delta.$ 

It is nothing   but Ind$_H^\cH 
(\BC_{t^{1/2}}),$ where $H\subset \cH$ is the subalgebra 
generated by $T_i, i=1, \ldots, n$ and $\BC_{t^{1/2}}$ 
is the one-dimensional representation of $H$ defined by 
$T_i\mapsto t^{1/2}.$ 

Due to Bernstein, Zelevinsky, and Lusztig 
(see, e.g., \cite{L}), we 
set $Y_a=T_a$ for $a\in P_+$ and extend it to the whole $P$ 
using $Y_{b-a}=Y_bY_a^{-1}$ for dominant $a,b.$ 
These elements are well defined and pairwise commutative. 
They form 
the subalgebra $\cY\cong \BC[P]$ inside $\cH.$ 

Using $Y,$ one can omit $T_0.$ Namely, 
the algebra $\cH$ is generated by $\{T_i,i>0,\, Y_b\}$ 
with the following relations: 
\begin {align} 
&T^{-1}_iY_b T^{-1}_i=Y_b Y_{\alpha_i}^{-1}\hbox{\ if\ } 
(b,\alpha^\vee_i)=1, \label{ftyt} 
\\ 
&T_iY_b=Y_b T_i \hbox{\ if\ } (b,\alpha^\vee_i)=0,\  1 \le i\le  n. 
\notag 
\end{align} 

The PBW theorem for $\cH$ gives that 
the spherical representation 
can be canonically identified with $\cX,$ as $\delta^+$ 
goes to $1.$ The problem is to calculate this 
isomorphism explicitly. It will be denoted by $\Phi.$ 
We come to the definition of Matsumoto's spherical 
functions: 
$$ 
\phi_{a}(\Lambda)=\Phi(\delta^+_a)(Y\mapsto\Lambda^{-1}), 
\ a\in P.$$ 
Here, by $Y\mapsto \Lambda^{-1},$ we mean that $\Lambda_b^{-1}$ 
substitutes for $Y_b.$ 
See \cite{Ma}. 
In the symmetric ($W$\~invariant) case, the spherical 
functions are due to Macdonald. 

By construction, 
$$\phi_a=t^{-l(a)}\Lambda_a^{-1}=t^{-(\rho^\vee,a)} 
\Lambda_{a}^{-1}\hbox{\ for\ } a\in P_+,\ 
\rho^\vee= (1/2)\sum_{\alpha>0}\alpha^\vee. 
$$ 
The formula $l(a)=(\rho^\vee,a)$ readily results 
from (\ref{falength}). 
So the actual problem 
is to calculate $\phi_a$ for non-dominant 
$a.$ 

{\bf Example.} 
Consider a root system of type $A_1$. In 
this case, $P=\BZ\omega,$ where $\omega\in \BZ$ is the fundamental 
weight, $Q=2\BZ,$ and $\omega=\pi s$ for $s=s_1.$ 
We identify $\Delta^+$ and $\cY,$ so $\delta^+=1.$ 
Letting $Y=Y_\omega,T=T_1,$ we get 
$Y_{m}\equal Y_{m\omega}=Y^m$ and 
$$\phi_m\equal \phi_{m\omega}= 
t^{-m/2}\Lambda^{-m},\hbox{\ for\ } m\ge 0. 
$$ 
Note that $TY^{-1}T=Y$ and $\pi=Y T^{-1}.$ 

Let us check that 
\begin{equation}\label{fxphi} 
\Lambda\phi_{-m}= 
t^{1/2}\phi_{-m-1}-(t^{1/2}-t^{-1/2})\phi_{m+1},\ m>0. 
\end{equation} 
Indeed, $\phi_{-m}=t^{-m/2}(T\pi)^m(1) 
\mid_{Y\mapsto \Lambda^{-1}}\,$ and 
\begin{align} 
&Y^{-1}\phi_{-m}\ =\ t^{-m/2}(T^{-1}\pi)(T\pi)^m(1)\\ 
=&t^{-m/2}(T-(t^{1/2}-t^{-1/2}))\pi(T\pi)^m(1)\notag\\ 
=&t^{-m/2}(T\pi)^{m+1}(1)-
-t^{-m/2}(t^{1/2}-t^{-1/2})(\pi T)^m\pi(1)\notag\\ 
=&t^{1/2}\phi_{-m-1}(Y^{-1}) 
-t^{-m/2}(t^{1/2}-t^{-1/2})(\pi T)^m 
(t^{-1/2}\pi T)(1)\notag\\ 
=&t^{1/2}\phi_{-m-1}(Y^{-1})-
(t^{1/2}-t^{-1/2})\phi_{m+1}(Y^{-1}).\notag 
\end{align} 

\subsection{Deforming {\em p}--adic formulas} 
The following chain of theorems represents 
a new vintage of the classical theory. 
We are not going to prove 
them here. Actually, all claims that are beyond the 
classical theory of affine Hecke algebras can be 
checked by direct and not very difficult 
calculations, with a reservation about 
Theorems \ref{tinverse} and \ref{tgeps}.

\begin{theorem}\label{t111} 
Let $\xi \in \BC^n$ be a 
fixed vector and let $q\in \BC^*$ be a fixed scalar. We 
represent $\hat w=bw,$ where $w\in W, b\in P$. 
In $\Delta_q^\xi\equal \oplus_{\hat w\in \hW} 
\BC\delta_w^\xi,$ 
the formulas 
$\pi \delta_{\hat w}^\xi =\delta_{\pi \hat w}^\xi$ and 
$$T_i\delta_{\hat w}^\xi=\left\{ \begin{array}{cl} 
\frac{t^{1/2}q^{(\alpha_i,w(\xi)+b)}-t^{-1/2}} 
{q^{(\alpha_i,w(\xi)+b)}-1}\delta_{s_i\hat w}^\xi -
\frac{t^{1/2}-t^{-1/2}}{q^{(\alpha_i,w(\xi)+b)}-1} 
\delta_{\hat w}^\xi&\mbox{if}\ i> 0, \\ 
\frac{t^{1/2}q^{1-(\vartheta,w(\xi)+b)}-t^{-1/2}} 
{q^{1-(\vartheta,w(\xi)+b)}-1}\delta_{s_0\hat w}^\xi -
\frac{t^{1/2}-t^{-1/2}}{q^{1-(\vartheta,w(\xi)+b)}-1} 
\delta_{\hat w}^\xi&\mbox{if}\ i=0\end{array}\right.$$ 
define a representation of the 
algebra $\cH,$ provided that all denominators are nonzero, 
i.e., $q^{(\alpha,b+\xi)}\neq 1$ for 
all $\alpha\in R,\, b\in P.$ 
\end{theorem} 

The regular representation $\Delta$ with the 
basis $\delta_{\hat w}$ is the limit of representation 
$\Delta_{q}^\xi$ as $q\to \infty,$ provided that 
$\xi$ lies in the fundamental 
alcove: 
\begin{equation}\label{falcove} 
(\xi,\alpha_i)>0 \hbox{\ for\ } i=1, 
\ldots, n,\ (\xi,\vartheta)<1. 
\end{equation} 

We see that the representation $\Delta_q^\xi$ 
is a flat deformation of $\Delta$ for such $\xi.$ 
Moreover, $\Delta_q^\xi\cong \Delta.$ 
This will readily follow from the next theorem. 

Note that taking $\xi$ in other alcoves, we 
get other 
limits of $\Delta_{q}^\xi$ as $q\to \infty.$ 
They are isomorphic to 
the same regular representation; however, the formulas 
do depend on the particular alcove. 
We see that the regular representation has rather many 
remarkable systems of basic vectors. They are not 
quite new in the theory of affine Hecke algebras, but 
such systems were not studied systematically. 

\begin{theorem}\label{tdaha} 
(i) We set  $X_b(\delta^\xi_{\hat w})= 
q^{(\alpha_i,w(\xi)+b)}\delta^\xi_{\hat w}\,$ for $\hat w=bw,$ 
where we use the notation 
$X_{[b,j]}=q^j X_b.$ 
These operators have a simple spectrum 
in $\Delta^\xi_q$ under the conditions of the 
theorem and 
satisfy the relations dual to (\ref{ftyt}) with $i=0$ added: 
\begin {align} 
&T_iX_b T_i=X_b X_{\alpha_i}^{-1}\hbox{\ if\ } 
(b,\alpha^\vee_i)=1, \label{ftxt} 
\\ 
&T_iX_b=X_b T_i \hbox{\ if\ } (b,\alpha^\vee_i)=0,\  0 \le i\le  n, 
\notag\\ 
&\hbox{and\ moreover,\ }\pi X_b \pi^{-1}= 
X_{\pi(b)} \hbox{ \ as\ } \pi\in \Pi.\notag 
\end{align} 

(ii) The {\em double affine Hecke algebra} 
$\HH\, $ is defined by imposing 
relations (\ref{ftt}) and (\ref{ftxt}). Then 
the representation $\Delta^\xi_q$ is nothing 
but the induced representation 
$$\Ind_{\cX}^{\HH\,}(\BC\delta^\xi_{id}),\ \cX\equal \BC[X_b], 
$$ 
which is isomorphic to $\Delta$ as an $\cH$\~module. 
\end{theorem} 

\rmk 
The operators $X_b$ are in a way the coordinates 
of the Bruhat-Tits buildings corresponding to the 
$p$\~adic group $G.$ In the classical theory, we use 
only their combinatorial variants, namely, the distances 
between vertices, which are integers. 
The $X$\~operators clarify dramatically 
the theory of the $p$\~adic spherical Fourier transform, 
because, as we will see, they are the ``missing" 
Fourier\~images of the $Y$\~operators. Obviously, the 
$X_b$ do not survive in the limit $q\to \infty$; 
however, they do not collapse completely. Unfortunately, 
the Gaussian, which is $q^{x^2/2}$ as $X_b=q^{(x,b)},$ 
does. 
\sq 

This theorem is expected to be connected 
with \cite{HO2} and via this paper with \cite{KL1}. 
Here we will stick to the spherical representations. 
The following theorem establishes a connection 
with $\Delta^+.$ 

\begin{theorem}\label{t113} 
We set $q^\xi=t^{-\rho},$ i.e., 
$q^{(\xi,b)}\mapsto t^{-(\rho,b)}$ for all 
$b\in P.$ The corresponding representation 
will be denoted by $\widetilde{\Delta}.$ 
It is well defined for generic $q,t.$ 
For any $b\in P,$ let $\pi_b$ 
be the minimal length representative in the set 
$\{bW\}.$ It equals $\,bu_b^{-1}$ for the length-minimum 
element $u_b\in W$ such that $u_b(b)\in -P_+.$ 
Setting $\delta_b^\#\equal 
\delta_{\pi_b}$ in $\widetilde{\Delta},$ 
the space 
$\Delta^\#\equal 
\oplus_{b\in P}\BC \delta_b^\#$ is an 
$\HH\,$\~submodule of $\widetilde{\Delta}.$ 
It is isomorphic to $\Delta^+$ as an $\cH$\~module. 
\end{theorem} 

The representation $\Delta^\#$ is described by the 
same formulas from Theorem \ref{t111}, which 
vanish {\em automatically} on $s_i\pi_b$ not in the 
form $\pi_{c}$, thanks to the special choice of $q^\xi.$ 
It results directly from the following: 
$$s_i\pi_b=\pi_c\Leftrightarrow 
(\alpha_i,b+d)\neq 0,\hbox{\ where\ } 
(\alpha_i,d)\equal \delta_{i0}. 
$$ 
Here $c=b-((\alpha_i,b+d)\alpha_i^\vee$ 
for $\alpha_0^\vee\equal -\vartheta.$ 

We define the action $(\!(\ )\!)$ of $\hW$ on $\BR^n$ 
by the formulas $wa(\!(x)\!)=w(a+x).$ 
The above $c$ is $s_i(\!(b)\!).$ 
This action is constantly used in the theory 
of Kac\~Moody algebras. It is very convenient 
when dealing with $\Delta^\#.$ Note that 
$$ 
\pi_{b}(\!(c)\!) =bu_b^{-1}(\!(c)\!) =u_b^{-1}(c)+b 
\hbox{\ for\ } b,c\in P. 
$$ 

Let us calculate the formulas from Theorem \ref{t111} 
upon $q^\xi\mapsto t^{-\rho}$ as 
$q\to \infty.$ This substitution changes the consideration, 
but not too much: 
\begin{equation}\label{fdeltas} 
T_i\delta_b^\#=\left\{ \begin{array}{cl} 
 t^{1/2}\delta_{s_i(\!(b)\!)}^\# 
&\mbox{\ if\ } (\alpha_i,b+d)>0, \\ 
t^{-1/2}\delta_{s_i(\!(b)\!)}^\# + 
(t^{1/2}-t^{-1/2})\delta_{b}^\#\ 
&\mbox{\ as\ } (\alpha_i,b+d)<0. 
\end{array}\right. 
\end{equation} 
Otherwise it is zero. 
The formulas $\pi \delta_{b}^\# =$ 
$\delta_{\pi(\!(b)\!)}^\#$ hold for arbitrary 
$\pi\in \Pi, b\in P.$ 

Since this calculation is different from that for generic 
$\xi,$ it is not surprising that (\ref{fdeltas}) 
does not coincide 
with (\ref{fthecke}) restricted to $\hat w=b$ 
and multiplied on the right by the $t$\~symmetrizer $\delta^+.$ 
The representations $\lim_{q\to\infty} \Delta^\#$ and 
$\Delta^+$ are equivalent, but the $T$\~formulas with 
respect to the limit of the basis 
$\{\delta^\# \}$ are 
different from those in terms of the classical basis 
$\{\delta_b^+=\delta_b\delta^+\}.$ 

\subsection{Fourier transform} 
In the first place, 
Macdonald's nonsymmetric polynomials generalize 
the Matsumoto spherical functions. 
We use $\pi_b=bu_b^{-1}=$ $\hbox{Min-length}\, \{bw,\,w\in W\}.$ 

\begin{theorem}\label{tpolyn} 
(i) Let $\cP$ be the representation of the double 
affine Hecke algebra $\HH\,$ in the space of 
Laurent polynomials $\cP=\BC[X_b]$: 
\begin{align}\label{ftonx} 
&T_i\  = \  t_i ^{1/2} s_i\ +\ 
(t_i^{1/2}-t_i^{-1/2})(X_{\alpha_i}-1)^{-1}(s_i-1), 
\ 0\le i\le n,\\ 
&X_b(X_c)=X_{b+c},\ \pi(X_b)=X_{\pi(b)},\, \pi\in \Pi, 
\hbox{\ where\ } 
X_{[b,j]}=q^jX_b.\notag 
\end{align} 

(ii) For generic $q,t,$ 
the polynomials $\eps_b$ are uniquely defined from 
the relations: 
\begin{align} \label{fepsb} 
&Y_a(\eps_b)=t^{(u(\rho),a)}q^{-(b,a)}\eps_b,\ 
\hbox{\ where\ } \pi_b=bu \hbox{\ for\ } u\in W \\ 
&\eps_b(t^{-\rho})=1, \hbox{\ where\ } 
X_b(t^{-\rho})=t^{-(b,\rho)}.\notag 
\end{align} 

(iii) Setting 
$X_b^*=X_{-b},\ q^*=q^{-1},$  and $ t^*=t^{-1},$ 
the limit of $\eps^*_b(X\mapsto \Lambda)$ as $q\to \infty$ 
coincides with $\phi_b(\Lambda)$ for $b\in P.$ 
\end{theorem} 

\rmk 
Note that the $\Delta^\xi$\~formulas from 
Theorem \ref{t111} are actually the evaluations 
of (\ref{ftonx}) at $q^\xi.$ To be more exact, there 
is an $\HH\,$\~homomorphism from $\cP$ to the $\HH\,$\~module 
of functions on $\Delta^\xi.$ 
For instance,  Theorem \ref{t111} can be 
deduced from Theorem \ref{tpolyn}. 
The formulas for the polynomial representation of the 
double affine Hecke algebra $\HH\,$ are nothing   but 
the {\em Demazure-Lusztig operators} 
in the affine setting. 
\ 

We are going to establish a Fourier-isomorphism 
$\Delta^\#\to \cP,$ which is a generalization 
of the Macdonald\~Mat\-su\-moto 
inversion formula. We use the constant term functional on Laurent 
series and polynomials  denoted by $\lr\,\rr.$ 

The first step is to make both representations unitary using 
\begin{align} 
&\mu\ =\ \prod_{\tilde{\alpha} \in \tilde{R}} 
\frac{1-X_{\tilde{\alpha}}}{1-tX_{\tilde{\alpha}}},\ 
\mu^0=\mu/\lr\mu\rr, 
\label{fmugen}\\ 
&\mu^1(\pi_b)=\mu(t^{-\pi_b(\!(\rho)\!)})/\mu(t^{-\rho}), 
\ \pi_b=bu_b^{-1}.\notag 
\end{align} 
Here we treat $\mu$ as a Laurent series to define $\mu^0.$ 
The coefficients of $\mu^0$ are rational functions in terms 
of $q,t.$ 
The values $\mu^1(\pi_b)$  are rational functions in terms 
of $q,t^{1/2}.$ 

The corresponding pairings are 
\begin{align} 
&\lr f\, ,\, g\rr_{\,pol}=\lr f\, T_{w_0}w_0(g(X^{-1})\, \mu^0\rr, 
\ f,g\in \cP,\notag\\ 
&\lr \sum f_b\delta^\#_b ,\sum g_b\delta^\#_b\rr_{\,Del}= 
\sum (\mu^1(\pi_b))^{-1}\,f_b g_b.\notag 
\end{align} 
Here $w_0$ is the longest element in $W.$ Note that 
the element $T_{w_0}^2$ 
is central in the nonaffine Hecke algebra $\H$ 
generated by $\{T_i,i>0\}.$ Both pairings 
are well defined and symmetric. Let us give the formulas for the 
corresponding anti-involutions: 
\begin{align} 
T_i\mapsto T_i,\, X_b\mapsto X_b,\, T_0\mapsto T_0,\, 
Y_b\mapsto T_{w_0} Y^{-1}_{w_0(b)}T_{w_0}^{-1}\ \ \ 
&\hbox{\ in\ }\Delta^+,\notag\\ 
T_i\mapsto T_i, \, Y_b\mapsto Y_b,\, 
T_0\mapsto T_{s_\vartheta}^{-1}Y_{\vartheta},\, 
X_b\mapsto T_{w_0}^{-1}X^{-1}_{w_0(b)}T_{w_0} 
&\hbox{\ in\ }\cP,\notag 
\end{align} 
where $1\le i\le n,\, b\in P.$ 

\begin{theorem}\label{tinverse} 
(i) Given $f=\sum_b f_b\delta^\#_b \in \Delta^\#,$ we 
set $\widehat{f}=\sum_b\, f_b\,\eps_c^*\in \cP,$ 
where $X_b^*=X_{b}^{-1}, \ q^*=q^{-1},t^*=t^{-1}.$ 
The {\em inversion} 
\index{inversion (DAHA)} 
 of this transform is as follows: 
\begin{equation}\label{fhinv} 
f_b\ =\ t^{l(w_0)/2}(\mu^1(\pi_b))^{-1} 
\lr\, \widehat{f},\, \eps^*\rr_{\,pol}. 
\end{equation} 

\index{Plancherel $\HH$} 
(ii) The {\em Plancherel formula} reads 
\begin{equation}\label{fhplan} 
\lr f\,,\, g\rr_{\,Del}\ =\ 
t^{l(w_0)/2}\lr \widehat{f}\,,\, \widehat{g}\rr_{\,pol}. 
\end{equation} 
Both pairings are positive definite over $\BR$ if 
$\,t=q^k,\,q>0$, and $k>-1/h$ for the Coxeter number 
$h=(\rho,\vartheta)+1.$ 

(iii) The transform $f=\sum_b f_b\delta^\#_b\mapsto 
\widetilde{f}= \sum_b f_b^*\delta^\#_b$ is an involution: 
$\widetilde{(\widetilde{f})}=f.$ To apply it for the second 
time, we need to replace 
$\eps_b^*$ by the corresponding $\delta$\~function, which is 
$\sum_c \eps_b^*(\pi_c)\mu^1(\pi_c)\delta^\#_c.$ 
\end{theorem} 

Recall that $\eps^*_b$ becomes the Matsumoto spherical 
function $\phi_b$ in the limit $q\to \infty$ 
upon the substitution $X\mapsto \Lambda.$ 
It is easy 
to calculate the limits of $\mu^0$ and $\mu^1(\pi_b).$ 
We come to a variant of the Macdonald\~Matsumoto formula. 
Claim (iii) has no counterpart in the $p$\~adic theory. 
Technically, it is because the conjugation $\ast$ sends 
$q\mapsto q^{-1}$ and is not compatible with the limit 
$q\to \infty.$ It is equivalent to the non$-p$\~adic 
self-duality $\eps_b(\pi_c)=\eps_c(\pi_b)$ of the 
nonsymmetric Macdonald polynomials. 
The following theorem also has no 
$p$\~adic counterpart because the Gaussian is missing. 

\begin{theorem} \label{tgeps} 
We set $\gamma(\pi_b)= 
q^{(\,\pi_b(\!(k\rho)\!)\,,\,\pi_b(\!(k\rho)\!)\,)/2},$ 
where  $t=q^k,$ use $\nu_\alpha=(\alpha,\alpha)/2,$ 
and $\rho=(1/2)\sum_{\alpha>0}\alpha.$ 
For arbitrary $b,c\in P,$ 
\begin{align}\label{fhatmu} 
&\lr \eps_b^*\, , \eps_c^*\,\gamma\rr_{Del}\ =\ 
\gamma(\pi_0)^2\gamma(\pi_b)^{-1}\gamma(\pi_c)^{-1} 
\eps_c^*(\pi_b)\lr 1\, ,\,\gamma\rr_{Del},\\ 
&\lr 1\, ,\, \gamma\rr_{Del} = 
(\sum_{a\in P} \gamma(\pi_a))\, 
\prod_{\alpha\in R_+}\prod_{ j=1}^{\infty}\Bigl( 
\frac{1- t^{(\rho,\alpha)}q^{j\nu_\alpha}} 
{1-t^{(\rho,\alpha)-1}q^{j\nu_\alpha}}\Bigr). 
\end{align} 
\end{theorem} 

One of the main applications of the double Hecke algebra 
is adding the Gaussian to the classical $p$\~adic theory. 
Technically, one does not  need $\HH\, $ to do this. 
The $\xi$\~deformation of the Iwahori-Matsumoto formulas 
(Theorem \ref{t111}) is the main tool. Its justification 
is elementary. It is surprising 
that it had not been discovered well before the 
double Hecke algebras were introduced. 

\medskip 
\setcounter{equation}{0} 
\section{Degenerate DAHA} \label{SEC:DegDA}
The $p$\~adic origin of the 
double affine Hecke algebra (DAHA) is the most natural 
to consider, however, the connections with real harmonic 
analysis, radial parts and Dunkl operators 
are equally important. They played a key 
role in the beginning of the theory of DAHA; the {\em exact} 
link to the $p$\~adic theory is relatively recent. 
We now discuss 
the trigonometric and rational {\em degenerate DAHAs}, 
which govern the applications in real harmonic analysis. 

This section is a continuation of the paper
\cite{C8}, where we interpreted 
the Harish-Chandra transform 
as a map from the trigonometric-differential 
polynomial representation of the degenerate DAHA 
(in Laurent polynomials) 
to the difference-rational polynomial 
representation. 

{\bf Affine roots.} 
Continuing from the previous section, 
let $R=\{\al\}\subset \R^n$ be a root system of type $A,B,...,F,G$ 
with respect to a euclidean form $(z,z')$ on $\R^n 
\ni z,z'$, 
$W$ the Weyl group  generated by the reflections $s_\al$, 
$R_{+}$ be the set of positive  roots 
corresponding to (fixed) simple roots 
roots $\al_1,...,\al_n,$ 
$\Ga$ the Dynkin diagram 
with $\{\al_i, 1 \le i \le n\}$ as the vertices, 
$R^\vee=\{\al^\vee =2\al/(\al,\al)\}$ the dual root system, 
and 
\begin{align} 
& Q=\oplus^n_{i=1}\Z \al_i \subset P=\oplus^n_{i=1}\Z \om_i, 
\notag \end{align} 
where $\{\om_i\}$ are the fundamental weights defined by 
$ (\om_i,\al_j^\vee)=\de_{ij}$ for the 
simple coroots $\al_i^\vee.$ 

Recall that the form is normalized 
by the condition  $(\al,\al)=2$ for the 
{\em short} roots. This normalization coincides with that 
from the tables in \cite{Bo} for $A,C,D,E,G.$ 
Hence $\nu_\al=(\al,\al)/2$ can be  $1,2,$ or $3.$ 
We write $\nu_{\lng}$ for long roots ($\nu_{\sht}= 1$) 
and also set $nu_i=\nu_{\al_i}.$ 
Let  $\vth\in R^\vee $ be the maximal positive 
{\em coroot} (it is maximal short in $R$) and 
$\rho= (1/2)\sum_{\al\in R_+} \al \ = 
\sum_i \om_i.$


The vectors $\ \tal=[\al,\nu_\al j] \in 
\R^n\times \R \subset \R^{n+1}$ 
for $\al \in R, j \in \Z $ form the 
{affine root system} 
$\tR \supset R$ ($z\in \R^n$ are identified with $ [z,0]$). 
We add $\al_0 \equal [-\vth,1]$ to the simple 
roots. 
The set $\tR$ of positive roots is 
$R_+\cup \{[\al,\nu_\al j],\ \al\in R, \ j > 0\}$. 
Let $\tal^\vee=\tal/\nu_\al,$ so $\al_0^\vee=\al_0.$ 

The Dynkin diagram $\Ga$ of $R$ 
is completed by $\al_0$ (by $-\vth$ to be more 
exact). The notation is $\tGa$. It is the 
completed (affine) Dynkin diagram for $R^\vee$ from \cite{Bo} 
with the arrows reversed. 

The set of 
the indices of the images of $\al_0$ by all 
the automorphisms of $\tGa$ will be denoted by $O$ 
($O=\{0\} \for E_8,F_4,G_2$). Let $O'={r\in O, r\neq 0}$. 
The elements $\om_r$ for $r\in O'$ are the so-called minuscule 
weights: $(\om_r,\al^\vee)\le 1$ for 
$\al \in R_+$.

The {affine Weyl group} $\tW$ is generated by the affine 
reflections $s_{\tal}$. This group is 
the semidirect product $W\lsmash Q$ of 
its subgroups $W$ and the lattice $Q.$ 

The extended Weyl group $ \hW$ is generated by $W$ and $P.$ 
It is isomorphic to $W\lsmash P$ and, also, isomorphic 
to $\Pi \lsmash \tW$ 
for the group $\Pi$ formed by 
the elements of $\hW$ leaving $\tGa$ invariant. 

The latter group is isomorphic to $P/Q$ by the natural 
projection $\{\om_r \mapsto \pi_r,\ r\in O\},$ where 
$\om_r=\pi_r u_r, u_r\in W,\ \pi_0=\,$id. 
The elements 
$\{ u_r\}$ 
preserve the set $\{-\vth,\al_i, i>0\}.$ 

Setting 
$\hw = \pi_r\tw \in \hW,\ \pi_r\in \Pi,\, \tw\in \tW,$ 
the length $l(\hw)$ 
is by definition the length of the reduced decomposition 
$\tw = s_{i_l}...s_{i_2} s_{i_1} $ 
in terms of the simple reflections 
$s_i, 0\le i\le n.$ 

\smallskip 
\subsection{Definition of DAHA} 
By  $m,$ we denote the least natural number 
such that  $(P,P)=(1/m)\Z.$  Thus 
$$m=2 \for D_{2k},\ m=1 \for B_{2k} \and C_{k},$$ 
otherwise $m=|\Pi|$. 

The double affine Hecke algebra depends 
on the parameters 
$q, t_\nu,\, \nu\in \{\nu_\al\}.$ The definition ring is 
$\Q_{q,t}\equal$ 
$\Q[q^{\pm 1/m},t^{\pm 1/2}]$, formed by the 
polynomials in terms of $q^{\pm 1/m}$ and 
$\{t_\nu^{\pm 1/2} \}.$ 
We set 
\begin{align} 
&   t_{\tal} = t_{\al}=t_{\nu_\al},\ t_i = t_{\al_i},\ 
q_{\tal}=q^{\nu_\al},\ q_i=q^{\nu_{\al_i}},\notag\\ 
&\where \tal=[\al,\nu_\al j] \in \tR,\ 0\le i\le n. 
\label{taljx} 
\end{align} 

It will be convenient in many formulas to 
switch to the parameters 
$\{k_\nu\}$ instead of  $\{t_\nu \},$ setting 
$$ 
t_\al=t_\nu=q_\al^{k_\nu} \for \nu=\nu_\al, \and 
\rho_k=(1/2)\sum_{\al>0} k_\al \al. 
$$ 
The notation $k_i=k_{\al_i}$ also will be used. 

For pairwise commutative $X_1,\ldots,X_n,$ 
\begin{align} 
& X_{\tb}\ =\ \prod_{i=1}^nX_i^{l_i} q^{ j} 
\iif \tb=[b,j],\ \hw(X_{\tb})\ =\ X_{\hw(\tb)}. 
\label{Xdex} 
\\ 
&\hbox{where\ } b=\sum_{i=1}^n l_i \om_i\in P,\ j \in 
\frac{1}{ m}\Z,\ \hw\in \hW. 
\notag \end{align} 
We set $(\tilde{b},\tilde{c})=(b,c)$ ignoring the affine extensions 
in the inner product unless $(\tilde{b},\tilde{c}+d)$= 
$(b,c)+j$ is considered. 

Later $Y_{\tb}=Y_b q^{-j}$ will be needed. Note the 
negative sign of $j$. 

We will also use that $\pi_r^{-1}$ is $\pi_{r^*}$ and 
$u_r^{-1}$ is $u_{r^*}$ 
for $r^*\in O\ ,$  $u_r=\pi_r^{-1}\om_r.$ 
The reflection $^*$ is 
induced by the standard 
{\em nonaffine involution} of the Dynkin diagram $\Gamma$.
\index{involution $r\mapsto r^*$} 
\index{reflection $r\mapsto r^*$} 

\begin{definition} 
The  double  affine Hecke algebra $\HH\,$ 
is generated over $ \Q_{ q,t}$ by 
the elements $\{ T_i,\ 0\le i\le n\}$, 
pairwise commutative $\{X_b, \ b\in P\}$ satisfying 
(\ref{Xdex}), 
and the group $\Pi,$ where the following relations are imposed: 

(o)\ \  $ (T_i-t_i^{1/2})(T_i+t_i^{-1/2})\ =\ 
0,\ 0\ \le\ i\ \le\ n$; 

(i)\ \ \ $ T_iT_jT_i...\ =\ T_jT_iT_j...,\ m_{ij}$ 
factors on each side; 

(ii)\ \   $ \pi_rT_i\pi_r^{-1}\ =\ T_j \iif 
\pi_r(\al_i)=\al_j$; 

(iii)\  $T_iX_b T_i\ =\ X_b X_{\al_i}^{-1} \iif 
(b,\al^\vee_i)=1,\ 
0 \le i\le  n$; 

(iv)\ $T_iX_b\ =\ X_b T_i$ if $(b,\al^\vee_i)=0 
\for 0 \le i\le  n$; 

(v)\ \ $\pi_rX_b \pi_r^{-1}\ =\ X_{\pi_r(b)}\ =\ 
X_{ u^{-1}_r(b)} 
 q^{(\om_{r^*},b)},\  r\in O'$. 
\label{doublex} 
\sq 
\end{definition} 

Given $\tw \in \tW, r\in O,\,$ the product 
\begin{align} 
&T_{\pi_r\tw}\equal \pi_r\prod_{k=1}^l T_{i_k},\where 
\tw=\prod_{k=1}^l s_{i_k}, 
l=l(\tw), 
\label{Twx} 
\end{align} 
does not depend on the choice of the reduced decomposition 
(because $\{T\}$ satisfy the same ``braid'' relations 
as $\{s\}$ do). 
Moreover, 
\begin{align} 
&T_{\hv}T_{\hw}\ =\ T_{\hv\hw}\  \hbox{ whenever}\ 
 l(\hv\hw)=l(\hv)+l(\hw) \for 
\hv,\hw \in \hW. \label{TTx} 
\end{align} 
In particular, we arrive at the pairwise 
commutative elements from the previous section: 
\begin{align} 
& Y_{b}\ =\  \prod_{i=1}^nY_i^{l_i} \iif 
b=\sum_{i=1}^n l_i\om_i\in P,\where 
 Y_i\equal T_{\om_i}. 
\label{Ybx} 
\end{align} 
They satisfy the relations 
\begin{align} 
&T^{-1}_iY_b T^{-1}_i\ =\ Y_b Y_{\al_i}^{-1} \iif 
(b,\al^\vee_i)=1, 
\notag\\ 
& T_iY_b\ =\ Y_b T_i \iif (b,\al^\vee_i)=0, 
 \ 1 \le i\le  n. 
\end{align} 
For arbitrary nonzero $q,t,$ any element $H \in \HH\,$ 
has a unique decomposition in the form 
\begin{align} 
&H =\sum_{w\in W }\,  g_{w}\, f_w\, T_w,\ 
g_{w} \in \Q_{q,t}[X],\ f_{w} \in \Q_{q,t}[Y], 
\label{hatdecx} 
\end{align} 
and five more analogous decompositions corresponding 
to the other orderings 
of $\{T,X,Y\}.$ It makes the polynomial representation 
(to be defined next) 
the $\HH\,$\~module induced from the one-dimensional 
representation $T_i\mapsto t_i^{1/2},\,$ $Y_i\mapsto Y_i^{1/2}$ 
of the affine Hecke subalgebra $\h_Y=\lan T,Y\ran.$ 

{\bf Automorphisms.} 
The following maps can be uniquely extended to 
automorphisms of $\HH\,$(see \cite{C4},\cite{C12}): 
\index{involution $\vep$} 
\index{eaaaaa@$\vep$ involution} 
\index{automorphisms $\tau_{\pm}$} 
\index{automorphism $\si$} 
\index{saaaaa@$\si$ automorphism} 
\index{taaaaa@$\tau_{\pm}$ automorphisms} 
\begin{align} 
&\vep:\ X_i \mapsto Y_i,\   Y_i \mapsto X_i,\ 
 T_i \mapsto T_i^{-1}\,(i\ge 1),\, 
t_\nu \mapsto t_\nu^{-1},\, 
 q\mapsto  q^{-1}, 
\label{vepx} 
\\ 
&\tau_+:  X_b \mapsto X_b, \ Y_r \mapsto 
X_rY_r q^{-\frac{(\om_r,\om_r)}{2}},\ 
T_i\mapsto T_i\,(i\ge 1),\ \ t_\nu \mapsto t_\nu,\ 
 q\mapsto  q, 
\notag\\ 
&\tau_+:\ Y_\vth \mapsto q^{-1}\,X_\vth T_0^{-1} 
T_{s_\vth},\, T_0\mapsto q^{-1}\,X_\vth T_0^{-1}, 
\and 
\label{taux}\\ 
&\tau_-\  \equal  \vep\tau_+\vep,\and 
\si\equal \tau_+\tau_-^{-1}\tau_+ = 
\tau_-^{-1}\tau_+\tau_-^{-1}= \vep\si^{-1}\vep, 
\label{tauminax} 
\end{align} 
where $r\in O'.$ 
In the definition of $\tau_\pm$ and $\si,$ 
we need to add $q^{\pm 1/(2m)}$ to 
$\Q_{q,t}.$ 

Here the quadratic relation (o) from Definition 
\ref{doublex} may be omitted. Only the group relations matter. 
Thus these automorphisms act in the corresponding braid 
groups extended by fractional powers of $q$ treated as 
central elements. 
 
The elements $\tau_\pm$ generate the projective $PSL(2,\Z),$ 
which is isomorphic to the braid group $B_3$ due to Steinberg. 
Adding $\vep,$ we obtain the projective $PGL(2,\Z).$ 

These and the statements below are from \cite{C2}. 
Note that $\HH\, ,$ its degenerations below, and the 
corresponding polynomial representations 
are actually defined over $\Z$ extended by the parameters of DAHA. 
We prefer to stick to the field $\Q,$ because the 
Lusztig isomorphisms will require $\Q.$ 

\subsection{Polynomials, intertwiners} 
The Demazure\~Lusztig 
operators are defined as follows: 
\begin{align} 
&T_i\  = \  t_i ^{1/2} s_i\ +\ 
(t_i^{1/2}-t_i^{-1/2})(X_{\al_i}-1)^{-1}(s_i-1), 
\ 0\le i\le n, 
\label{Demazx} 
\end{align} 
and obviously preserve $\Q[q,t^{\pm 1/2}][X]$. 
We note that only the formula for $T_0$ involves $q$: 
\begin{align} 
&T_0\  =  t_0^{1/2}s_0\ +\ (t_0^{1/2}-t_0^{-1/2}) 
( q X_{\vth}^{-1} -1)^{-1}(s_0-1),\notag\\ 
&\where 
s_0(X_b)\ =\ X_bX_{\vth}^{-(b,\vth)} 
 q^{(b,\vth)},\ 
\al_0=[-\vth,1]. 
\end{align} 

The map sending $ T_j$ to the formula in 
(\ref{Demazx}),  $\ X_b \mapsto X_b$ 
(see (\ref{Xdex})), and 
$\pi_r\mapsto \pi_r$ induces a 
$ \Q_{ q,t}$\~linear 
homomorphism from $\HH\, $ to the algebra of linear endomorphisms 
of $\Q_{ q,t}[X]$. 
This $\HH\,$\~module, which will be called the 
{\em polynomial representation}, 
\index{polynomial rep} 
is faithful 
and remains faithful when   $q,t$ take 
any nonzero complex values, assuming that 
$q$ is not a root of unity. 

The images of the $Y_b$ are called the 
{\em difference Dunkl operators}. 
\index{Dunkl operators (difference)} 
To be more exact, 
they must be called difference-trigonometric Dunkl operators, 
because there are also 
difference-rational Dunkl operators. 

\medskip 
{\bf Intertwining operators.} 
The $Y$\~intertwiners (see \cite{C1}) 
are introduced as follows: 
\begin{align} 
&\Psi_i\ =\ 
T_i + (t_i^{1/2}-t_i^{-1/2}) 
(Y_{\al_i}^{-1}-1)^{-1} \for 1\le i\le n, 
\notag\\ 
&\Psi_0\ =\ X_{\vth}T_{s_\vth}-
(t_0^{1/2}-t_0^{-1/2})(Y_{\al_0}-1)^{-1}\notag\\ 
&=Y_0T_0X_0 + (t_0^{1/2}-t_0^{-1/2})(Y_{\al_0}^{-1}-1)^{-1},\for 
Y_0=Y_{\al_0}=q^{-1}Y_{\vth}^{-1},\notag\\ 
& F_i=\Psi_i (\psi_i)^{-1},\ 
\psi_i= t_i^{1/2} + 
(t_i^{1/2} -t_i^{-1/2})(Y_{\al_i}^{-1}-1)^{-1}. 
\label{Phix} 
\end{align} 

These formulas are the $\vep$\~images of the formulas 
for the $X$\~intertwiners, which are a straightforward 
generalization 
of those in the affine Hecke theory. 

The intertwiners belong to $\HH\,$ extended by 
the rational functions in terms of $\{Y\}$. The $F$ are called 
the {\em normalized intertwiners}. 
\index{intertwiners normalized} 
The elements 
$$ 
F_i,\ P_r\equal\ X_r T_{u_r^{-1}},\  0\le i\le n,\ r\in O', 
$$ 
satisfy the same relations 
as $\{s_i,\pi_r\}$ do, so the map 
\begin{align} 
\hw\mapsto F_{\hw}\ =\ P_r F_{i_l}\cdots F_{i_1}, 
\where \hw=\pi_r s_{i_l}\cdots s_{i_1}\in \hW, 
\label{Phiprodx} 
\end{align} 
is a  well defined homomorphism  from $\hW.$ 

The intertwining property is 
$$ 
F_{\hw} Y_b F_{\hw}^{-1}=Y_{\hw(b)} 
\where Y_{[b,j]}\equal Y_b q^{-j}. 
$$ 
The $P_1$ in the case of $GL$ is due to Knop and Sahi. 

As to $\Psi_i$, they 
satisfy the  homogeneous Coxeter relations 
and those with $\Pi_r.$ So we may set 
$\Psi_{\hw} =$ $P_r \Psi_{i_l}\cdots \Psi_{i_1}$ for the reduced 
decompositions. They intertwine $Y$ as well. 

The formulas for $\Psi_i$ when $1\le i\le n$ are well known 
in the theory of affine Hecke algebras. 
The affine intertwiners, those for $\hw\in\hW$, are 
the raising operators for the Macdonald nonsymmetric polynomials, 
serve the Harish-Chandra spherical transform and 
Opdam's nonsymmetric transform, and are a 
key tool in the theory of semisimple representations of DAHAs.

\subsection {Trigonometric degeneration} \label{SECTDEGEN} 
We set $\Q_k\equal\Q[k_\al].$ If the integral coefficients 
are needed, we take $\Z_k\equal\Z[k_\al,1/m]$ as the 
definition ring. 

\index{degenerate DAHA $\HH'$} 
The {\em degenerate (graded) double affine Hecke algebra} 
$\HH'\,$ is the span of 
the group algebra $\Q_k \hW$ and the pairwise commutative 
\begin{align} 
&y_{\tb}\equal 
\sum ^n_{i=1}(b,\al_i^\vee)y_i +u \for 
\tb=[b,u]\in P\times\Z, 
\notag \end{align} 
satisfying  the following relations: 
\begin{align} 
&s_j y_{b}-y_{s_j(b)}s_j\ =\ - k_j(b,\al_j ),\ 
(b,\al_0)\equal -(b,\vth), 
\notag\\ 
&\pi_r y_{\tb}\ =\ y_{\pi_r(\tb)}\pi_r \for 
  0\le j\le n, \  r\in O. 
\label{sukax} 
\end{align} 

\rem 
Without $s_0$ and $\pi_r$, we arrive at the defining relations 
of the graded affine Hecke algebra from  \cite{L}. 
The algebra $\HH'$ has two natural 
polynomial representations via the differential-trigonometric and 
difference-rational Dunkl operators. There is also a third one, 
the representation in terms of infinite differential-trigonometric 
Dunkl operators, which generates (at trivial center charge) 
differential-elliptic $W$\~invariant operators 
generalizing those 
due to Ol\-sha\-net\-sky\~Pe\-re\-lo\-mov. 
See, e.g., \cite{C1}. 
We will need in this section  only 
the (most known) differential-trigonometric polynomial 
representations. 
\sq 

Let us establish a connection with the general DAHA. We set 
$$ 
q=\exp(\mathfrak{v}),\ t_j=q_i^{k_i}=q^{\nu_{\al_i}k_i}, 
\ Y_b=\exp(-\mathfrak{v} y_b),\ \mathfrak{v}\in \C. 
$$ 
Using $\vep$ from (\ref{vepx}), 
the algebra $\HH\,$ is generated by 
$Y_b,\,$ $T_i$ for $1\le i\le n,$ and 
$$ 
\vep(T_0)=X_{\vth}T_{s_{\vth}},\ 
\vep(\pi_r)=X_rT_{u_r^{-1}},\ r\in O'. 
$$ 
It is straightforward to see that the relations 
(\ref{sukax}) for the $y_b,\,s_i (i>0),\,$ 
$s_0,$ and the $\pi_r$ are the leading coefficients 
of the $\mathfrak{v}$\~expansions of the general relations 
for this system of generators. 
Thus $\HH'$ is $\HH\,$ in the limit $\mathfrak{v}\to 0.$ 

When calculating the limits of the 
$Y_b$ in the polynomial representation, 
the ``trigonometric" derivatives of $\Q[X]$ appear: 
$$ 
\partial_{a}(X_{b})\ =\ (a ,b)X_{b},\ a,b\in P,\ \, 
w(\partial _{b})= 
\partial_{w(b)}, \  w\in W. 
$$ 

The limits of the formulas 
for $Y_b$ acting in the polynomial representation 
are the {\em trigonometric Dunkl operators} 
\index{Dunkl operators (trig)} 
\begin{align} 
\d_{b} \equal\ 
&\partial_b + \sum_{\al\in R_+} \frac{ k_{\al}(b,\al)}{ 
(1-X_{\al}^{-1}) } 
\bigl( 1-s_{\al} \bigr)- (\rho_{k},b). 
\label{dunkx} 
\end{align} 
They act on the Laurent polynomials $f\in \Q_k[X],$ 
are pairwise commutative, and $y_{[b,u]}=\d_b+u$ 
satisfy (\ref{sukax}) for 
the following action of the group $\hW$: 
$$ w^x(f)=w(f)\for w\in W,\ b^x(f)=X_b f \for b\in P.$$ 
For instance, $s_0^x(f)=X_\vth s_\vth(f)$, and 
$\pi_r^x(f)=X_r u_r^{-1}(f).$ 

Degenerating $\{ \Psi\}$, one obtains the intertwiners of 
$\HH'\,$: 
\begin{align} 
&\Psi'_i= 
s_i + \frac{\nu_i k_i }{ y_{\al_i} },\ 0\le i\le n,\ 
\bigl(\,\Psi'_0= 
X_\vth s_\vth + \frac{k_0 }{ 1-y_{\vth}} \hbox{\ in\ } 
\Q_k[X]\,\bigr),\notag\\ 
&P_r' = \pi_r,\ \bigl(\,P_r'=X_r u_r^{-1} \hbox{\ in\ } 
\Q_k[X]\, \bigr),\ r\in O'. 
\label{Phiprimex} 
\end{align} 
The operator $P_1'$ in the case of $GL$ (it is of 
infinite order) plays the key role in \cite{KnS}.

Recall that the general normalized intertwiners are 
$$ 
F_i=\Psi_i\psi_i^{-1}, \ 
\psi_i=t^{1/2}+(t_i^{1/2}-t^{1/2})(Y_{\al_i}^{-1}-1)^{-1}. 
$$ 
Their limits are 
$$F'_i=\Psi'_i(\psi'_i)^{-1},\ 
\psi'_i=1+\frac{\nu_i k_i }{ y_{\al_i} }.$$ 
They satisfy the unitarity condition 
$(F_i')^2=1,$ and the products $F'_{\hw}$ can be defined for any 
decompositions of $\hw.$ 
One then has: 
$$F'_{\hw}\, y_b\, (F'_{\hw})^{-1}\ =\ y_{\hw(b)}.$$ 

Equating 
$$F_i=F'_i \for 0\le i\le n,\ P_r=P'_r \for r\in O, 
$$ 
we come to the formulas for $T_i\, (0\le i\le n),\, $ 
$X_r\, (r\in O')$ in terms of 
$s_i, y_b,$ and $Y_b=\exp(-\mathfrak{v}y_b).$ 

These formulas determine 
the {\em Lusztig homomorphism} 
\index{Lusztig map} 
$\kapp'$ from $\HH\, $ to 
the completion $\Z_{k,q,t}\HH'[[\mathfrak{v}y_b]]$ for 
$\Z_{k,q,t}\equal\Z_k\Z_{q,t}.$ 
See, e.g., \cite{C1}. 

For instance, 
$X_r\in \HH$ becomes $\pi_r T_{u_r^{-1}}^{-1}$ in $\HH',$ 
where the $T$\~factor has to be further expressed in 
terms of $s,y.$ 
In the degenerate polynomial representation, $\kapp'(X_r)$ acts 
as $X_r (\kapp'(T_{u_r^{-1}})u_r)^{-1},$ not as 
straightforward multiplication 
by $X_r.$  However, these two actions coincide 
in the limit $\mathfrak{v}\to 0,$ since $T_w$ become $w.$ 

Upon the $\mathfrak{v}$\~completion, 
we obtain an isomorphism 
$$ 
\kapp':\Q_k[[\mathfrak{v}]]\otimes\HH\to 
\Q_k[[\mathfrak{v}]]\otimes\HH'. 
$$ 

We will use the 
notation $(d,[\al,j])=j.$ For instance, 
$(b+d,\al_0)=1-(b,\vth).$ 

Treating $\mathfrak{v}$ 
as a nonzero number, 
an arbitrary 
$\HH'$\~module $V'$ that is a union 
of finite dimensional $Y$\~modules has a natural structure 
of an $\HH$\~module provided that we have 
\begin{align} 
& q^{(\al_i,\xi+d)}=t_i \Rightarrow 
(\al_i,\xi+d) = \nu_i k_i,\label{lusisox} 
\\ 
& q^{(\al_i,\xi+d)}=\ 1 \Rightarrow 
(\al_i,\xi+d) = 0, \where\notag \\ 
& 0\le i\le n,\ 
y_b(v')=(b,\xi)v' \for \xi\in \C^n,\ 0\neq v'\in V'.\notag 
\end{align} 

For the modules of this type, 
the map $\kapp'$ is over the ring $\Q_{k,q,t}$ 
extended by $(\al,\xi+d),\ $ $q^{(\al,\xi+d)}$ for $\al\in R,$ 
and $y$\~eigenvalues $\xi.$  Moreover, we need to localize 
by $(1- q^{(\al,\xi+d)})\neq 0$ 
and by $(\al,\xi+d)\neq 0.$ 
Upon such extension and localization, $\kapp'$  is defined over 
$\Z_{k,q,t}$ 
if the module is $y$\~semisimple. If there are nontrivial 
Jordan blocks, then the formulas will contain factorials in the 
denominators. 

For instance, let $\i'[\xi]$ be the $\HH'$\~module 
induced from the 
one-dimensional $y$\~module $y_b(v)=(b,\xi)v.$ Assuming that 
$q$ is not a root of unity, the mapping $\kapp'$ 
supplies it with a structure of $\HH$\~module if 
$$ 
q^{(\al,\xi)+\nu_\al j}=t_\al \hbox{\ implies\ } 
(\al,\xi)+\nu_\al j=\nu_\al k_\al 
$$ 
for every $\al\in R, j\in \Z,$ 
and the corresponding implications hold for $t$ replaced by $1.$ 
This means that 
\begin{align} 
&(\al,\xi)-\nu_\al k_\al\, \not\in \, 
\nu_\al\Z\, +\, \frac{2\pi\imath}{\mathfrak{v}} 
(\Z\setminus \{0\})\, \not\ni\, (\al,\xi)
\hbox{\ for\ all\ }\al\in R. 
\label{intercond} 
\end{align} 
Generalizing, we obtain that $\kapp'$ is well defined for any 
$\HH'$\~module generated by its 
$y$\~eigenvectors with the $y$\~eigenvalues $\xi$ satisfying this 
condition, assuming that $\mathfrak{v}\not\in \pi\imath\Q.$ 

\rem 
Actually, there are at least four different variants of $\kapp'$ 
because the 
normalization factors $\psi,\psi'$ may be associated with different 
one-dimensional characters 
of the affine Hecke algebra $\lan T,Y\ran$ 
and its degeneration. It is also 
possible to multiply the normalized intertwiners 
by the characters of $\hW$ 
before equating. Note that if we 
divide the intertwiners $\Psi$ and/or $\Psi'$ 
by $\psi,\psi'$ 
on the left in the definition of $F,F',$ it corresponds to 
switching from $T_i\mapsto t_i$ to 
the character $T_i \mapsto -t_i^{-1/2}$ together with 
the multiplication by the sign-character of $\tW.$ 
\sq 


\vskip 0.2cm 
\subsection{Rational degeneration} 
The limiting procedure to the {\em 
rational Dunkl operators} 
\index{Dunkl operators (rational)} 
is as follows. 
We set $X_b=e^{\mathfrak{w} x_b},\ d_b(x_c)=(b,c),$ 
so the above derivatives 
$\partial_b$ become $\partial_b=(1/\mathfrak{w})d_b.$ 
In the limit $\mathfrak{w}\to 0,$ 
$\mathfrak{w}\d_b$ tends to the operators 
\begin{align} 
D_{b} \equal 
&d_b + 
\sum_{\al\in R_+} \frac{ k_{\al}(b,\al)}{ 
x_\al } 
\bigl( 1-s_{\al} \bigr), 
\label{dunkorigx} 
\end{align} 
which were introduced by Dunkl and, as a matter 
of fact, gave birth to the DAHA direction, 
although KZ and the Macdonald theory must be mentioned 
too as origins of DAHAs. 

These operators are pairwise commutative and satisfy 
the cross-relations 
\begin{align} 
&D_bx_c-x_cD_b=(b,c)+\sum_{\al>0} k_\al (b,\al)(c,\al^\vee)s_\al, 
\for b,c\in P. 
\label{duncrossx} 
\end{align} 
These relations, 
the commutativity of $D,$ the commutativity of $x,$ and the 
$W$\~equivariance 
$$ w\, x_b\, w^{-1}\,=\,x_{w(b)},\ w\, D_b\, w^{-1}\,=\,D_b\ 
\for b\in P_+,\, w\in W, 
$$ 
are the defining relations of the {\em rational DAHA} $\HH''.$ 
\index{double Hecke algebra (rational)} 
\index{degenerate DAHA $\HH''$} 
\index{Haaa@$\HH',\HH''$} 
The references are \cite{CM} (the case of $A_1$) and 
\cite{EG}; however, the key part of the definition is the 
commutativity of $D_b$ due to Dunkl \cite{Du}. 
The Dunkl operators and the operators of multiplication 
by the $x_b$ 
form the 
{\em polynomial representation} of $\HH'',$ which is faithful. 
\index{polynomial rep (rational)} 
It readily justifies the PBW theorem for $\HH''.$ 

Note that in contrast to the $q,t$\~setting, the definition 
of the rational DAHA can 
be extended to finite groups generated by complex reflections 
(Dunkl, Opdam and Malle). 
There is also a generalization due to Etingof\~Ginzburg 
from \cite{EG} (the symplectic reflection algebras). 

\rem 
Following \cite{CO}, there is a one-step limiting procedure from 
$\HH $ to $\HH''.$ 
We set 
$$ 
Y_b=\exp(-\sqrt{\mathfrak{u}} D_b),\ 
X_b=e^{\sqrt{\mathfrak{u}} x_b}, 
$$ 
assuming that 
$q=e^{\mathfrak{u}}$ and $\mathfrak{u}\to 0.$ 
We come directly to the 
relations of the rational DAHA and the formulas for $D_b.$ 
The advantage of this direct 
construction is that the automorphisms $\tau_{\pm}$ 
survive in the limit. 
Indeed, $\tau_+$ in $\HH\, $ 
can be interpreted as the formal conjugation by 
the $q$\~Gaussian $q^{x^2/2},$ where $x^2=\sum_i x_{\om_i} 
x_{\al_i^\vee}.$ 
In the limit, it becomes the conjugation by $e^{x^2/2},$ 
preserving 
$w\in W,\, x_b,\, $ and taking $D_b$ to $D_b-x_b.$ 
Respectively, $\tau_-$ 
preserves $w$ and $D_b,$ and sends $x_b\mapsto x_b-D_b.$ 
These automorphisms 
do not exist in the $\HH'.$ 
\sq 
 
The {\em abstract Lusztig-type map} 
\index{Lusztig map} 
from $\HH'$ to $\HH''$ is as follows. 
Let $w\mapsto w$ and 
$X_b=e^{\mathfrak{w} x_b}.$ 
We expand $X_\al$ in terms 
of $x_\al$ in the formulas for the trigonometric 
Dunkl operators $\d_b$: 
\begin{align} 
&\d_b=\frac{1}{\mathfrak{w}} 
D_b -(\rho_k,b)+ \sum_{\al\in R_+}k_\al(b,\al) 
\sum_m^\infty \frac{B_m}{m!}(-\mathfrak{w} x_\al)^m\, (1-s_\al) 
\label{lusratx} 
\end{align} 
for the Bernoulli numbers $B_m.$ 
Then we use these formulas as abstract expressions for $y_b$ 
in terms of the generators of $\HH''$: 
$y_b= \frac{1}{\mathfrak{w}} D_b+\ldots\ .$ 

One obtains an isomorphism 
$\kapp'':\Q[[\mathfrak{w}]]\otimes\HH' 
\to \Q[[\mathfrak{w}]]\otimes\HH'',$ 
which maps $\HH'$ to the 
extension of $ \HH''$  by the formal series in terms of 
$\mathfrak{w} x_b.$ An arbitrary 
representation $V''$ of $\HH''\,$ that is a union of 
finite dimensional 
$\Q_k[x]$\~modules becomes an $\HH'$\~module provided that 
\begin{align} 
&\mathfrak{w}\ze_\al\not\in 2\pi\imath(\Z\setminus \{0\})\for 
x_b(v)=(\ze,b)v, 0\neq v\in V''. 
\label{intercondd} 
\end{align} 

Similar to (\ref{intercond}), 
this constraint simply restricts choosing $\mathfrak{w}\neq 0.$ 
The formulas for $y_b$ become locally finite 
in any representations of $\HH''$, where $x_b$ act 
locally nilpotent, 
for instance, in finite dimensional $\H''$\~modules. 
In this case, there are no restrictions 
for $\mathfrak{w}.$ 

\rmk 
Note that the ``identification" of $\HH'$ and $\HH''$ 
has a common source with the method used in the so-called 
localization due to Opdam and Rouquier (see \cite{GGOR}, 
\cite{VV1}). They separate the differentials $d_b$ from 
the formula for the Dunkl operators to define a 
KZ-type connection with values in double Hecke modules. 
We equate these differentials in the rational and 
trigonometric formulas for the Dunkl operators 
to connect $\HH'$ and $\HH''.$ 

\vskip -0.3cm \sq 

Finally, the composition 
$$ 
\kapp\equal\kapp''\circ\kapp': 
\HH[[\mathfrak{v},\mathfrak{w}]]\to 
\HH''[[\mathfrak{v},\mathfrak{w}]] 
$$ 
is an isomorphism. Without the completion, 
it makes an arbitrary finite dimensional 
$\HH''$\~module $V''$ a 
module over $\HH$ as $q=e^{\mathfrak{v}}, t_\al=q^{k_\al}$ 
for sufficiently general (complex) nonzero numbers 
$\mathfrak{v},\mathfrak{w}.$ 
This isomorphism was discussed in \cite{BEG} (Proposition 7.1). 

The finite dimensional representations are the most 
natural here because, on the one hand, 
$\kapp''$ lifts the modules that are unions of 
finite dimensional $x$\~modules to those for $X$, on 
the other hand, $\kapp'$ maps 
the $\HH'\,$\~modules that are 
unions of finite dimensional $y$\~modules to those for $Y.$ 
So one must impose these conditions for both $x$ and $y.$ 

We obtain a functor from 
the category of finite dimensional representations 
of $\HH''$ to that for $\HH\,.$ Indeed, 
it is known now that there are finitely 
many irreducible objects in the former category. 
Therefore we can find $\mathfrak{v},\mathfrak{w}$ serving 
{\em all} irreducible representations and their extensions. 
This functor is {\em faithful} 
provided that $q$ is generic in the following 
sense: 
$$q^at^b=1 \Rightarrow a+kb=0 \for a,b\in\BQ. 
$$ 
This condition  ensures the {\em equivalence 
of the categories} of finite dimensional representations 
for $\HH\, $ and $\HH''.$ 

Using $\kapp$ for infinite dimensional representations 
is an interesting problem. It makes the theory 
analytic. For instance, the triple composition 
$\kapp''\circ \g \circ \kapp'$ for the inverse 
Opdam transform $\g$ 
(see \cite{O2} and formula (6.1) from \cite{C8}) 
embeds $\HH\, $ in 
$\HH''$ and 
identifies the $\HH''$\~module $\C_c^\infty(\R^n)$ with 
the $\HH$\~module of PW-functions under the condition 
$\Re k>-1/h.$ 
See \cite{O2,C8} for more detail. 
\smallskip 

The degenerations above play the role of the Lie algebras 
in the theory of DAHA. The Lusztig isomorphisms are certain 
counterparts of the $\exp$-$\log$ maps. 
It is especially true for the composition map from 
$\HH\,$  to the 
rational degeneration, because the latter has the projective action 
of the $PSL(2,\BZ)$ and some other features that make it 
close to the general $\HH.$  The algebra $\HH'$ is 
not projective $PSL(2,\BZ)$\~invariant. 

Note that there are some special ``rational" 
symmetries and tools that 
have no $q,t$\~counterparts. On the other hand, 
the semisimple representations have no immediate analogs 
in the rational theory and the perfect 
representations are simpler to deal with in the $q,t$\~case. 
Thus it really resembles the relation between Lie groups 
and Lie algebras. 

\subsection{Diagonal coinvariants}\label{sec:DIAGCOI} 
It was conjectured by Haiman \cite{Hai} that the space of 
{\em diagonal coinvariants} 
\index{diagonal coinvariants} for a root system $R$ of rank $n$ 
has a ``natural" quotient of dimension $(1+h)^n$ for 
the Coxeter number $h$. This space 
is the quotient $\C[x,y]/(\C[x,y]\C[x,y]^W_o)$ 
for the algebra of polynomials $\C[x,y]$ with the diagonal action 
of the Weyl group on $x\in \C^n \ni y$ and the ideal 
$\C[x,y]^W_o \subset \C[x,y]^W$ of the $W$\~invariant 
polynomials without the constant term. 

In \cite{Go}, such a quotient was constructed. 
It coincides with the whole space of the diagonal coinvariants 
in the $A_n$\~case due to Haiman, but this does not hold 
for other root systems. 

Let $k_{\sht}=-(1+1/h)=k_{\lng}, $ 
$h$ be the Coxeter number $1+(\rho,\vth).$ 
The polynomial representation $\Q[x]$ of $ \HH''$ 
has a unique nonzero irreducible quotient-module. 
It is of dimension $(1+h)^n,$ which was checked in \cite{BEG} 
and \cite{Go}, and also 
follows from \cite{C12} via $\kapp.$ 

The application of this representation to 
the coinvariants of the ring of commutative 
polynomials $\Q[x,y]$ with the diagonal action of $W$ 
is as follows. 

\smallskip 
The polynomial representation $\Q[x]$ is naturally a quotient of 
the {\em linear space} $\Q[x,y]$ 
considered as an induced $\HH''$\~module from the one-dimensional 
$W$\~module $w\mapsto 1.$ So is $V''.$ The subalgebra 
$(\HH'')^W$ of 
the $W$\~invariant elements from $\HH''$ preserves the image 
of $\Q\de$ in $V''$ 
for $\de\equal\prod_{\al>0} x_\al.$ 

Let $I_o\subset (\H'')^W$ be the ideal of the elements vanishing 
at the image of $\de$ in $V''.$ Gordon proves that 
{\em $V''$ coincides 
with the quotient $\tV''$ of $\HH''(\de)$ by the $\HH''$\~submodule 
$\HH'' I_o(\de).$} It is sufficient to check that $\tV''$ is 
irreducible. 

The graded space gr$(V'')$ 
of $V''$ with 
respect to the total $x,y$\~degree of the polynomials 
is isomorphic as a linear space to the quotient of 
$\Q[x,y]\de$ by the 
graded image of $\HH''I_o(\de).$  The latter 
contains $\Q[x,y]^W_o\de$ for the ideal 
$\Q[x,y]^W_o\subset \Q[x,y]^W$ 
of the $W$\~invariant polynomials without the constant term. 
Therefore $V''$  becomes a certain quotient of 
$\Q[x,y]/(\Q[x,y]\Q[x,y]_o^W).$ 
See \cite{Go},\cite{Hai} about the connection with 
the Haiman theorem in the $A_n$\~case and related questions for 
other root systems. 

The irreducibility of the $\tV''$ above is the key fact. The proof 
from \cite{Go} requires considering a KZ-type local systems 
and the technique from \cite{GGOR}. 
It was demonstrated in  \cite{C29} that the irreducibility 
can be readily proved in the $q,t$\~case by using the passage to 
the roots of unity and therefore gives an entirely algebraic 
and simple proof of Gordon's theorem via the $\kapp$\~isomorphism. 
\index{Paa@$\p$:\ polynomial rep|)} 

\setcounter{equation}{0} 
\section{Harish-Chandra inversion} \label{sec:harish} 
In this section, we use the 
trigonometric degenerate double Hecke 
algebra and Dunkl operators  to calculate the images of 
the operators of multiplication by symmetric 
Laurent polynomials with respect 
to the Harish-Chandra transform. 
These formulas were known only in some cases, although 
they are directly connected with the important problem of 
making the convolution on the symmetric 
spaces as explicit as possible. 
We first solve a more general problem of 
calculating the transforms of the coordinates 
for the nonsymmetric Harish-Chandra transform 
defined by Opdam; this problem has a complete solution. 
Then we employ 
the symmetrization. 
The resulting formulas readily lead to a new simple proof 
of the Harish-Chandra inversion theorem (see \cite{HC,He}) 
and the corresponding theorem from \cite{O2}. 

We will assume that $k> 0$ and restrict ourselves to compactly 
supported functions. In this case, we can 
borrow the growth estimates 
from \cite{O2}. Hopefully  this approach can be extended  to any $k$ 
and to the elliptic radial parts, which were discussed in the 
previous section. It is also expected to be helpful 
for developing the analytic theory of the direct and inverse 
transforms. 

For a maximal real split torus $ A$ of a 
semisimple Lie group $G,$ 
the \index{Harish-Chandra transform} 
{\em Harish-Chandra transform} 
is the integration of 
symmetric ($W$\~invariant) functions in terms of $X\in A$ 
multiplied by the spherical zonal function $\phi(X,\la),$ 
where $\la \in (\hbox{Lie}A\otimes_{\R}\C)^*.$ 
The measure is the restriction of the invariant measure on $G$ 
to the space of double cosets $K\backslash G/K \subset A/W$ for 
the maximal compact subgroup $K\subset G$ and 
the restricted  Weyl group $W$. 
The function $\phi,$ the kernel of the transform, 
is a $W$\~invariant eigenfunction 
of the radial parts of the $G$\~invariant differential operators 
on $G/K$; $\la$ determines the set of eigenvalues. 
The parameter $k$ is given by the root multiplicities 
($k=1$ in the group case). 

There is a generalization to arbitrary $k$ 
due to Calogero, Sutherland, Ko\-orn\-win\-der, 
Mo\-ser, Olsha\-net\-sky, 
Pe\-re\-lo\-mov, 
Heck\-man, and Op\-dam; see the papers 
\cite{HO1,H,O1} for a systematic 
(``trigonometric") theory. 
In the nonsymmetric variant due to Opdam \cite{O2}, the Dunkl-type 
operators from 
\cite{Ch4} replace the radial parts of $G$\~invariant operators 
and their 
$k$\~generalizations. 

The analytic part of the theory 
is in extending the direct and inverse transforms 
to various classes of functions. Not much is known here. Hopefully 
the difference theory will be more promising analytically. 

\smallskip 
In the papers \cite{C12,C5}, a difference counterpart of the 
Harish-Chandra 
transform was introduced, which is also a deformation of the 
Fourier transform in the $p$\~adic theory of spherical functions. 
Its  kernel (a $q$\~generalization of $\phi$) 
is defined as  an eigenfunction 
of the $q$\~difference ``radial parts" 
(Macdonald's ``minuscule-weight" operators and their 
generalizations). There are ``algebraic" applications 
in combinatorics (the Macdonald polynomials), 
representation theory 
(for instance, for quantum groups at roots of unity), 
and mathematical physics. This section is toward 
analytic applications. 
\smallskip 

The difference Fourier transform  is  self-dual, i.e., its kernel 
is $x\leftrightarrow \la$ symmetric for $X=q^x.$ 
This holds in the differential theory only for either 
the so-called rational degeneration with the tangent space 
$T_e(G/K)$ instead of $G/K$ (see \cite{He,Du2,J}) or for 
a special $k=0,1$. 
For such a special $k$, the classical differential and 
the new difference transforms coincide up to a normalization. 

The differential case corresponds to the limit  $q\to 1.$ 
At the moment, the analytic methods are not mature enough to 
manage the limiting procedure in detail, which would give 
a complete solution of the inversion problem. So 
we develop the corresponding technique  without any reference 
to the difference Fourier transform 
and the general double affine Hecke algebra. 
We use only the theory of the degenerate one, which 
generalizes Lusztig's graded affine Hecke algebra \cite{L} 
(the $GL_n$\~case is due to Drinfeld). 
Mainly we need the intertwining operators 
from \cite{C1} (see \cite{C15} and \cite{KnS} for $GL_n$). 

\smallskip 
The key result of the section is that 
the Opdam transforms  of the operators of 
multiplication by the coordinates 
coincide with the operators from 
 \cite{Ch5} (see (\ref{Desb}) below). 
Respectively, the Harish-Chandra transform sends the operators of 
multiplications by $W$\~invariant Laurent polynomials 
to the {\em difference} operators from (\ref{Laop}). 

It is not very surprising 
that the images of these important operators 
haven't been found before. 
The calculation involves  the following ingredients, 
which are new in the 
harmonic analysis on symmetric spaces: 

(a) the differential Dunkl-type operators from \cite{Ch4,C13}; 

(b) the Opdam transform \cite{O2}; 

(c) the difference Dunkl operators \cite{Ch5,C15}. 


\subsection{Affine Weyl groups} 
Recall the definition of the affine roots and the extended affine 
Weyl group from the previous section.
Let 
$R=\{\al\}   \subset \R^n$ be the root system 
of type $A,B,...,F,G,$ 
normalized by the standard condition 
$(\al,\al)=2$ for {\em short} $\al;$ 
$\al_1,...,\al_n$ simple roots; 
$a_1=\al^\vee_1,...,a_n=\al^\vee_n$ simple coroots, where 
$\al^\vee =2\al/(\al,\al)=\al/\nu_\al$ for 
$\nu_\al\equal(\al,\al)/2$; 
$\om_1,...,\om_n$ the fundamental weights 
determined by the relations  $ (\om_i,\al_j^\vee)= 
\de_i^j $ for the Kronecker delta. 
We will also use 
\begin{align} 
&Q=\oplus^n_{i=1}\Z \al_i \subset P=\oplus^n_{i=1}\Z \om_i,\ 
P_+=\oplus^n_{i=1}\Z_+ \om_i\for \Z_+=\{m\ge 0\}. 
\notag \end{align} 
 
The vectors $\ \tal=[\al,\nu_\al j] \in 
\R^n\times \R = \R^{n+1}$ 
for $\al \in R, j \in \Z $ 
form the \index{affine root system} {\em affine root system} 
$\tR \supset R$ 
($z\in \R^n$ are identified with $ [z,0]$). 
We add  $\al_0 \equal [-\vth,1]$ to the  simple roots 
for the \index{maximal short root $\vth$} 
{\em maximal short root} 
$\vth \in R$. 
The corresponding set $\tR_+$ of positive roots coincides 
with $R_+\cup \{[\al,\nu_\al j],\  \al\in R, \  j > 0\}$. 


The set of indices of the orbit of the zero vertex 
in the affine Dynkin diagram by its automorphisms 
will be denoted by $O$ ($O=\{0\} 
\for E_8,F_4,G_2$). Let $O'=\{r\in O, r\neq 0\}$. 
We identify the indices with the corresponding simple affine 
roots. 
The elements $\om_r$ for $r\in O'$ are the so-called minuscule 
weights. 

Given $\tal=[\al,\nu_\al j]\in \tR,  \ b \in P$, let 
\begin{align} 
&s_{\tal}(\tz)\ =\  \tz-(z,\al^\vee)\tal,\ 
\ b'(\tz)\ =\ [z,\ze-(z,b)] 
\label{saction} 
\end{align} 
for $\tz=[z,\ze] \in \R^{n+1}$. 

The 
{\em affine Weyl group} $\tW$ is generated by all $s_{\tal}$ 
(simple reflections $s_i=s_{\al_i}$ for $0 \le i \le n$ 
are enough). 
It is 
the semidirect product $W\lsmash Q^\vee$, where the nonaffine Weyl 
group $W$ is the span of $s_\al, 
\al \in R_+$. 
We will identify $b\in P$ with the corresponding 
translations. For instance, 
\begin{align} 
& \al' =\ s_{\al}s_{[\al,\nu_\al]}=\ s_{[-\al,\nu_\al]}s_{\al} \for 
\al\in R. 
\label{baction} 
\end{align}

The 
{\em extended Weyl group} $ \widehat{W}$ generated by 
$W\and P$ 
 is isomorphic to $W\lsmash P$: 
\begin{align} 
&(wb)([z,\ze])\ =\ [w(z),\ze-(z,b)] \for w\in W, b\in P. 
\notag \end{align} 
For $b_+\in P_+$, let 
\begin{align} 
&u_{b_+} = w_0w^+_0  \in  W,\ \pi_{b_+} = 
b_+(u_{b_+})^{-1} 
\ \in \ \widehat{W}, \ u_i=u_{\om_i},\pi_i=\pi_{\om_i}, 
\label{wo} 
\end{align} 
where $w_0$ (respectively, $w^+_0$) is the longest element in $W$ 
(respectively, in $ W_{b_+}$ generated by $s_i$ preserving $b_+$) 
relative to the 
set of generators $\{s_i\}$ for $i >0$. 

The elements $\pi_r\equal\pi_{\om_r}, 
r \in O',$ and $\pi_0=\hbox{id}$ 
leave $\{\al_i,i\ge 0\}$ invariant 
and form a group denoted by $\Pi$, 
 which is isomorphic to $P/Q$ by the natural 
projection $\om_r \mapsto \pi_r$. As with $\{u_r\}$, 
they preserve the set $\{-\vth,\al_i, i>0\}$. 
The relations $\pi_r(\al_0)= \al_r= (u_r)^{-1}(-\vth) 
$ distinguish the 
indices $r \in O'$. Moreover (see, e.g., \cite{Ch5}): 
\begin{align} 
& \widehat{W}  = \Pi \lsmash \tW, \where 
  \pi_rs_i\pi_r^{-1}  =  s_j \iif \pi_r(\al_i)=\al_j,\ 
0\le i,j\le n. 
\notag \end{align} 
The \index{length $l(\hw)$} {\em length} $l(\hw)$  of 
$\hw = \pi_r\tw \in \widehat{W}$ is the length of 
any of the reduced 
decompositions of $\tw\in \tW$ with respect to 
$\{s_i, 0\le i\le n\}$: 
\begin{align} 
l(\hw)=|\la(\hw)|\for \la(\hw)& \equal  \tR_+ \cap 
\hw^{-1}(-\tR_+) 
 \notag\\ 
& =\ \{\tal\in \tR_+, \ l( \hw s_{\tal}) < l(\hw) \}. 
\label{lambda} 
\end{align} 
For instance, given $b_+\in P_+$, 
\begin{align} 
&\la(b_+)\  =\ \{\tal= [\al,\nu_\al j],\ \ \al\in R_+\and 
( b_+, \al )>j\ge 0 \},\notag\\ 
& l(b_+)\ =\ 2(b_+,\rho), 
\where \rho=(1/2)\sum_{\al\in R_+}\al. 
\label{lambi} 
\end{align} 
We will also need  the dominant 
\index{affine Weyl chamber@affine Weyl chamber $\CC^a$} 
{\em affine Weyl chamber} 
\begin{align} 
&\CC^a\ =\ \{ z\in \R^n \hbox{ \ such\ that\  } 
(z,\al_i)>0 \for i>0,\ 
(z,\th)<1\}. 
\label{afchamber} 
\end{align} 

\subsection{Differential case} 
Let us fix a collection $k=\{k_\al\}\subset \C$ such that 
$k_{w(\al)}=k_\al$ for $w\in W,$ i.e., $k_\al$ depends 
only on the $(\al,\al).$ 
We set $k_i=k_{\al_i}$ and $k_0=k_\theta.$  Let 
$ \rho_k= (1/2)\sum _{\al\in R_+}k_\al \al.$ 
One then has: $(\rho_k,\al_i^\vee)=k_i$ for $1\le i\le n.$ 

It will be convenient in this section to switch
to the parameters 
$$\ka_\al\equal \nu_\al k_\al\for \nu_\al\equal(\al,\al)/2;\  
\ka_i=k_{\al_i},\, \ka_0=k_0.$$
For instance, 
$ \rho_k= (1/2)\sum _{\al\in R_+}\ka_\al \al^\vee$ 
and $(\rho_k,\al_i)=\ka_i$ for $1\le i\le n.$ 

\begin{definition} 
\index{degenerate DAHA $\HH'$} 
\index{Haaa@$\HH',\HH''$} 
\index{double Hecke algebra} \index{DAHA} 
The {\em degenerate double affine Hecke algebra} 
$\HH' $ is 
formed by the group algebra $\C \widehat{W}$ and the pairwise 
commutative $y_b$ for $b\in P $ or in $\R^n$ 
satisfying  the following relations: 
\begin{align} 
&s_i y_{b}-y_{s_i(b)}s_i =- \ka_i(b,\al^\vee_i ) 
,\ \pi_r y_{\tb}\ =\ y_{\pi_r(\tb)}\pi_r, 
\notag\\ 
&\hbox{for\ } (b,\al_0)=-(b,\th),\ y_{[b,u]}\equal y_b+u,\ 
0\le i\le n, \  r\in O. 
\label{suka} 
\end{align} 
\label{DEGHE} \sq
\end{definition} 

Without $s_0$ and $\pi_r$, these are the defining relations 
of the graded affine Hecke algebra from  \cite{L}. 
This algebra is a degeneration as $q\to 1,\ t=q^{\ka}$ 
of the general double affine 
Hecke algebra. 

We will need the 
\index{intertwiners aa@intertwiners $\Psi'_{\hw}$ (degenerate)} 
{\em intertwiners} of 
$\HH'\,$: 
\begin{align} 
&\Psi'_i\ =\ 
s_i + \frac{\ka_i }{ y_{\al_i} },\ 
\Psi'_0\ =\ 
X_\th s_\th + \frac{\ka_0 }{ 1-y_{\th} },\notag\\ 
&P_r'\ =\ X_r u_r^{-1}, \for  1\le i\le n,\ r\in O'. 
\label{Phiprime} 
\end{align} 
The elements 
$\{\Psi'_i,\ P'_r\}$ satisfy the homogeneous Coxeter relations 
for $\{s_i, \pi_r\},$ so the elements 
\begin{align} 
\Psi'_{\hw}\ =\ P'_r \Psi'_{i_l}\cdots \Psi'_{i_1}, 
\where \hw=\pi_r s_{i_l}\cdots s_{i_1}\in \widehat{W}, 
\label{Phiprod} 
\end{align} 
are well defined 
for reduced decompositions of $\hw$, 
and $\Psi'_{\hw}\Psi'_{\hu}=\Psi'_{\hw\hu}$ whenever 
$l(\hw\hu)=l(\hw)+l(\hu)$. They satisfy the relations:
\begin{align} 
&\Psi'_{\hw} y_b\ =\ y_{\hw(b)}\Psi'_{\hw},\ \hw\in \widehat{W}, 
\label{Phixx} 
\end{align} 
where $y_{[b,j]}\equal y_b+j,$ for instance, $y_{\al_0}=1-y_\th.$ 

The formulas for nonaffine $\Psi'_i$ when $1\le i\le n$ 
are well known 
in the theory of degenerate (graded) affine Hecke algebras. 
See \cite{L,Ch4} 
and \cite{O2}, Definition 8.2. 
However, the applications to Jack polynomials and the 
Harish-Chandra theory do require  the affine intertwiners 
$\Phi_0'$ and $P'_r.$ 
\smallskip

{\bf Differential representation.} 
Introducing the variables 
$X_b=e^{x_b}$ for $b\in P ,$ we 
extend the following formulas to the derivations of $\C[X]=\C[X_b]$: 
$$ 
\partial_{a}(X_{b})\ =\ (a ,b)X_{b},\ a,b\in P . 
$$ 
Note that $w(\partial _{b})= 
\partial _{w(b)}, \  w\in W$. 

The key measure-function in the Harish-Chandra theory of spherical 
functions and its $k$\~generalizations is as follows: 
\begin{align} 
&\tau\equal\tau(x;k)= 
\prod_{\al\in R_+}|2\sinh(x_{\al}/2)|^{2\ka_\al}. 
\label{tauform} 
\end{align} 

\begin{ proposition} 
(a) The following 
{\em trigonometric Dunkl operators} 
\index{Dunkl operators (trig)} 
acting on the Laurent polynomials 
$f\in \C[X]=\C[X_b],$ 
\begin{align} 
D_{b} \equal 
&\pa_b + 
\sum_{\al\in R_+} \frac{ \ka_{\al}(b,\al^\vee)}{ 
(1-X_{\al}^{-1}) } 
\bigl( 1-s_{\al} \bigr)- (\rho_{k},b) 
\label{dunk} 
\end{align} 
are pairwise commutative, and $y_{[b,u]}=D_b+u$ 
satisfy (\ref{suka}) for 
the following action of the group $\widehat{W}$: 
$$ w^x(f)=w(f)\for w\in W,\ b^x(f)=X_b f \for b\in P .$$ 
For instance, $s_0^x(f)=X_\th s_\th(f),\ 
\pi_r^x(f)=X_r u_r^{-1}(f).$ 

(b) The operators $D_b$ are formally self-adjoint 
with respect to the inner 
product 
\begin{align} 
&\{f,g\}_\tau\equal \int f(x)g(-x)\tau dx, 
\label{tauprod} 
\end{align} 
i.e., $\tau^{-1}D_b^+\tau=D_b$ for the anti-involution $^+$ 
sending $\hw^x\mapsto (\hw^x)^{-1}$ and 
$\partial_b\mapsto -\partial_b,$ 
where $\hw\in \widehat{W}$. \sq 
\label{DUNDEG} 
\end{ proposition} 
\smallskip 

\subsection{Difference-rational case} 
We introduce the variables 
\begin{align} 
&\la_{a+b}=\la_a+\la_b, \ \la_{[b,u]}=\la_b+u,\ b\in P , 
u\in \C, 
\label{twolim} 
\end{align} 
and define 
the \index{Demazure\~Lusztig operators $S_i$} 
{\em rational Demazure\~Lusztig operators} from \cite{Ch5} 
as follows: 
\begin{align} 
&S_i\  = \   s_i^\la\ +\ 
\frac{\ka_i }{ \la_{\al_i}}(s_i^\la-1),\ 
\ 0\le i\le n, 
\label{hatS} 
\end{align} 
where by $\hw^\la$ we mean the action on $\{\la_b\}$: 
$$ 
\hw^\la(\la_b)=\la_{\hw(b)},\ b^\la(\la_c)=\la_c-(b,c) 
\for b,c\in P . 
$$ 
For instance, 
$ 
S_0  =   s_0^\la + 
\frac{\ka_0 }{ 1-\la_\th}(s_0^\la-1). 
$ 

We set 
\begin{align} 
& S_{\hw}\equal\pi_r^\la S_{i_1}\ldots S_{i_l} \for 
\hw=\pi_r s_{i_1}\ldots s_{i_l}, 
\label{Shwde} 
\\ 
&\De_b\equal S_b \for b\in P ,\ \De_i=\De_{\om_i}. 
\label{Desb} 
\end{align} 
The definition does not depend on the particular 
choice of the decomposition of $\hw\in \widehat{W}$, and 
the map $\hw\mapsto S_{\hw}$ is a homomorphism. 
The operators  $\De_b$ are pairwise 
commutative. They are called {\em difference Dunkl operators.} 
\index{Dunkl operators (difference)} 

The counterpart of $\tau$ is 
the \index{Harish-Chandra function $\si$} 
{\em asymmetric Harish\~Chandra $c\bar{c}$\~function}: 
\begin{align} 
&\si\equal \si(\la;k)=\prod_{a \in R_+} 
\frac{\Ga(\la_\al+\ka_a)\Ga(-\la_\al+\ka_\al+1) 
}{ 
\Ga(\la_\al)\Ga(-\la_\al+1) 
}, 
\label{sigma} 
\end{align} 
where  $\Ga$ is 
the classical $\Ga$\~function. 

\begin{ proposition} 
(a) The operators $S_{\hw}$ are well defined and  preserve 
the space of polynomials $\C[\la]$ in terms of $\la_b$. 
They form the degenerate double Hecke algebra, namely, 
the map $\hw\mapsto S_{\hw},\ 
y_b\mapsto \la_b$ is a faithful representation of $\HH' $. 

(b) The operators $S_{\hw}$ are formally unitary with 
respect to the inner 
product 
\begin{align} 
&\{f,g\}_\si\equal \int f(\la)g(\la)\si d\la, 
\label{siprod} 
\end{align} 
i.e., $\si^{-1}S_{\hw}^+\si=S_{\hw}^{-1}$ for 
the anti-involution $^+$ 
sending $\la_b\mapsto \la_b$ and $\hw\mapsto \hw^{-1},$ 
where $\hw\in \widehat{W}$. 
\label{DESI} 
\sq 
\end{ proposition} 

The operators 
\begin{align} 
&\La_p=p(\De_{1},\ldots,\De_{n}) \for p\in 
\C[X_1^{\pm 1},\ldots, X_n^{\pm 1}]^W 
\label{Laop} 
\end{align} 
are $W$\~invariant 
and preserve $\C[\la]^W$ (usual symmetric polynomials in $\la$). 
So if we restrict them to $W$\~invariant functions, we will 
get $W$\~invariant difference operators. However, the formulas 
for the restrictions of $\La_p$ 
are simple enough in special cases only. 

\begin{proposition} 
Given $r\in O'$, let $m_r = \sum _{w \in W/W_r} X_{w(-\om_r)}$, 
where $W_r$ is 
the stabilizer of $\om_r$ in $W$. 
Then $\la(\om_r)$ belongs to $R_+$ and 
the restriction of $\La_{m_r}(\De_b)$ 
onto $\C[\la]^W$ is as follows: 
\begin{align} 
& \La_r= \ \sum_{w\in W/W_r} 
 \prod _{\al \in \la(\om_r)} 
\frac{\la_{w(\al)}+\ka_\al }{ 
     \la_{w(\al)}   } \, w(-\om_r), 
\label{macdeop} 
\end{align} 
where $w(-\om_r)=-w(\om_r)\in P $ is considered as a 
difference operator: 
$$w(-\om_r)(\la_c)= \la_c+(w(\om_r),c). 
$$ 
\label{MACOPER} 
\sq 
\end{proposition} 

The formula for $\La_r$ is very close to the corresponding formulas 
for the general double Hecke algebra from \cite{C2}. Note that 
in the differential case, the formulas for the 
invariant operators are much more complicated than in the difference 
theory. It was (and still is) a convincing argument in favor of 
the difference theory. 

\subsection{Opdam transform} 
From now on, we assume that $\R\ni \ka_\al>0$ for all $\al$ 
and keep the notation from the previous section; 
$i$ is the imaginary 
unit; $\Re,\Im$ are the real and imaginary parts. 
\begin{theorem} 
There exists a solution $G(x,\la)$ of the eigenvalue 
problem 
\begin{align} 
& D_b (G(x,\la))=\la_b G(x,\la),\  b\in P ,\  G(0,\la)=1, 
\label{eigenop} 
\end{align} 
 holomorphic for all $\la$ 
and for $x$ in $\R^n+iU$ for a neighbourhood $U\subset \R^n$ of 
zero. If $x\in \R^n$, then 
$|G(x,\la)|\le |W|^{1/2}\exp(\max_w(w(x),\Re\la))$, so 
$G$ is bounded for $x\in \R^n$ when $\la\in i\R^n$. 
The solution of (\ref{eigenop}) is unique in the class of 
continuously differentiable functions on $\R^n$ (for a given 
$\la$). 
\label{OPDAM} 
\end{theorem} 

This theorem is from  \cite{O2} 
(Theorem 3.15 and Proposition 6.1). 
Opdam uses the fact that 
\begin{align} 
& F\equal |W|^{-1}\sum_{w\in W}G(w(x),\la) 
\label{fhyper} 
\end{align} 
 is a 
\index{hypergeometric function} 
{\em generalized hypergeometric function}, i.e., a 
$W$\~symmetric eigenfunction of the QMBP operators or 
{\em Heckman\~Opdam operators}, 
\index{Heckman\~Opdam operators $L'_p$} which are the 
restrictions $L'_p$ of the operators $p(D_{b_1},\ldots,D_{b_n})$ 
to symmetric functions: 
$$ 
L_p'F(x,\la)\ =\ p(\la_1,\ldots, \la_n)F(x,\la), 
$$ 
where $p$ is any 
$W$\~invariant polynomial of $\la_i=\la_{\om_i}$. 
The operators $L'_p$  generalize the radial parts 
of Laplace operators on the corresponding symmetric space. 
The normalization is the same: $F(0,\la)=1$. 
It fixes $F$ uniquely, 
so it is $W$\~invariant with respect to $\la$ 
as well. 

A systematic algebraic 
and analytic theory of $F$\~functions is 
due to Heck\-man and Op\-dam 
(see \cite{HO1,H,O1,HS}). There is 
a formula for $G$ in terms of $F$ (at least for generic $\la$) 
via the operators $D_b$ from (\ref{dunk}). 
The positivity of $k$ implies that 
it holds for all $\la\in \C^n$. 
See \cite{O2} for a nice and simple argument (Lemma 3.14). 
This formula and  the relation to the affine 
Knizhnik\~Zamolodchikov equation \cite{Ch4,M,C13,O2} are applied to 
establish the growth 
estimates for $G$. Actually, it  gives more than was 
formulated in the theorem 
(see Corollary 6.5, \cite{O2}). 

We introduce the \index{Opdam transform $\f$} 
{\em Opdam transform} 
(the first component of what he 
called ``Cherednik's transform") as follows: 
\begin{align} 
&\f(f)(\la)\equal \int_{\R^n}f(x)G(-x,\la)\tau dx 
\label{Optran} 
\end{align} 
for the standard measure $dx$ on $\R^n$.

\begin{ proposition} 
(a) Let us assume that the $f(x)$ are taken from the space 
$\C^\infty_c(\R^n)$ 
of $\C$\~valued  compactly supported $\infty$\~differentiable 
functions 
on $\R^n$. The inner product 
\begin{align} 
&\{f,f'\}_\tau\equal \int_{\R^n}f(x)f'(-x)\tau dx 
\label{ffprime} 
\end{align} 
satisfies the conditions of part (b), 
Proposition \ref{DUNDEG}. Namely, $\{D_b\}$ 
preserve $\C^\infty_c(\R^n)$ 
and are self-adjoint with respect to the pairing (\ref{ffprime}). 

(b) The Opdam transforms  of such functions are analytic in $\la$ 
on the whole $\C^n$ and satisfy the Paley-Wiener condition. 
A function 
$g(\la)$ is of the PW-type ($g\in PW(\C^n)$) if 
there exists a constant $A=A(g)>0$  such that, for any $N>0,$ 
\begin{align} 
& g(\la)\le C(1+|\la|)^{-N}\exp (A|\Re\la|) 
\label{PWcond} 
\end{align} 
for a proper constant $C=C(N;g)$. 
\label{OPPROP} 
\end{ proposition} 
{\em Proof.} The first claim is obvious. The Paley-Wiener condition 
follows from Theorem 8.6 \cite{O2}. The transform under 
consideration is actually 
the first component of Opdam's transform from Definition 7.9 (ibid.) 
without the complex conjugation and for the opposite 
sign of $x$ (instead of $\la$). Opdam's estimates remain valid in 
this case. 
\sq 

\subsection{Inverse transform} 
In this section we discuss the {\em inversion} (for positive $k$). 
\index{inversion ($\h'$)}  
\index{Opdam transform $\g$ (inverse)}  
The {\em inverse Opdam transform} is defined for Paley-Wiener 
functions $g(\la)$ on $\C^n$ by the formula 
\begin{align} 
&\g(g)(x)\equal \int_{i\R^n}g(\la)G(x,\la)\si d\la 
\label{invtran} 
\end{align} 
for the standard measure $d\la$. 
The transforms of such $g$ belong to $\C^\infty_c(\R^n)$. 

 The existense readily follows from 
(\ref{PWcond}) and the known properties 
of the  Ha\-rish-Chan\-dra $c$\~function (see  below). 
The embedding $\g(PW(\C^n))\subset\C^\infty_c(\R^n)$  is due to 
Opdam; it is similar to the classical 
one from \cite{He} (see also \cite{GV}). 

Let us discuss the shift of the integration contour 
in (\ref{invtran}). 
There exists an open neighborhood 
$U_+^a\subset \R^n$ of the \index{affine Weyl chamber $\bar{\CC^a}$} 
{\em closure} $\bar{\CC}^a_+\in \R^n$ of the 
affine Weyl chamber $\CC^a$ from (\ref{afchamber}) such that 
\begin{align} 
&\g(g)(x)\ =\  \int_{\xi+ i\R^n}g(\la)G(x,\la)\si d\la 
\label{contour} 
\end{align} 
for $\xi\in U_+^a$. 
Indeed,  $\ka_\al>0$ and $\si$ has no singularities in $U_+^a+i\R^n$. 
Then we use 
the classical formulas for $|\Ga(x+iy)/\Ga(x)|$ 
for real $x,y, x>0$. It gives (cf. \cite{He,O2}): 
\begin{align} 
&|\si(\la)| \ \le \ C(1+|\la|)^K, \where K=2\sum_{ \al>0} \ka_\al,\ 
\la\in U^a_+ +i\R^n, 
\label{cestim} 
\end{align} 
for sufficiently large $C>0$. Thus the products 
of PW-functions by $\si$  tend to  $0$ for $|\la|\mapsto \infty$, 
and we can switch to  $\xi$. 
Actually, we can do this  for any integrand analytic in 
$ U^a_+ +i\R^n$ 
and approaching $0$ at $\infty$. 
We come to the following. 

\begin{ proposition} 
The conditions of part (b), Proposition \ref{DESI} are satisfied 
for 
\begin{align} 
&\{g,g'\}_\si\equal\int_{i\R^n}g(\la)g'(\la)\si d\la 
\label{ggprime} 
\end{align} 
in the class of PW-functions, i.e., the operators $S_{\hw}$ are 
well defined on such functions and unitary. 
\end{ proposition} 
{\em Proof.} 
It is sufficient to check the unitarity  for the generators 
$S_i= S_{s_i}\ (0\le i\le n)$ and $\pi_r^\la\ (r\in O')$. 
 For instance, let us consider $s_0$. We will 
integrate over $\xi+i\R$, assuming that 
$\xi'\equal s_\th(\xi)+\th\in U_+^a$ 
and avoiding the wall $(\th,\xi)=1$. 

We will apply  formula (\ref{siprod}), the formula 
$$ 
\int_{\xi+i\R^n}s_0(g(\la)\si) d\la \ =\ 
\int_{\xi'+i\R^n}g(\la)\si d\la \ =\ \int_{\xi+i\R^n}g(\la)\si d\la, 
$$ 
and a similar formula for $(1-\la_\th)^{-1}g\si$, where $g$ is of 
the PW-type in $ U^a_+ +i\R^n$. 
Note that $(1-\la_\th)^{-1}\si$ is regular in this domain. 
One then has: 
\begin{align} 
&\int_{\xi+i\R^n}S_{0}(g(\la))\ g'(\la)\si d\la\notag\\ 
=\ &\int_{\xi+i\R^n}\big(s_0+\ka_0(1-\la_\th)^{-1}(s_0-1)\big) 
(g(\la))\ g'(\la) 
\si d\la 
\notag\\ 
=\ &\int_{\xi+i\R^n}g(\la)\ \big(s_0+\ka_0(s_0-1)(1-\la_\th)^{-1}\big) 
(g'(\la) 
\si) d\la 
\notag\\ 
=\ &\int_{\xi+i\R^n}g(\la)\ S_{0}^{-1}(g'(\la))\si d\la. 
\label{intint} 
\end{align} 
Since (\ref{intint}) holds for one $\xi,$ it is valid for all 
of them in $U_+^a,$ 
including $0$. 
The consideration of the other generators is  the same. 
\sq 

\begin{theorem} 
Given $\hw\in \widehat{W},\ b\in P ,\ 
f(x)\in \C^\infty_c(\R^n),\ g(\la)\in PW(\C^n),$ 
\begin{align} 
&\hw^x(G(x,\la))\ =\ S_{\hw}^{-1}G(x,\la), \hbox{\ e.g.,\ } 
 X_b G\ =\ \De_b^{-1}(G), 
\label{gidesym} 
\\ 
&\f(\hw^x(f(x)))\ =\ S_{\hw}\f(f(x)),\ 
\g(S_{\hw}(g(\la)))\ =\ \hw^x(\f(f(x))), 
\label{gisym} 
\\ 
&\f(X_b f(x))\ =\ \De_b(\f(f(x))),\ 
\g(\De_b(g(\la)))\ =\ X_b\g(g(\la)), 
\label{fxsym} 
\\ 
&\f(D_b(f(x)))\ =\ \la_b\f(f(x)),\ 
\g(\la_b g(\la))\ =\ D_b(\g(g(\la))). 
\label{fdsym} 
\end{align} 
\label{GISYM} 
\end{theorem} 
{\em Proof.} 
Applying the intertwiners from 
(\ref{Phiprime}) to $G(x,\la)$, we see that 
\begin{align} 
&(1 + \frac{\ka_i }{ \la_{\al_i} })^{-1} 
(s_i + \frac{\ka_i }{ \la_{\al_i} })(G)\ = 
s_i^\la(G), \for 1\le i\le n, 
\label{Phig} 
 \\ 
&(1 + \frac{\ka_0 }{ 1-\la_{\th} })^{-1} 
(X_\th s_\th + \frac{\ka_0 }{ 1-\la_{\th} })(G)\ = 
s_0^\la(G),\notag\\ 
&P_r'(G)\  =\  X_r u_r^{-1}(G)=(\pi_r^\la)^{-1}(G) 
 \for \ r\in O'. 
\notag \end{align} 

Due to the main property 
of the intertwiners (\ref{Phixx}), we obtain these equalities up 
to $\lambda$\~multipliers. 
The scalar factors on the left are necessary to preserve the 
normalization $G(0,\la)=1,$ so we can use the uniqueness of 
$G(x,\la)$ 
from Theorem \ref{OPDAM}. Expressing $s_j^x$ 
in terms of $s_j^\la$ (when applied to $G$!), we obtain 
(\ref{gidesym}) for $\hw=s_j$. 
It is obvious for $\hw=\pi_r$. Using the commutativity of $\hw^x$ 
and $S_{\hu}$ for $\hw,\hu\in \widehat{W}$, we establish this 
relation 
in the general case. For $w\in W$ it 
is due to Opdam. 

Formula (\ref{gisym}) results directly from 
(\ref{gidesym}) because we already  know that the $\hw^x$ 
are unitary for $\{f,f'\}_\tau$ and the $S_{\hw}$ are unitary 
with respect 
to $\{g,g'\}_\si$ for the considered classes of functions. 
See Theorem \ref{DUNDEG}, Theorem \ref{DESI}, and 
formulas (\ref{ffprime}) and (\ref{ggprime}). 

For instance, let us check (\ref{fxsym}) for $\f$, which is a 
particular case 
of (\ref{gisym}): 
\begin{align} 
&\f(X_bf(x))\ =\ \int_{\R^n} f(-x) X_b^{-1}G(x,\la)\tau dx \ 
\notag\\ 
=\ &\int_{\R^n} f(-x) \De_b(G(x,\la))\tau dx \ =\ 
\De_b(\f(f(x))). 
\label{trxb} 
\end{align} 
Thus the 
Opdam transforms of the operators of multiplication by 
$X_b$ are the 
operators $\De_b$. 

Since $D_b$ and $\la_b$ are self-adjoint for the 
corresponding inner products (see Theorem \ref{DUNDEG}), we 
obtain (\ref{fdsym}), which is in fact the defining property 
of $\f$ and $\g$. 
\sq 

\begin{corollary} 
The compositions $\g\f: C_c^\infty(\R^n)\to  C_c^\infty(\R^n)$, 
and $\f\g: PW(\C^n)\to  PW(\C^n)$ are multiplications by nonzero 
constants. The transforms  $\f,\g$ establish isomorphisms 
between the corresponding space identifying 
$\{\ ,\ \}_\tau \and \{\ ,\ \}_\si$ up to proportionality. 
\label{INVERSE} 
\end{corollary} 
{\em Proof.} 
The first statement readily follows from Theorem \ref{GISYM}. 
The composition $\g\f$ sends $\C^\infty_c(\R^n)$ into itself, 
is continuous (due to Opdam), 
and commutes with the operators $X_b$ (multiplications by $X_b$). 
Thus it is a multiplication by a function 
$u(x)$ of $\C^\infty$\~type. 
It must also commute 
with $D_b$. Hence $G(x,\la)u(x)$ is another solution 
of the eigenvalue problem (\ref{eigenop}), and we conclude that 
$u(x)$ has to be a constant. 

Let us  check that $\f\g$, which is 
a continuous operator on $PW(\C^n)$ (for any fixed $A$) 
commuting with multiplications by any $\la_b$, has to be 
a multiplication by an analytic function $v(\la)$. 
Indeed, the image of $\f\g$ with respect to the 
standard  Fourier transform ($k=0$) 
is a continuous operator on $C_c^\infty(\R^n)$ commuting with 
the derivatives $\partial/\partial x_i$. 
Thus it is a convolution 
with some function, and its inverse Fourier transform is the 
multiplication by a certain $v(\la)$. 
Since $\f\g(g)=\hbox {Const\ } g$ 
for any $g(\la)$ 
from $\f(C_c^\infty(\R^n))$,  $v$ is  constant. 

The claim about the inner products is obvious because 
 $\{\f (f) , g \}_\si$ =  $\{f ,\g(g) \}_\tau$ 
for $f\in C_c^\infty(\R^n)$, $g\in PW(\C^n)$. 
\sq 

The corollary is due to Opdam (\cite{O2},Theorem 9.13 (1)). 
His proof is different. 
He uses  Peetre's 
characterization of differential operators (similar to what 
Van Den Ban and Schlichtkrull did in \cite{BS}). 
In his approach, a certain nontrivial analytic argument 
is needed to check that 
$\g\f$ is multiplication by a constant. 
The proof above is analytically elementary. 
\smallskip 

{\bf The symmetric case.} The transforms can be readily 
reduced to the symmetric level. If $f$ is 
$W$\~invariant, then 
\begin{align} 
&\f(f(x))\ =\  \int_{\R^n} f(x) F(-x,\la)\tau dx 
\label{hartran} 
\end{align} 
for  $F$ from (\ref{fhyper}). 
Here we have applied the $W$\~sym\-met\-rization 
to the integrand in 
(\ref{Optran}) and used that $\tau$ is $W$\~invariant. 
So $\f$ {\em coincides} with  the $k$\~deformation of the 
Harish-Chandra transform 
on $W$\~invariant functions 
up to a minor technical detail. The $W$\~invariance of  $F$ in $\la$ 
results in the $W$\~invariance of $\f(f)$.

As for $\g$, we $W$\~symmetrize the integrand in the definition 
with respect to $x$ and $\la.$ Namely, we 
replace  $G$ by  $F$ and $\si$ by 
its $W$\~symmetrization. The latter 
is the genuine 
Harish-Chandra ``$c\bar{c}-$ function" 
\begin{align} 
&\si'\ =\ \prod_{a \in R_+} 
\frac{\Ga(\la_\al+\ka_a)\Ga(-\la_\al+\ka_\al) 
}{ 
\Ga(\la_\al)\Ga(-\la_\al) 
} 
\label{sigmap} 
\end{align} 
up to a coefficient of proportionality. 

Finally, given a $W$\~invariant function $f\in\C^\infty_c(\R^n),$ 
\begin{align} 
& f(x)\ =\ \hbox{Const} \int_{i\R^n}g(\la)F(x,\la)\si' d\la 
\for g=\f(f). 
\label{harish} 
\end{align} 
A similar formula holds for $\g$. 
See \cite{HC,He2} and \cite{GV}(Ch.6) for the classical theory. 

As an application, we are able to 
 calculate the Fourier transforms of the operators 
$p(X)\in \C[X_b]^W$ 
(symmetric Laurent polynomials acting by multiplication) 
in the  Harish-Chandra theory and its $k$\~deformation. 
They are exactly the 
operators $\La_p$ above. 
In the minuscule case, we obtain  formulas (\ref{macdeop}). 
We mention that in  \cite{O2} and other papers 
the pairings serving the Fourier transforms are hermitian. 
Complex conjugations can be added to ours. 

I hope that the method we used can be generalized to 
negative $k,$ to other classes of functions, 
and to the $p$\~adic theory (cf. \cite{BS,HO2}). 
The relations (\ref{gidesym}) considered as 
difference equations for $G(x,\la)$ 
with respect to $\la$ may also help with the growth estimates via 
the theory of difference equations and the equivalence with 
the difference Knizhnik\~Zamolodchikov equations. 
However, I believe 
that the main advantage here will be connected with 
the difference Fourier transform.

\bibliographystyle{unsrt}

\end{document}